\documentclass[11pt]{article}

\newcommand{\nex}{x}

\newcommand{\de}{\,\mathrm{d}}                               
\newcommand{\e}{\operatorname{e}}                               
\newcommand{\im}{\operatorname{i}}

\newcommand{\p}{\partial}

\newcommand{\real}{\mathrm{Re}\,}

\newcommand{\R}{\mathbb{R}}       
\newcommand{\C}{\mathbb{C}}       
\newcommand{\N}{\mathbb{N}}

\newcommand{\mgf}{ H}                                        
\newcommand{\elf}{ E}  
  
\newcommand{\dive}{\operatorname{div}}

\newcommand{\nor}{{\nu}} 
\newcommand{\curl}{\operatorname{curl}}

\newcommand{\oP}{\operatorname{P}}  
\newcommand{\Pper}{{\oP_{\nu}}}
\newcommand{\Ppar}{{\oP_{t}}}
\newcommand{\ED}{\R^3\setminus\overline\Omega}
\newcommand{\oED}{\R^3\setminus\Omega}
\newcommand{\der}{{\rm D}}

\usepackage{xcolor}
\definecolor{capa}{RGB}{50, 90, 160}
\definecolor{delftblue}{RGB}{0,61,165}
\usepackage[colorlinks=true,
  linkcolor=delftblue,
  urlcolor=delftblue,
  citecolor=delftblue]{hyperref}
\newcommand{\Id}{{\rm I}}  
\newcommand{\oA}{{\rm A}}  
\newcommand{\oB}{{\rm B}}  
\newcommand{\oH}{{\rm H}}  
\newcommand{\oR}{{\rm R}}  
\newcommand{\oS}{{\rm S}}  
\newcommand{\oK}{{\rm K}}  
\newcommand{\oT}{{\rm T}}  
\newcommand{\oQ}{{\rm Q}}  
\newcommand{\oL}{{\rm L}}  
\newcommand{\oM}{{\rm M}}  
\usepackage{enumitem}

\usepackage{multirow}
\usepackage{amsmath}
\usepackage{amsfonts}
\usepackage{amsmath}
\usepackage{amssymb}  
\usepackage{graphicx}
\usepackage{caption}
\usepackage{mathrsfs}
\usepackage{upgreek}
\usepackage{amsthm}
\usepackage{subfig}
\usepackage{booktabs}
\usepackage{authblk}
\usepackage{cite}
\usepackage{bbm}
\usepackage{url}

\usepackage[]{algorithm2e}

\newtheorem{theorem}{Theorem}[section]
\newtheorem{lemma}[theorem]{Lemma}
\newtheorem{proposition}[theorem]{Proposition}
\newtheorem{corollary}[theorem]{Corollary}
\newtheorem{remark}[theorem]{Remark}

\title{Maxwell \`a la Helmholtz: 
Electromagnetic scattering by 3D perfect electric conductors via Helmholtz integral operators}

\author[1]{Juan Burbano-Gallegos\thanks{\href{mailto:j.s.burbano@utwente.nl}{j.s.burbano@utwente.nl}}}
\author[1]{Carlos P\'erez-Arancibia\thanks{\href{mailto:c.a.perezarancibia@utwente.nl}{c.a.perezarancibia@utwente.nl}}}
\author[2]{Catalin Turc\thanks{\href{mailto:catalin.turc@njit.edu}{catalin.turc@njit.edu}}}
\affil[1]{\small{Department of Applied Mathematics, University of Twente, Netherlands}}
\affil[2]{\small{Department of Mathematical Sciences, New Jersey Institute of Technology, USA}}
\date{\today}

\begin{document}
\maketitle

\begin{abstract} This paper introduces a novel class of indirect boundary integral equation (BIE) formulations for the solution of electromagnetic scattering problems involving smooth perfectly electric conductors (PECs) in three-dimensions. These combined-field-type BIE formulations rely exclusively on classical Helmholtz boundary operators, resulting in provably well-posed, frequency-robust,  Fredholm second-kind BIEs. Notably, we prove that the proposed formulations are free from spurious resonances, while retaining the versatility of Helmholtz integral operators. The approach is based on the equivalence between the Maxwell PEC scattering problem and two independent vector Helmholtz boundary value problems for the electric and magnetic fields, with boundary conditions defined in terms of the Dirichlet and Neumann traces of the corresponding vector Helmholtz solutions. While certain aspects of this equivalence (for the electric field) have been previously exploited in the so-called field-only BIE formulations, we here rigorously establish and generalize the equivalence between Maxwell and Helmholtz problems for both fields. Finally, a variety of numerical examples highlights the robustness and accuracy of the proposed approach when combined with Density Interpolation-based Nyström methods and fast linear algebra solvers, implemented in the open-source Julia package \texttt{Inti.jl}.
 \end{abstract}

\section{Introduction}

The numerical solution of frequency-domain Maxwell scattering problems via boundary integral equations (BIEs) is widely acknowledged to be significantly more challenging than its arguably simpler scalar precursor in acoustics, governed by the Helmholtz equation. This gap in difficulty is perhaps already evident in the development of the most fundamental component of BIE formulations---the field representation formula. Notably, the integral representation for electromagnetics—the Stratton–Chu formula~\cite{Stratton:1939ha}—appeared nearly half a century after the analogous Kirchhoff formula for acoustics~\cite{kirchhoff1882}. Yet, this gap is sometimes underappreciated  in the BIE community, where advances in singular integration techniques, fast algorithms, preconditioning, and other areas are often made for the Helmholtz problem, with the tacit assumption that they can be somehow extended to Maxwell’s equations. In practice, however, transferring these advances to full 3D electromagnetic problems is often hindered by challenges intrinsic to traditional Maxwell BIE formulations, such as the electric field integral equation (EFIE), magnetic field integral equation (MFIE), and combined field integral equation (CFIE)~\cite{volakis2012integral}. Amongst these challenges we mention the need to work with surface-tangent vector fields (e.g., surface currents) with more subtle regularity properties, the presence of complicated kernels (particularly in the case of the EFIE operator), difficulties associated with Calderón preconditioning, and the occurrence of low-frequency breakdown. The latter is a phenomenon whereby the conditioning of the discretized BIE system matrix deteriorates dramatically as the frequency approaches zero, leading to unstable or inaccurate numerical solutions. Nevertheless, there have been important contributions in the literature in which these difficulties have been successfully addressed, including the high-order accurate Nyström methods developed in~\cite{Bruno2009electromagnetic,ganesh2006spectrally,garza2020boundary,hu2021chebyshev}, as well as the extensive literature on boundary element methods and the method of moments---most of which rely on low-order approximations, both in terms of accuracy and geometric representation---such as~\cite{rao1982electromagnetic,andriulli2008multiplicative,seo2005fast,bleszynski1996aim,song1995multilevel,song1997multilevel,adrian2019refinement}, to mention just a few. Also, the low frequency breakdown has been addressed in the novel framework of current and charge BIE formulations~\cite{epstein2010debye,epstein2013debye,taskinen2006current,bendali2012extension} which incorporate Maxwell alongside Helmholtz boundary integral operators.

This paper aims to streamline the numerical solution of full 3D electromagnetic scattering problems involving smooth perfectly electrically conducting (PEC) obstacles by bridging the gap in difficulty between the original Maxwell PEC scattering problem and more tractable vector Helmholtz formulations. Specifically, it introduces frequency-robust BIE formulations for PEC scattering problems that rely solely on the classical Helmholtz boundary integral operators from the Calderón calculus, applied componentwise to vector surface densities that are not constrained to be tangential to the obstacle's surface. Remarkably, the proposed approach allows any existing Helmholtz BIE solver to be seamlessly adapted for solving Maxwell PEC scattering problems, provided that a sufficiently accurate representation of the surface is available.

We develop two types of indirect BIE formulations---one for the electric field and one for the magnetic field---based on the classical combined-field potential from acoustics~\cite{panich1965question,leis1965dirichletschen,brakhage1965dirichletsche,Burton1971Application}, and  prove their well-posedness in Hölder spaces via the Fredholm alternative~\cite{COLTON:1983}. Notably, the proposed approach entirely avoids the use of differential operators acting on either the layer potentials or the surface densities---a common and often burdensome feature of traditional electromagnetic BIEs. By building solely on Helmholtz operators, the proposed BIE formulations are naturally compatible and easily interfaced with the extensive numerical techniques developed for scalar Helmholtz problems---including advanced quadrature schemes, fast solvers, hierarchical matrix compression, and preconditioning strategies---thereby enabling efficient and scalable solutions to full 3D PEC scattering problems.

The key in the derivation of our novel indirect BIE formulations is found in two theorems that establish the equivalence between the PEC scattering problem---posed in terms of Maxwell’s equations with the standard boundary condition on the tangential component of the total electric field---and vector Helmholtz problems for the electric and magnetic fields (Theorems~\ref{thm:elf} and~\ref{lem:equiv_int}, respectively)---with boundary conditions involving both Dirichlet and Nuemann traces of each electric/magnetic field component. Specifically, one of these conditions enforces the PEC boundary condition in standard form (i.e., the condition on the tangential/normal component of the total electric/magnetic field~\cite[Sec.~1.4]{volakis2012integral}), while the other, involving the Neumann trace of the field together with either the mean curvature or the curvature operator depending on the formulation, ensures the divergence-free character of the vector Helmholtz equation solutions. While the latter condition for the electric field has appeared in several previous works~\cite{desanto1993new,desanto2006new,yuffa20183}, we provide here a simple derivation based on elementary vector calculus (Lemma~\ref{lem:neumann_trace}), following a classical tubular neighborhood differential geometry approach which can be found in electromagnetic textbooks~(e.g.~\cite[Sec.~2.5.6]{NEDELEC:2001}). Furthermore, we exploit this non-standard condition to establish the equivalence between the vector Helmholtz boundary value problem for the electric field  and the original Maxwell PEC scattering problem. Notably, the corresponding formulation for the magnetic field features a novel boundary condition that, to the best of our knowledge, has not appeared previously in the literature. We also extend these results to the problematic zero-frequency limit, where the electromagnetic problem decouples into two independent vector Laplace boundary value problems for the electric and magnetic fields (see Theorems~\ref{rm:zero_frequency} and~\ref{rm:mgf_zero_frequency}). 

A class of methods, known as \emph{field-only surface integral equations}~\cite{klaseboer2016nonsingular,sun2017robust,sun2020field,sun2020field_die}, has previously explored reformulations of 3D electromagnetic problems in terms of Helmholtz equations, bearing similarities to our approach. In particular, the work~\cite{sun2020field} employs the aforementioned non-standard boundary condition derived in~\cite{yuffa20183}, involving the normal component of the normal derivative of the electric field, to construct a direct BIE formulation. However, the resulting system is not frequency robust, as it suffers from spurious resonances. Furthermore, the formulation is not of the second-kind Fredholm type, and neither the conditioning of the resulting BIE discretization nor any preconditioning strategy is properly discussed, so that good conditioning cannot be guaranteed in general. Although the authors state that their method avoids low-frequency (long-wavelength) breakdown (see~\cite{andriulli2008multiplicative} for a detailed discussion of the issue in the context of classical BIE formulations for Maxwell's equations), no investigation is provided regarding the conditioning of the linear system in that regime.
In fact, as shown in~\cite{werner1963perfect_reflection}, in the zero-frequency limit ($k = 0$), the electrostatic problem with the PEC boundary conditions lacks uniqueness unless the the values of the surface charge integrals~{\cite[Eq.~1.32]{werner1963perfect_reflection}} are prescribed. Therefore, similar issues are expected to arise in any direct BIE formulation of that limit--yet none of these considerations are addressed in~\cite{sun2020field}. The problem of spurious resonances was partially addressed in~\cite{maltez2024combined}, where a direct combined-field formulation was proposed and fully analyzed for the case of PEC spheres. However, neither~\cite{sun2020field} nor~\cite{maltez2024combined} discusses a corresponding formulation for the magnetic field based on the boundary condition involving the normal derivative of the magnetic field.
In contrast, our approach is provably frequency robust and, after suitable Calderón-type preconditioning, yields second-kind Fredholm formulations for both the electric and magnetic vector Helmholtz problems. Moreover, we address the manifestation of the low-frequency issues in our electric-field formulation---while our magnetic-field formulation remains unaffected by low-frequency breakdown provided that the scatterer’s surface is simply connected. Indeed, as the frequency tends to zero, the standard Helmholtz combined-field ansatz gives rise to electric-field BIEs that become contaminated by a nearly non-trivial kernel in the system matrix, associated with harmonic fields satisfying homogeneous Dirichlet boundary conditions on the surface~\cite{werner1963perfect_reflection}. {To address this, we introduce and analyze a modified combined-field potential that effectively restores accuracy arbitrarily close to the zero-frequency limit}. 
Despite the important differences between our methodology and the previous field-only approaches, we adopt the terminology originally introduced in~\cite{klaseboer2016nonsingular}, referring to our formulations as the \emph{Electric Field-Only Integral Equation} (EFOIE) and \emph{Magnetic Field-Only Integral Equation} (MFOIE), respectively.

From a numerical viewpoint, our approach allows us to take full advantage of technology developed by some of the authors over the past decade for the solution of Helmholtz BIEs. In particular, leveraging high-order Nyström methods based on the Density Interpolation technique~\cite{HDI3D,perez2019planewave,perez2020planewave,faria2021general}, we present a variety of numerical examples to validate the proposed BIE formulations. All numerical experiments have been carried out using the \texttt{Inti.jl} software package~\cite{IntegralEquations_Inti_2025}, which provides direct interfaces with both H-matrix and Fast Multipole methods, and includes full implementations of the General-Purpose Density Interpolation Methods~\cite{faria2021general}. 

This paper is organized as follows: After setting up the problem in Section~\ref{sec:problem_setup}, we present in Section~\ref{sec:equiv_helm} two main theorems establishing the equivalence between the electromagnetic PEC scattering problem and certain vector Helmholtz boundary value problems for the electric and magnetic fields. Next, Section~\ref{sec:equiv_lap} explores the particular cases of electrostatics and magnetostatics, where we present additional theorems connecting the classical electromagnetic/magnetostatic boundary value problems to equivalent vector Laplace problems for the electric and magnetic fields.
In Section~\ref{sec:BIE_formulations}, we develop indirect combined-field-only boundary integral equations for the vector Helmholtz boundary value problems, formulated solely in terms of Helmholtz and Laplace boundary integral operators. We then prove the well-posedness of these equations in classical Hölder spaces and establish the Fredholm second-kind character of their Calderón-regularized versions.
Finally, Section~\ref{sec:numerics} presents a variety of numerical examples, demonstrating the performance of our methods. These examples are implemented using a density-interpolation-based Nyström discretization of Helmholtz integral operators, available in the software package Inti.jl.

\section{Scattering problem\label{sec:problem_setup}}

{In this paper, we consider a bounded PEC obstacle embedded in an unbounded, homogeneous, isotropic medium characterized by dielectric and magnetic constants $\epsilon > 0$ and $\mu>0$, respectively. We let $\Gamma$ denote the boundary of the PEC obstacle, which occupies an open, bounded region $\Omega \subset \mathbb{R}^3$, with the exterior domain $\mathbb{R}^3 \setminus \overline{\Omega}$ assumed to be connected. The surface $\Gamma = \p \Omega$ is assumed to consist of $J$ disjoint, closed, and bounded components $\Gamma_j$, for $j \in {1, \ldots, J}$, each of which is at least $C^{2,\alpha}$-smooth for some $\alpha \in (0,1]$.  This regularity ensures that geometric quantities such as the curvature and the unit normal vector are well-defined and Hölder continuous on $\Gamma$ (see Appendix~\ref{app:curv_tens}). Unless stated otherwise, $\Gamma$ is not assumed to be simply connected. }

We consider an incident electromagnetic field $E^i, H^i \in C^{1,\alpha}(U, \mathbb{C}^3)$\footnote{Let $X \subset \mathbb{R}^m$ and $Y \subset \mathbb{C}^n$. Given $k \in \mathbb{N}_0$ and $\alpha \in (0,1]$, the Hölder space $C^{k,\alpha}(X, Y)$ is defined as 
$$
C^{k,\alpha}(X, Y) = \left\{ F=(f_1,\cdots,f_n): X \to Y \ \middle| \ f_i \in C^k(X), \ \forall i\in\{1,\ldots,n\}, \ \text{and} \ \forall |\beta| = k, \ \der^\beta u_i \in C^{0,\alpha}(X) \right\},
$$
with the norm
\[
\|F\|_{C^{k,\alpha}(X,Y)} := \sum_{i=1}^n \left( \sum_{|\beta| \le k} \|\der^\beta f_i\|_{C(X)} + \sum_{|\beta| = k} \left[\der^\beta f_i\right]_{C^{0,\alpha}(X)} \right).
\]
For the definition of the standard Hölder space $C^{0,\alpha}(X)$, as well as the norm $\|\cdot\|_{C(X)}$ and the seminorm $[\ \cdot\ ]_{C^{0,\alpha}(X)}$, we refer the reader to~\cite[Sec. 4.1]{gilbarg2001elliptic}.}, defined in an open set $U \subset \R^3$ containing~$\Gamma$, which satisfies the time-harmonic Maxwell equations:
\begin{equation}\label{eq:maxwell_inc}
\curl E^i - \im\omega \mu H^i =0 \quad \text{and} \quad \curl H^i + \im\omega \epsilon E^i = 0 \quad\text{in}\quad U.
\end{equation}
Note that standard incident fields---such as plane waves and dipole sources located away from the obstacle surface---satisfy these conditions.

When this field interacts with the PEC obstacle, it generates scattered electric and magnetic fields, denoted by $\elf^s$ and $\mgf^s$, respectively. These scattered fields satisfy the time-harmonic Maxwell equations:
\begin{subequations}\begin{equation}\label{eq:maxwell}
\curl \elf^s - \im\omega \mu\mgf^s = 0 \quad \text{and} \quad \curl \mgf^s + \im\omega \epsilon \elf^s = 0 \quad \text{in} \quad \ED,
\end{equation}
with the electric field subject to the PEC boundary condition~\cite{kirsch2016mathematical}: 
\begin{equation}\label{eq:PEC_BC}
\nor\times\elf^s =  -\nor\times \elf^i \quad\text{on}\quad\Gamma,
\end{equation}
where $\nu \in C^{1,\alpha}(\Gamma, \R^3)$ is the unit normal vector to $\Gamma$, pointing outward from $\Omega$. The scattered field is further required to satisfy the Silver--Müller radiation condition at infinity, given by   either
\begin{equation}\label{eq:rad_cond_sm}\begin{split}
\lim _{|\nex| \rightarrow \infty} |\nex| \left( \curl\elf^s(\nex) \times \frac{\nex}{|\nex|} - \im k\elf^s(\nex) \right) = 0&\quad\text{or}\\
\lim _{|\nex| \rightarrow \infty} |\nex| \left( \curl H^s(\nex) \times \frac{\nex}{|\nex|} - \im kH^s(\nex) \right) = 0&,
\end{split}\end{equation}\label{eq:scatt_problem}\end{subequations}
where the limit holds uniformly with respect to all directions $\nex/|\nex|$ and where $k:=\omega\sqrt{\epsilon\mu}>0$ is the wavenumber.  The time dependence $\e^{-\im \omega t}$, with $\omega > 0$ denoting the angular frequency, is assumed throughout the paper.
\begin{figure}[htbp]
  \centering
  \includegraphics[scale=1]{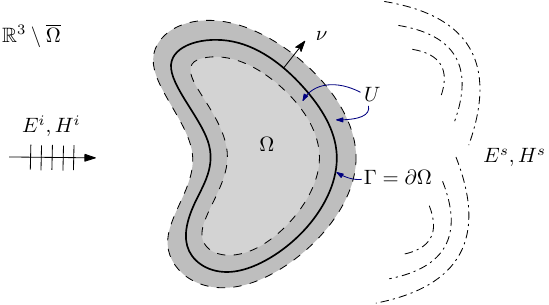}
  \caption{Illustration of the domains and fields involved in the PEC scattering problem setup.}
  \label{fig:your_figure_label}
\end{figure}

As is well known, under the aforementioned smoothness assumption of the incident field, the problem of scattering~\eqref{eq:scatt_problem}  admits a unique solution $E^s,H^s\in C^1(\ED) \cap C(\R^3 \setminus \Omega)$~\cite[Thm, 3.40]{kirsch2016mathematical}. Furthermore, the scattered field can be represented by the Stratton--Chu formula, as given by~\cite[Thm. 3.30]{kirsch2016mathematical}:
\begin{subequations}\begin{align}
E^s(x) =& \curl \int_{\Gamma}[\nor(y) \times E^s(y)] G(x, y) \de s(y) -\frac{1}{\im\omega \epsilon} \curl\curl \int_{\Gamma}[\nor(y) \times H^s(y)] G(x, y) \de s(y),\label{eq:SC_int_elf}\\
H^s(x)=& \curl \int_{\Gamma}[\nor(y) \times H^s(y)] G(x, y) \de s(y) +\frac{1}{\im\omega \mu} \curl\curl \int_{\Gamma}[\nor(y) \times E^s(y)] G(x, y) \de s(y),\label{eq:SC_int_mgf}
\end{align}\label{eq:stratton-chu}\end{subequations}
for all $x\in\ED$, where $$
G(x,y) := \frac{\e^{\im k|x-y|}}{4\pi |x-y|},\quad x\neq y,
$$ is the free-space Green's function for the Helmholtz equation.

While the classical well-posedness of problem~\eqref{eq:scatt_problem} is established in $C^{1}(\ED,\C^3) \cap C(\R^3 \setminus \Omega,\C^3)$, the assumed $C^{2,\alpha}$-smoothness of the surface $\Gamma$, together with the $C^{1,\alpha}(U,\C^3)$ regularity of the incident field, can be leveraged to establish higher $C^{2}(\ED,\C^3)\cap C^{1,\alpha}(\R^3\setminus\Omega,\C^3)$-regularity of the scattered electromagnetic field, $E^s$ and $H^s$, as we prove below in Corollary~\ref{eq:reg_solution}.

\section{Maxwell’s equations as a vector Helmholtz boundary value problem\label{sec:equiv_helm}}

In this section, we establish the connection between the Maxwell PEC scattering problem~\eqref{eq:scatt_problem} and associated vector Helmholtz boundary value problems, where the boundary conditions are expressed through the Dirichlet and Neumann traces of the electromagnetic field. These reformulations will serve as the foundation for the BIE formulations developed in Section~\ref{sec:BIE_formulations}, where Helmholtz integral operators are employed to represent the electromagnetic field.

We begin by introducing the notation for the exterior traces of sufficiently regular $\mathbb{C}^n$-valued vector fields, with $n \in \{1, 3\}$. 
Let $U_e \subset \mathbb{R}^3 \setminus \overline{\Omega}$ be an open set such that $\Gamma \subset \overline{U_e}$. We recall here that $\Gamma = \partial \Omega$ is assumed to be a closed $C^{2,\alpha}$ smooth regular surface with unit normal vector $\nu \in C^{1,\alpha}(\Gamma, \mathbb{R}^3)$.
The exterior Dirichlet trace is then defined as
\begin{subequations}\begin{equation}\label{eq:Dir_trace_def}
\gamma: C^{m,\beta}(\overline{U_e},\C^n)\to C^{m,\beta}(\Gamma,\C^n),\quad \gamma F(\nex) = \lim_{\delta\to 0+} F(\nex+\delta\nor(\nex)),\quad \nex\in\Gamma,
\end{equation}
where $(m,\beta)\in(\{0,1\}\times [0,1])\cup(\{2\}\times [0,\alpha])$. 

Similarly, we introduce the following notation for the exterior Neumann trace: 
\begin{equation}\label{eq:Neu_trace_def}
\p_\nu: C^{m+1,\beta}(\overline{U_e},\C^n)\to C^{m,\beta}(\Gamma,\C^n),\quad\p_\nu F(\nex) := \lim_{\delta\to 0+}\der F(x+\delta\nu(x)) \nor(\nex),\quad\nex\in\Gamma,
\end{equation}
for $(m,\beta)\in(\{0,1\}\times [0,1])\cup(\{2\}\times [0,\alpha])$,
\label{eq:ext_traces}\end{subequations} where 
$\der: F\mapsto(\der F)_{i,j}= \frac{\p F_{i}}{\p x_j}$,  $i\in\{1,\ldots,n\}$, $ j\in\{1,2,3\},$
 denotes the total derivative (which is represented by the Jacobian matrix $(\p F_i/\p x_j)_{i,j}$). Both trace operators $\gamma$ and $\p_\nu$ are well-defined and bounded by virtue of the assumed $C^{2,\alpha}$-smoothness, $\alpha\in(0,1)$, of the surface $\Gamma$.  In some cases, we also use the notation $\frac{\p F}{\p \nu}$ to denote the Neumann trace~\eqref{eq:Neu_trace_def} of a vector field $F$. 

\begin{remark}
Note that, when $n = 3$, the above definition~\eqref{eq:Neu_trace_def} implies that the Neumann trace~$\p_\nu F$ is a vector field whose coordinate components correspond to the scalar Neumann trace of each component $(F)_i=F_i$ ($i \in {1, 2, 3}$) of the vector field $F$, i.e., $(\p_\nu F)_i=\p_\nu F_i$.
\end{remark} 

In what follows, it will be useful to introduce the notation for the projectors on the tangent plane of $\Gamma$ and along its normal. Specifically, we define the operators:
\begin{subequations}\begin{align}
\oP_\nu:C^{m,\beta}(\Gamma,\C^3)\to C^{m,\beta}(\Gamma,\C^3),&\quad  \Pper\varphi(x)=(\nu(x)\cdot\varphi(x))\nu(x)\quad\text{and}\label{eq:P_perp_op}\\
 \oP_t:C^{m,\beta}(\Gamma,\C^3)\to C^{m,\beta}(\Gamma,\C^3),&\quad  \Ppar\varphi(x)=\varphi(x)-(\nu(x)\cdot\varphi(x))\nu(x),\quad x\in\Gamma,\label{eq:P_par_op}
\end{align}\label{eq:projectors}\end{subequations}
for $(m,\beta)\in(\{0\}\times [0,1])\cup(\{1\}\times [0,\alpha])$,  which map the vector $\varphi(x)$ onto the normal and the tangent plane at $x\in\Gamma$, respectively. Clearly, since $\nu\in C^{1,\alpha}(\Gamma,\R^3)$, both projectors are bounded.  \\

Our first lemma establishes an expression for the normal derivative of a sufficiently smooth vector field in terms of its curl and divergence on the surface $\Gamma$. 
\begin{lemma}\label{lem:neumann_trace} Let $U_e \subset \mathbb{R}^3 \setminus \overline{\Omega}$ be an open set such that $\Gamma \subset \overline{U_e}$, where $\Gamma = \partial \Omega$ is a closed, $C^{2,\alpha}$-smooth regular surface with unit normal vector $\nu \in C^{1,\alpha}(\Gamma, \mathbb{R}^3)$. Suppose $F \in C^{1,\beta}(\overline{U_e}, \mathbb{C}^3)$ with $\beta \in [0, \alpha]$. Then, the following identity holds on $\Gamma$:  
\begin{equation}\label{eq:field_on_boundary}
\p_\nu F = \left\{\nabla_{\Gamma}(\nu \cdot \gamma F) - \mathscr{R}(\gamma F) - \nu \times \gamma(\curl F)\right\} 
+ \left\{\gamma(\dive F) - \dive_{\Gamma} (\Ppar \gamma F) - 2 \mathscr{H}(\nu \cdot \gamma F)\right\} \nu,
\end{equation}  
where $\p_\nu F \in C^{0,\beta}(\Gamma,\C^3)$, with $\mathscr{R} \in C^{0,\alpha}(\Gamma, \mathbb{R}^{3 \times 3})$, and $\mathscr{H} \in C^{0,\alpha}(\Gamma, \mathbb{R})$ denoting the surface curvature operator and mean curvature of $\Gamma$, respectively.
\end{lemma}
\begin{proof} 
We begin by ``lifting" the normal vector field $\nor$ off of the surface $\Gamma$ following the presentation in~\cite[Sec. 2.5.6]{NEDELEC:2001}. We do so by considering a tubular neighborhood $U_\delta(\Gamma)=\{x+s \nu(x)\in\R^3: x \in \Gamma,|s|<\delta\}=\bigcup_{|s|<\delta} \Gamma_s \subset\R^3$ of~$\Gamma$ which comprises the parallel surfaces $\Gamma_s := \{ x + s\nu(x): x \in \Gamma\}$ for $|s|<\delta$ and $\delta > 0$ sufficiently small. The lifted normal vector field $\nu_s: U_\delta(\Gamma)\to \mathbb{R}^3$ is defined by means of the formula $\nu_s(y) = \nu(\nex)$, where $y = \nex + s\nor(\nex) \in U_\delta(\Gamma)$, $\nex \in \Gamma$. Clearly, $\nu_s\in C^{1,\alpha}(U_\delta(\Gamma),\R^3)$ is the unit normal vector field on $\Gamma_s$ for $s\in(-\delta,\delta)$. The curvature operator of $\Gamma_s$, which is denoted by $\mathscr R_s\in C^{0,\alpha}(U_\delta(\Gamma),\R^{3\times 3})$, is then given by $\mathscr R_s=\der \nu_s$. We recall that $\mathscr R_s$ is symmetric and  acts on the tangent plane, meaning that $\mathscr R_s\nu_s= 0\in\C^3$  and $\mathscr R_s v\cdot\nu_s=0$ for all $v\in\C^3$ (see Appendix~\ref{app:curv_tens}). 

With these observations at hand, we first consider the identity
$$
\nabla_{\Gamma}(\nor\cdot\gamma F) =\lim_{s\to 0+}\left\{\nabla (\nu_s\cdot F)-\nu_s\p_s(\nu_s\cdot F)\right\},
$$
which follows from the fact that the surface gradient on $\Gamma_s$, denoted by $\nabla_{\Gamma_s}$, can be expressed as $\nabla_{\Gamma_s}=\nabla - \nu_s \p_s$. Expanding the first term inside the limit we get  
$$
\nabla(\nu_s\cdot F) =  (\der \nu_s)^\top F+(\der F)^\top\nu_s=\mathscr R_s F+(\der F)^\top\nu_s\quad\text{in}\quad  U_\delta(\Gamma)\cap U_e,
$$
where we used the above mentioned symmetry of curvature operator $\mathscr R_s=\der\nu_s=(\der\nu_s)^\top$.
The transpose term above, on the other hand, can be expressed in terms of $\curl F$. Indeed,  from the identity
$$
\nu_s\times\curl F= (\der F)^\top \nu_s-(\der F)\nu_s= (\der F)^\top \nu_s-\p_s F,
$$
we arrive at
\begin{equation}\label{eq:id_1}
\nabla(\nu_s\cdot F)= \mathscr R_s F+\nu_s\times\curl F+\p_sF\quad\text{in}\quad  U_\delta(\Gamma)\cap U_e.\end{equation}

Using again the product rule, we obtain 
$$
\p_s(\nu_s\cdot  F) =\p_s\nu_s\cdot F+ \nu_s\cdot\p_s F=\nu_s\cdot\p_s F,
$$
since $\p_s\nu_s=0$ due to the fact that $\nu_s$ remains constant along the $s$-direction. Therefore, taking the limit as $s\to0+$, we get
\begin{align}\label{eq:id_5}
\nabla_{\Gamma}(\nor\cdot\gamma F) &~=\lim_{s\to 0+}\left\{\mathscr R_s F+\nu_s\times\curl F+\p_sF-(\nu_s\cdot\p_s F)\nu_s\right\}\nonumber\\
 &~=\mathscr R(\gamma F)+\nor\times\gamma(\curl F)+\p_\nu F-(\nor\cdot\p_\nu F)\nor,
\end{align}
where the continuity of $F$, $\der F$, and $\mathscr R_s$ up to and including $\Gamma$ were used. In particular, we employed  the limit $\mathscr R = \lim_{s\to 0+}\mathscr R_s$.

Similarly, taking divergence to $F$ using the decomposition $F=\nu_s\times F\times\nu_s+(\nu_s\cdot F)\nu_s$, we obtain
\begin{equation}\label{eq:id_3}
\dive F=\dive_{\Gamma_s}(\nu_s\times F\times\nu_s) + \dive((\nu_s\cdot F)\nu_s)\quad\text{in}\quad U_\delta(\Gamma)\cap U_e,
\end{equation}
where $\dive_{\Gamma_s}$ denote the surface divergence operator acting on vector fields tangent to $\Gamma_s$ (see, e.g.,~\cite[Sec. 2.5.6]{NEDELEC:2001}). 
Applying, again, standard vector calculus identities, we get 
\begin{align}
\dive((\nu_s\cdot F)\nu_s) =&~ \nabla(\nu_s\cdot F)\cdot\nu_s+(\nu_s\cdot  F)\dive\nu_s\nonumber\\
 =&~ \p_sF\cdot\nu_s+2(\nu_s\cdot F)\mathscr H_s\label{eq:id_2}
\end{align}
in $U_\delta(\Gamma)\cap U_e$, where~\eqref{eq:id_1} and the properties of the curvature operator we used; in particular, the fact that $\mathscr R_sF\cdot\nu_s=0$ and $\dive\nu_s=2\mathscr H_s$. 

Therefore, replacing~\eqref{eq:id_2} in~\eqref{eq:id_3} and taking the limits as $s\to 0+$ we arrive at
\begin{equation}\label{eq:id_4}
\gamma(\dive F)=\dive_{\Gamma}(\Ppar\gamma F) + \p_\nu F\cdot\nu+2(\nu\cdot\gamma F)\mathscr H,
\end{equation}
where the limits $\Ppar\gamma F = \lim_{s\to 0+}\nu_s\times F\times\nu_s$ and $\mathscr H= \lim_{s\to 0+}\mathscr H_s$ were used.

Combining~\eqref{eq:id_5} and~\eqref{eq:id_4} to form $\p_\nu F$, the desired identity~\eqref{eq:field_on_boundary} follows. This completes the proof.
\end{proof}

A result analogous to Lemma~\ref{lem:neumann_trace} can be formulated for the interior Neumann trace of a sufficiently regular vector field. Such result will play a crucial role in Section~\ref{sec:electric_field_BIE}, where it will be used to establish the uniqueness of a boundary integral equation.

\begin{corollary}\label{cor:neumann_trace_interior}
Let
\begin{subequations}\begin{equation}
\gamma^-: C^{m,\beta}(\overline{U_i},\C^n)\to C^{m,\beta}(\Gamma,\C^n),\quad \gamma^- F(\nex) = \lim_{\delta\to 0+} F(\nex-\delta\nor(\nex)),\quad \nex\in\Gamma,
\end{equation}
and
\begin{equation}
\p^-_\nu: C^{m+1,\beta}(\overline{U_i},\C^n)\to C^{m,\beta}(\Gamma,\C^n),\quad\p^-_\nu F(\nex) := \lim_{\delta\to 0+}\der F(x-\delta\nu(x)) \nor(\nex),\quad\nex\in\Gamma,
\end{equation}\label{eq:int_traces}\end{subequations}
for $(m,\beta)\in(\{0,1\}\times [0,1])\cup(\{2\}\times [0,\alpha])$, where $U_i\subset\Omega$ is an open set whose closure contains~$\Gamma$, denote the interior Dirichlet and Neumann trace operators on the $C^{2,\alpha}$-smooth surface $\Gamma=\p\Omega$, respectively. Then, for a field $F\in C^{1,\alpha}(\overline{U_i},\C^3)$ it holds  that  
\begin{equation}\label{eq:field_on_boundary_int}
\begin{aligned}
\partial_\nu^- F =& \Big\{ \nabla_{\Gamma} (\nu \cdot \gamma^- F) - \mathscr{R}(\gamma^- F) - \nu \times \gamma^-(\curl F) \Big\} +  \Big\{ \gamma^-(\dive F) - \dive_{\Gamma} (\Ppar \gamma^- F) - 2 \mathscr{H}(\nu \cdot \gamma^- F) \Big\} \nu,
\end{aligned}
\end{equation}  
with $\p_\nu^-F\in C^{0,\alpha}(\Gamma,\C^3).$
\end{corollary}
\begin{proof} The proof is analogous to Lemma~\ref{lem:neumann_trace}, so it is omitted for brevity.
\end{proof}

The following theorem establishes an equivalence between the electromagnetic scattering problem and a vector Helmholtz boundary value problem, where the boundary conditions are expressed solely in terms of the traces  $\gamma$  and $\p_\nu$ of the scattered electric field  $E^s$. This equivalence will serve as the foundation for deriving a robust combined-field formulation in the next section, using the standard Helmholtz boundary integral operators from Calderón calculus.

\begin{theorem} \label{thm:elf}  Fields $E^s\in C^2(\ED,\C^3)\cap C^1(\R^3\setminus\Omega,\C^3)$  and $H^s\in C^1(\ED,\C^3)\cap C(\R^3\setminus\Omega,\C^3)$ form a solution of the scattering problem~\eqref{eq:scatt_problem} 
if and only if $H^s = (\im\omega\mu)^{-1}\curl E^s$ and $E^s$ satisfies the vector Helmholtz equation:
\begin{subequations}\begin{equation}\label{eq:vec_Helm_eqn}
\Delta E^s + k^2E^s=  0\quad\text { in }\  \ED,
\end{equation}
the vector Sommerfeld radiation condition:
\begin{equation}\label{eq:sommerfeld_condition}
\lim _{|\nex| \rightarrow \infty}|x|\left\{(\der E^s)\frac{\nex}{|x|}  - i k E^s\right\} = 0\ \ \text{uniformly in}\ \ \frac{\nex}{|\nex|},
\end{equation}
and the boundary conditions:
\begin{eqnarray}\label{eq:BC_elf}
\Ppar(\gamma E^s)=-\Ppar(\gamma E^i)\quad\text{and}\quad\Pper(\p_\nu E^s+2\mathscr H \gamma E^s) = -\Pper(\p_\nu E^i+2\mathscr H \gamma E^i).
\end{eqnarray}\label{eq:elf_equiv}\end{subequations}
\end{theorem}
\begin{proof}Suppose the fields $E^s$ and $H^s$ solve the scattering problem~\eqref{eq:scatt_problem}. Since $\elf^s$ and $\mgf^s$ satisfy Maxwell’s equations~\eqref{eq:maxwell} in $\ED$ and the Silver--M\"uller radiation condition~\eqref{eq:rad_cond_sm}, they can  be represented by the Stratton--Chu formula~\eqref{eq:SC_int_elf} and~\eqref{eq:SC_int_mgf}, respectively. It then follows the well-known fact that $\elf^s$ satisfies both the Helmholtz equations~\eqref{eq:vec_Helm_eqn} in~$\ED$ and the Sommerfeld radiation condition~\eqref{eq:sommerfeld_condition}, as these properties are inherited from the free-space Helmholtz Green's function~$G$ in each of the integral kernels; see, e.g.,~\cite[Rem. 3.21]{kirsch2016mathematical}. 

The first part of~\eqref{eq:BC_elf}  follows directly from the PEC boundary condition $\nu\times E=0$ on $\Gamma$ which implies  $\Ppar\gamma E=\nu\times\gamma E\times\nu=0$, where $E=E^i+E^s\in C^1(U\setminus\Omega,\C^3)$ denotes the total field. To obtain the second boundary condition in~\eqref{eq:BC_elf}, we  apply Lemma~\ref{lem:neumann_trace} to derive an expression for the normal derivative of $\elf$. This lemma yields the identity
\begin{equation}\label{eq:nor_der_elf}
\oP_\nu\p_\nu E= \gamma(\dive E)\nu-\dive_{\Gamma} (\Ppar\gamma E)\nor-2 \mathscr H(\Pper \gamma E).
\end{equation}

Utilizing the first boundary condition $\Ppar\gamma E=0$ and the fact that $\gamma(\dive E) = 0$, 
we obtain
\begin{align}
\oP_\nu\p_\nu E=&~-2 \mathscr H(\Pper \gamma E).\label{eq:normal_der_elf}
\end{align}
Note that $\gamma(\dive E) = 0$ holds true by virtue of the fact that $\dive E\in C(U\setminus\Omega,\C)$ and $\dive E=-(\im\omega\epsilon)^{-1}\dive\curl H=0$ in $U\setminus\overline\Omega$, which in turn follows from the fact that incident and scattered fields satisfy Maxwell's equations,~\eqref{eq:maxwell_inc} and~\eqref{eq:maxwell}, respectively.
The sought expression for second boundary condition is then obtained directly from  in~\eqref{eq:normal_der_elf} since $\p_\nu E=\p_\nu(E^i+E^s)$.

To prove the converse implication, we first note that the first part of~\eqref{eq:BC_elf} is clearly equivalent to the PEC boundary condition. Furthermore, we observe that, in view of Green's representation formula~\cite[Thm. 3.6]{kirsch2016mathematical}:
\begin{equation}\label{eq:GF_Es}
E^s(x) = \int_\Gamma\left\{\frac{\p G(x,y)}{\p \nu(y)} E^s(y)-G(x,y)\frac{\p E^s(y)}{\p \nu}\right\} \de s(y),\quad x\in\ED.
\end{equation}
---which holds component-wise for $E^s\in C^2(\ED,\C^3)\cap C^1(\R^3\setminus\Omega,\C^3)$  satisfying~\eqref{eq:vec_Helm_eqn} and~\eqref{eq:sommerfeld_condition}---$E^s$ satisfies the Silver-M\"uller radiation condition~\eqref{eq:rad_cond_sm}. This can be straightforwardly verified by observing that the vector fields  
$$
 \frac{\p G(x,y)}{\p \nu(y)} E^s(y) \quad \text{and} \quad  G(x,y) E^s(y),  
\quad x \in \ED, \, y \in \Gamma,
$$  
satisfy the specified condition. We then conclude that to complete the proof, it suffices to show that $\curl H^s = (\im\omega\mu)^{-1}\curl\curl E^s = -\im\omega\epsilon E^s$ in $\ED$, as the other Maxwell's equation is automatically satisfied by construction. 

In view of the vector calculus identity:
  $$\curl\curl E^s = -\Delta E^s + \nabla\dive E^s=k^2E^s+\nabla\dive E^s\quad\text{in }\ED,$$ 
and the fact that $k^2=\omega^2\epsilon\mu$,  it becomes clear that we only need to show that $u:=\dive E^s=0$ in $\ED$.  To prove this, in what follows we establish that $u$ is the unique radiative solution of a homogeneous Dirichlet boundary value problem for the Helmholtz equation in the exterior domain $\ED$. 
Taking the divergence on both sides of~\eqref{eq:GF_Es}  and leveraging the analyticity of $G(\cdot,y)$ in $\R^3\setminus\{x\neq y\}$, we get 
$$
u(x) = \int_\Gamma\left\{\nabla_x\frac{\p G(x,y)}{\p \nu(y)}\cdot E^s(y)-\nabla_xG(x,y)\cdot\frac{\p E^s(y)}{\p \nu}\right\} \de s(y),\quad x\in\ED.
$$
We then conclude from the identity above that $u$ satisfies both the scalar Helmholtz equation 
$$\Delta u+k^2u=0\quad\text{in}\quad \ED,$$ and  the Sommerfeld radiation condition
$$
\lim_{|x|\to\infty}|x|\left\{\nabla u\cdot\frac{x}{|x|}-\im k u\right\} = 0\text{ uniformly in  } \frac{x}{|x|},
$$
by virtue of the fact that  first or higher-order derivatives of $G$ satisfy both properties.  
Therefore, in order to prove that $u=0$ in $\ED$, by the uniqueness of solutions to the exterior Dirichlet problem~\cite[Thm. 3.21]{COLTON:1983}, it suffices to show that the second boundary condition in~\eqref{eq:BC_elf} implies that $\gamma u =\gamma(\dive E^s)=0$. Indeed, from~\eqref{eq:nor_der_elf}, using the first boundary condition $\Ppar\gamma E=0$, we readily get 
$$
\Pper(\p_\nu E+2\mathscr H\gamma E)=\gamma(\dive E)\nu-\dive_{\Gamma} (\Ppar\gamma E)\nor=\gamma(\dive E^s)\nu=0,
$$
where again we used the fact that $\gamma(\dive E^i)=0$ which follows from~\eqref{eq:maxwell_inc}. The proof is now complete.
\end{proof}

Our next theorem establishes the equivalence between the Maxwell PEC scattering problem and an exterior Helmholtz boundary value problem for the magnetic field.

\begin{theorem} \label{lem:equiv_int}  The fields $E^s\in C^1(\ED,\C^3)\cap C(\R^3\setminus\Omega,\C^3)$ and $H^s\in C^2(\ED,\C^3)\cap C^1(\R^3\setminus\Omega,\C^3)$ form a solution of the scattering problem~\eqref{eq:scatt_problem} 
if and only if $E^s = -(\im\omega\epsilon)^{-1}\curl H^s$ and $H^s$ satisfies the vector Helmholtz equation:
\begin{subequations}
    \begin{equation}\label{eq:vec_Helm_eqn_mgf}
\Delta H^s + k^2H^s=  0\quad\text { in }\  \ED,
\end{equation}
the vector Sommerfeld radiation condition:
\begin{equation}\label{eq:sommerfeld_condition_mgf}
\lim _{|\nex| \rightarrow \infty}|x|\left\{(\der H^s)\frac{\nex}{|x|}  - i k H^s\right\} = 0\ \ \text{uniformly in}\ \ \frac{\nex}{|\nex|},
\end{equation}
and the boundary conditions:
\begin{eqnarray}\label{eq:BC_mgf}
\Pper(\gamma H^s)=-\Pper(\gamma H^i)\qquad\text{and}\qquad\Ppar(\p_\nu H^s+\mathscr R \gamma H^s) = -\Ppar(\p_\nu H^i+\mathscr R \gamma H^i).
\end{eqnarray}\label{eq:mgf_equiv}\end{subequations}
\end{theorem}

\begin{proof} As in the proof of Theorem~\ref{thm:elf}, we obtain from Stratton--Chu formula~\eqref{eq:SC_int_mgf} that the scattered magnetic field $\mgf^s$ satisfies the Helmholtz equation~\eqref{eq:vec_Helm_eqn_mgf} in $\ED$, as well as the radiation condition~\eqref{eq:sommerfeld_condition_mgf}.

To establish the validity of the boundary conditions in~\eqref{eq:BC_mgf}, we first note that, from the PEC condition $\nu \times \elf = 0$ on~$\Gamma$ in~\eqref{eq:PEC_BC}, Maxwell's equations~\eqref{eq:maxwell}, and the fact that $H = H^i + H^s \in C^{1}(U \setminus \Omega, \C^3)$, it follows that
\begin{equation}\label{eq:PEC_H}
\nu \times \gamma(\curl \mgf) = 0,
\end{equation}  
and consequently,  
\begin{equation}\label{eq:PEC_BC_H}
\dive_\Gamma(\nu \times \gamma(\curl \mgf)) = 0.
\end{equation}

By taking the trace of the vector calculus identity $\nabla \dive \mgf = \curl^2 \mgf + \Delta \mgf$ which holds in $U\setminus\overline\Omega$, and applying the dot product with the normal, we obtain  
\begin{equation}\label{eq:vector_id_mgf_dot_nu}
\p_\nu(\dive \mgf) = \nu \cdot \gamma(\curl^2 \mgf) - k^2 \nu \cdot \gamma \mgf,
\end{equation}
where we have used that $\Delta H = -k^2 H \in C(U\setminus\overline\Omega, \C^3)$ from where it follows that $\gamma(\Delta H) = -k^2 \gamma H$.
Note also that $\gamma(\curl^2 H) = -\im\omega\epsilon \gamma(\curl E)$ is well-defined, owing to the fact that $\curl E=\im\omega H\in C(U \setminus \Omega, \C^3)$. Next, we employ the surface divergence identity:  
\begin{equation}\label{eq:div_surf_id}
\nu \cdot \gamma(\curl F) = -\dive_\Gamma(\nu \times \gamma F),
\end{equation}
for vector fields $F\in C^1(\oED,\C^3)$, 
which can be derived  from standard vector calculus identities. Specifically, adopting the notation from the proof of Lemma~\ref{lem:neumann_trace}, consider the vector calculus identity 
$\nu_s \cdot \curl F = -\dive\left(\nu_s \times F\right) + (\curl \nu_s) \cdot F,
$
which holds in the tubular neighborhood $U_\delta(\Gamma)$.  Since $\nu_s$ is a gradient field (i.e., the gradient of the signed distance function to the boundary), we have $\curl \nu_s = 0$. Substituting this into the identity above yields 
$
\nu_s \cdot \curl F = -\dive\left(\nu_s \times F\right).
$ The desired surface identity follows by taking the limit as the tubular neighborhood collapses to the surface $\Gamma$.

Therefore, substituting $\p_\nu(\dive \mgf) = 0$ and $\nu \cdot \gamma(\curl^2 H) =-\dive(\nu\times\gamma(\curl H))= 0$—the latter being a consequence of~\eqref{eq:div_surf_id} and~\eqref{eq:PEC_BC_H}—into~\eqref{eq:vector_id_mgf_dot_nu}, we obtain  
\begin{equation} \label{eq:nu_X_H_zero}
\nu \cdot \gamma \mgf = 0,
\end{equation}  
since $k \neq 0$. Thus, we obtain the first part of~\eqref{eq:BC_mgf}, as $\Pper(\gamma \mgf) = (\nu \cdot \gamma (H^s + H^s)) \nu = 0$.

To prove the second part of~\eqref{eq:BC_mgf}, we apply $\Ppar$ to $\p_\nu H\in C(\Gamma,\C^3)$ using Lemma~\ref{lem:neumann_trace}. In detail,  we get 
\begin{equation}
    \Ppar(\p_\nu\mgf) = \Ppar\{\nabla_\Gamma(\nu\cdot\gamma\mgf) - \mathscr R(\gamma\mgf) - \nu\times\gamma(\curl\mgf)\},
\end{equation}
from where the desired result follows by plugging~\eqref{eq:PEC_H} and~\eqref{eq:nu_X_H_zero} into the above identity.

To prove the converse implication we proceed similarly to the proof of Theorem~\ref{thm:elf}. Clearly~$H^s\in C^2(\ED,\C^3)\cap C^1(\oED)$ satisfies the Silver--M\"uller radiation condition~\eqref{eq:rad_cond_sm} since
\begin{equation}\label{eq:mgf_single_double_layer_representation}
H^s(x) = \int_\Gamma\left\{\frac{\p G(x,y)}{\p \nu(y)} H^s(y)-G(x,y)\frac{\p H^s(y)}{\p \nu}\right\} \de s(y),\quad x\in\ED.
\end{equation}
and both vector fields inside the integral above satisfy such condition.  

On the other hand, from the first boundary condition in~\eqref{eq:BC_mgf}, we have $\nu \cdot \gamma H = 0$, which directly implies $\nabla_\Gamma(\nu \cdot \gamma H) = 0$. Applying Lemma~\ref{lem:neumann_trace} to $H$ and utilizing the second boundary condition in~\eqref{eq:BC_mgf}, we obtain  
\begin{equation}\label{eq:aca1}
\Ppar\left(\p_\nu H + \mathscr R \gamma H\right) = \nabla_\Gamma(\nu \cdot \gamma H) - \nu \times \gamma (\curl H) = -\nu \times \gamma (\curl H) = 0.
\end{equation}
From this, the PEC boundary condition $\nu \times \gamma E = -(\im\omega\epsilon)^{-1} \nu \times \gamma (\curl H) = 0$ follows immediately.

It remains to show that $\curl E^s=\im\omega\mu H^s$ in $\ED$.  Consider the identity  
\begin{equation}
    \curl \curl \mgf^s = -\Delta \mgf^s + \nabla \dive \mgf^s = k^2 \mgf^s + \nabla \dive \mgf^s \quad \text{in} \quad \ED,
\end{equation}  
from which the Maxwell equation $\curl \elf^s = -(\im\omega\epsilon)^{-1} \curl^2 H^s = \im\omega\mu \mgf^s$ in $\ED$ follows, provided that $v := \dive \mgf^s = 0$ in $\ED$. To show that this holds, we apply the divergence to~\eqref{eq:mgf_single_double_layer_representation} to get that $v\in C^2(\ED)$ is a radiative solution of the Helmholtz equation: $\Delta v+k^2v=0$ in $\ED$. On the other hand, combining~\eqref{eq:vector_id_mgf_dot_nu} and~\eqref{eq:div_surf_id} we obtain  
$$
\p_\nu\dive \mgf = -\dive_\Gamma(\nu \times \gamma(\curl \mgf)) - k^2 \nu \cdot \gamma \mgf = 0,
$$  
which results since both terms in the middle vanish due to~\eqref{eq:aca1} and the first boundary condition in~\eqref{eq:BC_mgf}. 
Therefore, $\p_\nu v=\p_\nu\dive H^s=\p_\nu\dive(H-H^i)=\p_\nu\dive H=0$, since by assumption the incident field satisfies the Maxwell equations in a neighborhood of $\Gamma$. From the uniqueness of the solution to the exterior Neumann problem~\cite[Thm. 3.25]{COLTON:1983}, we finally conclude that $v=0$ in $\ED$ and thus $\curl E^s=\im\omega\mu H^s$ in $\ED$. This completes the proof.

\end{proof}

\section{The zero frequency limit\label{sec:equiv_lap}}
In the previous section, we established the equivalence between the Maxwell PEC problem and the exterior Helmholtz boundary value problems for the electric and magnetic fields, as given in~\eqref{eq:elf_equiv} and~\eqref{eq:mgf_equiv}, respectively, under the assumption that $k > 0$ (or, equivalently, $\omega > 0$). In the zero-frequency limit ($k \downarrow 0$ or, equivalently, $\omega \downarrow 0$), Maxwell’s equations decouple, giving rise to a different set of equations and boundary conditions. It is therefore interesting to determine whether the equivalences between the Maxwell PEC scattering problem in the zero-frequency limit and certain vector Laplace problems still hold in this setting.

{As it turns out, the zero-frequency limit introduces additional challenges. In~\cite{werner1963perfect_reflection}, it is shown that the time-harmonic scattered electric field $E^s$ admits an analytic extension from $k > 0$ to $k = 0$, denoted by $E_0^s$, which corresponds to the physically correct electrostatic solution. Among the infinitely many possible extensions, the physically relevant one is the solution that satisfies the charge integral conditions:
\begin{equation}\label{eq:int_condition}
\int_{\Gamma_j}\nu\cdot\gamma E_0^s\de s=0,\quad j\in\{1,\ldots,J\},
\end{equation}
where $\Gamma_j$ denote the maximal connected components of $\Gamma$.  In fact, the limiting solution is uniquely determined by these conditions. Without them, the resulting electrostatic problem---considered in what follows---becomes ill-posed: for each of the $J$ connected components, one can construct a nontrivial harmonic field that satisfies the homogeneous electric boundary conditions~\eqref{eq:BC_elf} and still decays at infinity~\cite{werner1963perfect_reflection,werner1983spectral}.} 

The following theorem, which is analogous to Theorem~\ref{thm:elf} in the zero-frequency limit, establishes the equivalence between an electrostatic problem and a vector Laplace boundary value problem. This equivalence holds independently of the integral conditions~\eqref{eq:int_condition}, provided that the surface $\Gamma$ is simply connected.

\begin{theorem}\label{rm:zero_frequency}  Let $U$ be an open set containing $\Gamma = \partial \Omega \subset \mathbb{R}^3$, where $\Gamma$ is assumed to be a $C^{2,\alpha}$-smooth and simply connected surface. Suppose that a vector field $E^i_0 \in C^1(U, \mathbb{C}^3)$ is given, satisfying $\curl E^i_0 = 0$ and $\dive E^i_0 = 0$ in $U$. Let $E^s_0 \in C^2(\ED, \C^3) \cap C^1(\R^3 \setminus \Omega, 
\C^3)$ be a vector field that satisfies the decay condition
\begin{equation}\label{eq:dec_staticE}
\lim_{|x|\to\infty} E^s_0(x) = 0\ \ \text{uniformly in }\frac{x}{|x|}.
\end{equation}

Then, $E_0^s$ is a solution of the electrostatic problem:
\begin{subequations}
\begin{equation}
\curl E^s_0 = 0\quad\text{and}\quad \dive E^s_0=0\quad\text{in}\quad \ED,    
\end{equation}
with the boundary condition
\begin{equation}\label{eq:electrostatic_bc}
\nu\times E^s_0 = -\nu\times E^i_0\quad\text{on}\quad\Gamma,    \end{equation}
 \label{eq:electrostatic_problem_elf}\end{subequations}
 if and only if $E_0^s$ is a solution of the vector Laplace boundary value problem:
\begin{subequations}
    \begin{equation}\label{eq:vec_Laplace}
    \Delta E_0^s=0\quad\text{in}\quad\ED,
\end{equation}
\begin{equation}\label{eq:b_conds_Laplace}
\Ppar(\gamma E^s_0)=-\Ppar(\gamma E^i_0)\quad\text{and}\quad\Pper(\p_\nu E^s_0+2\mathscr H \gamma E^s_0) = -\Pper(\p_\nu E^i_0+2\mathscr H \gamma E^i_0).
\end{equation}\label{eq:vector_laplace}\end{subequations}

    \end{theorem}

\begin{proof}

To prove the direct implication we leverage the formula in~\cite[Eq. 5.12]{COLTON:1983}, which shows that a field $E^s_0$ satisfying~\eqref{eq:electrostatic_problem_elf} can be represented via  
\begin{align}\label{eq:CK_0_freq_rep}
E^s_0(x) =  -\nabla \int_{\Gamma}\nu(y)\cdot E_0^s(y) G_0(x, y) \de s(y) +\curl \int_{\Gamma}\nu(y)\times E_0^s(y) G_0(x, y) \de s(y),\ \ x\in\ED,
\end{align}
where $G_0(x,y) := (4\pi|x-y|)^{-1}$, $x\neq y$, is the fundamental solution of the Laplace equation.  
Therefore, from the form of the integral kernel $G_0$ it follows that $E^s_0$ given by~\eqref{eq:CK_0_freq_rep},  satisfies both the vector Laplace equation~\eqref{eq:vec_Laplace}.

Regarding the boundary conditions, the first part of~\eqref{eq:b_conds_Laplace} is the same as~\eqref{eq:electrostatic_bc}. To derive the second boundary condition in~\eqref{eq:b_conds_Laplace}, we proceed as in the proof of Theorem~\ref{thm:elf} by applying Lemma~\ref{lem:neumann_trace} to the normal derivative of the total field $\elf_0 = \elf^i_0 + \elf^s_0 \in C^1(U\setminus\overline{\Omega}, \C^3)$. Leveraging the first part of~\eqref{eq:b_conds_Laplace} and the properties of $E^i_0$, which ensure that $\Ppar \gamma E_0 = 0$ and $\gamma(\dive E_0)=\gamma\dive (E^s_0+E^i_0) = 0$, the desired identity follows.

To prove the converse implication, assume that $E^s_0$ satisfies all the conditions in~\eqref{eq:vector_laplace}. Applying Green's representation formula for the Laplace equation component-wise, we get~\cite[Thm. 7.12]{Mclean2000Strongly}
\begin{equation}\label{eq:laplace_rep_for}
E^s_0(x) = \int_{\Gamma} \left\{\frac{\p G_0(x,y)}{
\p\nu(y)}E^s_0(y)-G_0(x,y)\frac{\p E^s_0(y)}{\p \nu}\right\}\de s(y),\quad x\in\ED.
\end{equation}
We used here both $\Delta E_0=0$ in $\ED$ and the decay condition~\eqref{eq:dec_staticE}.

Now, since the derivatives of the expression inside the curly brackets in~\eqref{eq:laplace_rep_for} satisfy the Laplace equation component-wise, it follows that $u_0 := \dive E^s_0$ satisfies $\Delta u_0 = 0$ in $\ED$. Moreover, using the asymptotics $|\nabla_x G(x,y)| = O(|x|^{-2})$ and $|\der^2_x G(x,y)| = O(|x|^{-3})$ as $|x| \to \infty$ uniformly in $x/|x|$ and for all $y \in \Gamma$, we find from the Green's formula~\eqref{eq:laplace_rep_for} that $ u_0 \to 0$ uniformly in $x/|x|$ as $|x| \to \infty$. Additionally, as shown in the proof of Theorem~\ref{thm:elf}, the second part of~\eqref{eq:b_conds_Laplace} ensures that $\gamma(\dive E_0) = 0$. Consequently, $\gamma u_0 = \gamma(\dive E^s_0) = \gamma(\dive E_0 - \dive E^i_0) = 0$, by the properties of~$E^i_0$. Therefore, by the uniqueness of solutions to the exterior Laplace boundary value problem with Dirichlet boundary conditions~\cite[Thm. 8.10]{Mclean2000Strongly}, we conclude that $u_0 = \dive E^s_0 = 0$ in $\ED$.
 
It remains to show that $\curl E^s_0=0$ in $\ED$. It readily follows from the identity $0=\Delta E^s_0 = -\curl^2 E^s_0 + \nabla(\dive E^s_0)= -\curl^2 E^s_0 + \nabla u_0$
 that $\curl^2 E^s_0 = 0$  in $\ED.$ {As is well-known, since $\Gamma$ is assumed to be simply connected}, this identity implies the existence of a bounded scalar field potential $v_0\in C^{2}(\ED)$ such that $\curl E^s_0 = \nabla v_0$ in $\ED$. Therefore, $ \Delta v_0 = \dive(\curl E^s_0) = 0$ in  $\ED$.
Moreover, from the boundary conditions~\eqref{eq:b_conds_Laplace} and identity~\eqref{eq:div_surf_id}, it follows that 
\begin{align*}
    \partial_\nu v_0 &= \nu\cdot\nabla v_0 =\nu\cdot\gamma(\curl E^s_0)=\nu\cdot\gamma(\curl (E_0-E^i_0)) = -\nabla_\Gamma(\nu\times \gamma E_0) = 0.
\end{align*}
On the other hand, from Green’s representation formula~\cite[Thm.~7.12]{Mclean2000Strongly}  applied to $v_0$, we obtain that $ v_0(x)-v_\infty=O(|x|^{-1})$ as $|x|\to\infty$, uniformly in $x/|x|$, where $v_\infty$ is a constant. Therefore, by the uniqueness of solutions to the exterior Laplace boundary value problem with Neumann boundary condition~\cite[Thm.~8.18]{Mclean2000Strongly}, it follows that $v_0=v_\infty$  in $\ED$. Consequently, we conclude that $\curl E_0=\nabla v_\infty = 0$ in $\ED$. The proof is now complete.
\end{proof}

{\begin{remark}
It is shown in~\cite{werner1983spectral} that the number of linearly independent solutions to a problem equivalent to the homogeneous version of~\eqref{eq:dec_staticE}–\eqref{eq:vector_laplace} equals the number $J$ of connected components of $\Gamma$, regardless of whether the individual components $\Gamma_j$ are simply connected.
\end{remark}}

Our final theorem in this section establishes the validity of the magnetic vector Helmholtz problem~\eqref{eq:mgf_equiv} in the zero-frequency limit. 
\begin{theorem}\label{rm:mgf_zero_frequency}
Let $U$ be an open set containing $\Gamma = \partial \Omega \subset \mathbb{R}^3$, where $\Gamma$ is assumed to be a $C^{2,\alpha}$-smooth surface. Suppose that a vector field $H^i_0 \in C^1(U, \mathbb{C}^3)$ is given, satisfying $\curl H^i_0 = 0$ and $\dive H^i_0 = 0$ in $U$. Let $H^s_0 \in C^2(\ED, \C^3) \cap C^1(\R^3 \setminus \Omega, \C^3)$ be a vector field that satisfies the decay condition
\begin{equation}\label{eq:rad_cond_Laplace_H}
\lim_{|x|\to\infty} H^s_0(x) = 0\ \ \text{uniformly in }\frac{x}{|x|}.
\end{equation}

 Then, $H_0^s$ is a solution of the magnetostatic problem:
\begin{subequations}
\begin{equation}
\curl H^s_0 = 0\quad\text{and}\quad \dive H^s_0=0\quad\text{in}\quad \ED,    
\end{equation}
with the boundary condition
\begin{equation}\label{eq:magnetostatic_bc}
\nu\cdot H^s_0 = -\nu\cdot H^i_0\quad\text{on}\quad\Gamma,    
\end{equation}
\label{eq:magnetostatic_problem_mgf}
\end{subequations}
if and only if $H_0^s$ is a solution of the vector Laplace boundary value problem:
\begin{subequations}
    \begin{equation}\label{eq:vec_Laplace_H}
    \Delta H_0^s=0\quad\text{in}\quad\ED,
\end{equation}
\begin{equation}\label{eq:b_conds_Laplace_H}
\Pper(\gamma H^s_0)=-\Pper(\gamma H^i_0)\quad\text{and}\quad\Ppar(\p_\nu H^s_0+\mathscr R \gamma H^s_0) = -\Ppar(\p_\nu H^i_0+\mathscr R \gamma H^i_0).
\end{equation}\label{eq:H_vec}
\end{subequations}
\end{theorem}
\begin{proof}
To prove the direct implication can proceed as the proof of Theorem~\ref{rm:zero_frequency} to show via the representation formula~\cite[Eq. 5.12]{COLTON:1983}
\begin{align}\label{eq:CK_0_freq_rep_H}
H^s_0(x) =  -\nabla \int_{\Gamma}\nu(y)\cdot H_0^s(y) G_0(x, y) \de s(y) +\curl \int_{\Gamma}\nu(y)\times H_0^s(y) G_0(x, y) \de s(y),\  x\in\ED,
\end{align}
that $H^s_0$ satisfies both the vector Laplace equation~\eqref{eq:vec_Laplace_H}  and the decay condition~\eqref{eq:rad_cond_Laplace_H}.

Now, regarding the boundary conditions, the first part of~\eqref{eq:b_conds_Laplace_H} is the same as~\eqref{eq:magnetostatic_bc}. To derive the second boundary condition in~\eqref{eq:b_conds_Laplace_H}, we have from  Lemma~\ref{lem:neumann_trace} that to the normal derivative of the total field $H_0 = H^i_0 + H^s_0 \in C^1(U\setminus\overline{\Omega}, \C^3)$ is given by 

\begin{align}
\oP_t\p_\nu H_0 =&~ \nabla_{\Gamma}(\nu \cdot \gamma H_0) - \mathscr{R}(\gamma H_0) - \nu \times \gamma(\curl H_0)\label{eq:bc_H_0_0}\\
=&~ - \mathscr{R}(\gamma H_0).\nonumber
\end{align}

Therefore, applying $\oP_t$ to the identity above using the fact that  $\oP_t(\dive_{\Gamma} (\Ppar \gamma H_0)  \nu) = 0$ we obtain the second boundary condition in~\eqref{eq:b_conds_Laplace_H}.

To prove the converse implication, assume that $H^s_0$ satisfies all the conditions in~\eqref{eq:H_vec}. These conditions allows us to applying Green's representation formula to $H_0^s$ to get~\cite[Thm. 7.12]{Mclean2000Strongly}
\begin{equation}\label{eq:laplace_rep_for}
H^s_0(x) = \int_{\Gamma} \left\{\frac{\p G_0(x,y)}{
\p\nu(y)}H^s_0(y)-G_0(x,y)\frac{\p H^s_0(y)}{\p \nu}\right\}\de s(y),\quad x\in\ED.
\end{equation}

It follows from this formula that $u_0 := \dive H^s_0$ satisfies $\Delta u_0 = 0$ in $\ED$ and that $ u_0 \to 0$ uniformly in $x/|x|$ as $|x| \to \infty$. Moreover, from~\eqref{eq:bc_H_0_0} and the boundary conditions~\eqref{eq:H_vec} we have for the total field
$$0=\oP_t(\p_\nu H_0+ \mathscr{R}(\gamma H_0)) = \nabla_{\Gamma}(\nu \cdot \gamma H_0)  - \nu \times \gamma(\curl H_0)\label{eq:bc_H_0}=- \nu \times \gamma(\curl H_0)\label{eq:bc_H_0}.$$
Then, using the identity~\eqref{eq:vector_id_mgf_dot_nu} with $k=0$,  we arrive at
$$
\p_\nu \dive H_0=-\dive_\nu(\nu \times \gamma(\curl H_0))=0.
$$
Consequently, $\p_\nu u_0 = \p_\nu(\dive H^s_0) = \p_\nu(\dive H_0 - \dive H^i_0) = 0$ using the assumed properties of the incident field. Therefore, by the uniqueness of solutions to the exterior Laplace boundary value problem with Neumann boundary condition~\cite[Thm.~8.18]{Mclean2000Strongly}, we conclude that $u_0 = \dive H^s_0 = 0$ in $\ED$.
 
It remains to show that $\curl H^s_0=0$ in $\ED$.  
{Let $D_r = B_r \setminus \overline{\Omega}$, where $B_r =\{x \in \R^3 : |x| < r\}$ is a ball of radius $r > 0$, assumed to be large enough so that $\Omega \subset B_r$. We denote by $S_r$ the boundary $\partial B_r$ of $B_r$. Applying the first vector Green’s identity~\cite[Eq.~4.11]{COLTON:1983}, we get
\begin{equation}\label{eq:vec_green_id}
\int_{D_r} \overline{H_0^s}\cdot \Delta H_0^s+|\curl H_0^s|^2+|\dive H_0^s|^2 \de x = \int_{\p D_r}  \overline{H_0^s}\cdot(\curl H_0^s\times\nu_r)+(\nu_r\cdot \overline{H_0^s})\dive H_0^s\de s,
\end{equation}
where $\nu_r$ is the unit normal to $\partial D_r = \Gamma \cup S_r$, oriented outward from $D_r$.  Then, using the facts: $\Delta H_0^s=0$ and $\dive H_0^s=0$ in $\ED$ we obtain 
$$
\int_{D_r} \overline{H_0^s}\cdot \Delta H_0^s+|\curl H_0^s|^2+|\dive H_0^s|^2 \de x =\int_{D_r}|\curl H_0^s|^2 \de x.
$$
On the other hand, from the identities  
\begin{align*}
\gamma(\curl H^s_0)\times \nu=&~ \gamma(\curl H_0)\times\nu=-\oP_t(\p_\nu H_0+\mathscr{R}(\gamma H_0))  =0 \quad
\text{and}\\
\gamma(\dive H_0^s) =&~0,
\end{align*}
which follow from~\eqref{eq:b_conds_Laplace_H}, the properties of the incident field $H_0^i$, and the fact that $\dive H_0^s\in C(\R^3\setminus\Omega)$ vanishes in $\ED$, it follows that
\begin{align*}
 \int_{\Gamma} \overline{H_0^s}\cdot(\curl H_0^s\times\nu)+(\nu\cdot \overline{H_0^s})\dive H_0^s\de s=0.
\end{align*}
Therefore,~\eqref{eq:vec_green_id} becomes 
\begin{equation}\label{eq:before_limit}
\int_{D_r} |\curl H_0^s|^2 \de x = \int_{S_r}   \overline{H_0^s}\cdot(\curl H_0^s\times\nu_r)\de s.
\end{equation}
Noting that from the representation formula~\eqref{eq:laplace_rep_for} it holds that $\overline{H_0^s}\cdot(\curl H_0^s\times x/|x|)=O(|x|^{-3})$ as $|x|\to\infty$ uniformly in $x/|x|$, we take the limit as $r \to \infty$ in~\eqref{eq:before_limit}. It then follows that $\|\curl H_0^s\|^2_{L^2(\ED)} = 0$, from which we conclude that $\curl H_0^s = 0$ in $\ED$. This completes the proof.}
\end{proof}

{\begin{remark}\label{rem:uniqueness_magnetostatic}Note that, unlike in Theorem~\ref{rm:zero_frequency}, the theorem above does not assume that $\Gamma$ is simply connected.
As shown in~\cite{werner1983spectral}, however, there exist $P$ linearly independent solutions to a static homogeneous problem equivalent to~\eqref{eq:rad_cond_Laplace_H}–\eqref{eq:H_vec}, where $P = \sum_{j=1}^J P_j$ and $P_j$ denotes the topological genus of the connected component $\Gamma_j$ of~$\Gamma$. Consequently, for multiply connected surfaces, both static magnetic problems lack uniqueness.\end{remark}}

\section{Electromagnetic boundary integral equation formulations using Helmholtz integral operators}\label{sec:BIE_formulations}
In this section, we develop novel wavenumber-robust, second-kind regularized boundary integral equations for solving the electromagnetic scattering problem~\eqref{eq:scatt_problem}.

We begin by introducing the vectorial Helmholtz layer potentials and their associated boundary integral operators. We let $\mathcal S: C^{0,\alpha}(\Gamma,\C^3)\to C^{2}(\R^3\setminus\Gamma,\C^3)$ and $\mathcal D:C^{0,\alpha}(\Gamma,\C^3)\to C^{2}(\R^3\setminus\Gamma,\C^3)$, $\alpha\in(0,1)$, denote the Helmholtz single- and double-layer potentials defined as 
\begin{align}
(\mathcal S\varphi)(\nex) :=&\int_\Gamma G(x,y)\varphi(y)\de s(y)\quad\text{and}\label{eq:SL_pot}\\
(\mathcal D\varphi)(\nex):=&\int_{\Gamma} \frac{\p G(x,y)}{\p \nu(y)}\varphi(y)\de s(y),\quad x\in\R^3\setminus\Gamma,\label{eq:DL_pot}    
\end{align}
for vector density functions $\varphi\in C^{0,\alpha}(\Gamma,\C^3)$. Here and throughout this work, the (vectorial) layer potentials~\eqref{eq:SL_pot} and~\eqref{eq:DL_pot}, along with the associated (vectorial) boundary integral operators introduced below, are to be interpreted component-wide.

It follows from the standard jump conditions for the scalar layer potentials that the Helmholtz vectorial layer potentials satisfy similar jump conditions. Specifically, for a sufficiently regular vector density function $\varphi:\Gamma\to\C^3$, the following relations hold for the single-layer and double-layer potentials as the target point $x\in \R^3\setminus\Gamma$ approaches the boundary $\Gamma$ from the exterior or interior:
\begin{subequations}\label{eq:jumps_LP}\begin{align}
\gamma^\pm (\mathcal S\varphi) =&~\oS\varphi,&
\gamma^\pm (\mathcal D\varphi) =&\pm\frac{1}{2}\varphi+\oK\varphi,\label{eq:jump_Dir}\\
\p_\nu^\pm (\mathcal S\varphi) =&\mp\frac{1}{2}\varphi+ \oK'\varphi,&
\p_\nu^\pm (\mathcal D\varphi) =&~\oT\varphi,\label{eq:jump_Neu}
\end{align}\end{subequations}
where $\gamma^+$ and $\p_\nu^+$ (resp. $\gamma^-$ and $\p_\nu^-$) denote the exterior (resp. interior) Dirichlet and Neumann traces defined in~\eqref{eq:ext_traces} (resp.~\eqref{eq:int_traces}), and $\oS$, $\oK$, $\oK'$ and $\oT$ are the (vector) boundary integral operators: 
\begin{subequations}\begin{align}
\oS:C^{0,\alpha}(\Gamma,\C^3)\to C^{1,\alpha}(\Gamma,\C^3),&\quad \oS\varphi(\nex):=\int_{\Gamma}G(x,y)\varphi(y)\de s(y),\label{eq:single_op}\\
\oK:C^{0,\alpha}(\Gamma,\C^3)\to C^{1,\alpha}(\Gamma,\C^3),&\quad \oK\varphi(\nex) :=\int_{\Gamma}\frac{\p G(x,y)}{\p \nu(y)}\varphi(y)\de s(y),\label{eq:double_op}\\
\oK':C^{0,\alpha}(\Gamma,\C^3)\to C^{0,\alpha}(\Gamma,\C^3),&\quad \oK'\varphi(\nex):= \int_{\Gamma}\frac{\p G(x,y)}{\p \nu(x)}\varphi(y)\de s(y),\label{eq:adj_op}\\
\oT:C^{1,\alpha}(\Gamma,\C^3)\to C^{0,\alpha}(\Gamma,\C^3),&\quad \oT\varphi(\nex) :={\rm f.p.}\!\int_{\Gamma}\frac{\p^2 G(x,y)}{\p\nu(x)\p \nu(y)}\varphi(y)\de s(y),\label{eq:hyper_op}
\end{align}\label{eq:BIOs}\end{subequations}
for $\nex\in\Gamma$, which are well defined component-wise and bounded under the assumption that $\Gamma$ is $C^{2,\alpha}$-smooth. This result directly follows from the properties of the corresponding scalar operators~\cite[Thm. 3.4]{COLTON:2012}. 

It is important to note that the hypersingular integral defining $\oT$ is understood in the sense of the Hadamard finite part (see, e.g., \cite{hsiao2008boundary,hackbusch1995integral}). Note also that $\oK'$ is compact and that the compact embedding of $C^{1,\alpha}(\Gamma)$ into $C^{0,\alpha}(\Gamma)$~\cite[Thm. 3.2]{COLTON:2012} ensures that the operators $\oS$ and $\oK$, defined as mappings from $C^{0,\alpha}(\Gamma, \mathbb{C}^3)$ to $C^{0,\alpha}(\Gamma, \mathbb{C}^3)$, are also compact. Lastly, we observe that $\oT - \oT_0 : C^{0,\alpha}(\Gamma, \mathbb{C}^3) \to C^{0,\alpha}(\Gamma, \mathbb{C}^3)$ is compact, where $\oT_0$ denotes the Laplace hypersingular operator, defined in~\eqref{eq:hyper_op} with the wavenumber $k = 0$~\cite[Thm. 2.31]{COLTON:1983}. 

In the sequel, we also use the subscript 0 for $\oS_0$, $\oK_0$, and $\oK’_0$, which denote the single-layer, double-layer, and adjoint double-layer operators, respectively, as defined in~\eqref{eq:BIOs} for $k = 0$. Moreover, the same notation is used to distinguish the single- and double-layer potentials,~\eqref{eq:SL_pot} and~\eqref{eq:DL_pot}, respectively, corresponding to $k = 0$.
%

\subsection{Electric combined-field-only formulation}\label{sec:electric_field_BIE}
In this section, we derive a novel indirect boundary integral equation formulation for the Maxwell PEC scattering problem~\eqref{eq:scatt_problem}, leveraging its equivalence to the vector Helmholtz problem~\eqref{eq:elf_equiv} for the electric field, as established in Theorem~\ref{thm:elf}. Building on the layer potential framework introduced earlier, we consider a low-rank perturbation of the classical acoustic combined single- and double-layer ansatz~\cite{panich1965question,leis1965dirichletschen,brakhage1965dirichletsche,Burton1971Application}. While this perturbation is not necessary for the well-posedness of the resulting system of BIEs at $k>0$, it plays a crucial role in numerically stabilizing the BIE solutions in the low-frequency regime, as demonstrated by numerical examples in Section~\ref{sec:numerics}. This stabilization enables accurate computations even at arbitrarily small frequencies.

We begin by considering the following modified combined-field ansatz:
\begin{equation}\label{eq:cfie_e_ansatz}
\elf^s(\nex) = (\mathcal D - \im\eta\mathcal S)\varphi(\nex)+\xi\sum_{j=1}^J\nabla\Phi_j(x)\ell_j(\varphi),\quad \nex \in \ED,\ \ \eta\in\R\setminus\{0\},\ \ \xi\in\C,
\end{equation}
where $\varphi \in C^{1,\alpha}(\Gamma, \C^3)$ is a surface vector field not necessarily tangential to $\Gamma$. In what follows we let $\{\Gamma_j\}_{j=1}^J$ denote the set of its maximal connected components. Besides the standard acoustic combined-field potential, the ansatz~\eqref{eq:cfie_e_ansatz} involves the linear functionals
\begin{equation}\label{eq:functional}
\ell_j: C^{0,\alpha}(\Gamma,\C^3) \to \C^3,\quad\ell_j(\varphi) :=  \int_{\Gamma_j}\nu \cdot \left( \frac{1}{2}\Id + \oK - \im\eta\oS \right)\varphi  \de x,
\end{equation}
and
\begin{equation}\label{eq:int_source}
\Phi_j :=  G(\cdot - x^\star_j)\in C^\infty(\R^3\setminus\{x^\star_j\}),
\end{equation}
where the points $x^\star_j \in \Omega$ are placed, respectively, inside the regions enclosed by $\Gamma_j$, $j\in\{1,\ldots,J\}$.

The magnetic field is obtained from the electric field~\eqref{eq:cfie_e_ansatz} via: 
$$H^s=(\im\omega\mu)^{-1}\curl E^s\quad (\omega>0).$$

In the sequel, we let $\Xi\in\C^{J\times J}$ be the matrix with coefficients
\begin{equation}\label{eq:constant_vector}
(\Xi_{i,j}) := \int_{\Gamma_i}  \p_\nu\Phi_j\de s,\quad i,j\in\{1,\ldots,J\}.
\end{equation}

The discussion on the specific role of the parameter $\xi \in \mathbb{C}$ and insights into the modification of the classical acoustic combined-field ansatz~\eqref{eq:cfie_e_ansatz} will be postponed until Remark~\ref{rem:purpose_of_xi} at the end of this section. Meanwhile, we provide the following remark:

\begin{remark}\label{rem:integral_cond_ansatz} Note that, for $k = 0$, the ansatz~\eqref{eq:cfie_e_ansatz} with $\xi=1$ ensures that the integral condition~\eqref{eq:int_condition}---which is necessary for the well-posedness of the exterior electrostatic problem~\eqref{eq:dec_staticE}-\eqref{eq:electrostatic_problem_elf}---is exactly satisfied. Indeed, it is easy to see, using Green's representation,  that for $k=0$ we have $\Xi_{i,j} = -\delta_{i,j}$ for $i,j\in\{1,\ldots,J\}$, so that a static field $E^s_0$ of the form~\eqref{eq:cfie_e_ansatz} satisfies$$
\int_{\Gamma_i}\nu \cdot \gamma E^s_0  \de s = \int_{\Gamma_i}\nu \cdot \left(\frac{1}{2}\Id + \oK_0 - \im\eta \oS_0 \right)\varphi  \de s +\xi\sum_{j=1}^J\Xi_{i,j}\ell_j(\varphi)=\left(1-\xi\right)\ell_i(\varphi)=0
$$
for all $i\in\{1,\ldots,J\}$.
\end{remark}

To streamline notation, we introduce the following operators, which will be used throughout the subsequent sections:
\begin{align}
\oA(\xi): C^{0,\alpha}(\Gamma, \mathbb{C}^3) \to C^{1,\alpha}(\Gamma, \mathbb{C}^3),&\quad \oA(\xi)\varphi=2\Big(\oK - \im\eta \oS+\xi\sum_{j=1}^J\gamma(\nabla\Phi_j)\ell_j\Big)\,\varphi\label{eq:op_A},\\
\oB(\xi): C^{1,\alpha}(\Gamma, \mathbb{C}^3) \to C^{0,\alpha}(\Gamma, \mathbb{C}^3),&\quad \oB(\xi)\varphi= -4\Big(\oT - \im\eta \oK'+\xi\sum_{j=1}^J\p_\nu(\nabla\Phi_j)\ell_j\Big)\, \varphi,\label{eq:op_B}
\end{align}
for $\xi\in\C$. In view of the mapping properties of the Helmholtz integral operators discussed earlier, along with the fact that $\gamma(\nabla\Phi_j) \in C^{2,\alpha}(\Gamma,\C^3)$ and $\p_\nu(\nabla\Phi_j) \in C^{1,\alpha}(\Gamma,\C^3)$ for all $j\in{1,\ldots,J}$, it follows that these operators are bounded.

\begin{remark}
In the sequel, we will often omit the $\xi$-dependence of the operators $\oA$ and $\oB$ for clarity of notation when it is evident from the context. The special case $\xi = 0$, which plays a crucial role in the derivation of the BIE for the magnetic field, will be denoted by $\oA^{(0)} := \oA(0)$ and $\oB^{(0)} := \oB(0)$.
\end{remark}

We also introduce the (point-wise multiplication) operator  defined by
\begin{align}
\oH:C^{m,\alpha}(\Gamma,\C^3)\to C^{0,\alpha}(\Gamma,\C^3),& \quad \oH\varphi= \mathscr H\varphi,\label{eq:H_op}
\end{align}
for $m\in\N_0$ and $\alpha\in[0,1)$, which is bounded by virtue of the fact that  $\mathscr{H} \in C^{0,\alpha}(\Gamma)$ which is in turn a  direct consequence of the fact that $\Gamma$ is of class $C^{2,\alpha}$, $\alpha\in(0,1)$.
Note that $\oH$ commutes with the projectors, i.e., $\oP_\nu \oH=\oH\oP_\nu$ and $\oP_t \oH=\oH\oP_t$.

We then proceed to enforce the boundary conditions for the vectorial Helmholtz problem~\eqref{eq:elf_equiv}, which are the same as those for the vector Laplace problem~\eqref{eq:vector_laplace}. From the boundary condition $\Ppar\gamma\elf^s = -\Ppar\gamma\elf^i$ in~\eqref{eq:BC_elf} and the jump conditions of the Dirichlet trace of  the layer potentials~\eqref{eq:jump_Dir}, we obtain 
\begin{equation}\label{eq:tan_e_cfie}
\Ppar\varphi +\Ppar\oA\varphi = -2\Ppar\gamma\elf^i.
\end{equation}

On the other hand, from the boundary condition in~\eqref{eq:BC_elf} for the normal component, $\Pper(\p_\nu \elf^s + 2\mathscr H \gamma \elf^s) = -\Pper(\p_\nu\elf^i + 2\mathscr H \gamma \elf^i)$, and the jump conditions of both Dirichlet and Nuemann traces of the layer potentials~\eqref{eq:jump_Dir} and~\eqref{eq:jump_Neu}, respectively,  we derive 
\begin{equation}\label{eq:per_e_cfie}
\left(\im\eta\Id + 2\oH\right)\Pper\varphi + \Pper\left(-\frac12\oB + 2\oH \oA\right)\varphi = -2\Pper(\p_\nu\elf^i + 2\oH \gamma\elf^i).
\end{equation}
To combine these two equations into a single expression, we apply to~\eqref{eq:tan_e_cfie} the operator $\im\eta\Id + 2\oH$ and then add the result to~\eqref{eq:per_e_cfie}. This yields 
\begin{subequations}\begin{equation}
\begin{split}
\left(\im\eta\Id + 2\oH\right)\varphi + (\im\eta\Ppar+2\oH)\oA\varphi -\frac12 \Pper\oB\varphi 
=f_E,
\end{split}
\end{equation}
where \begin{equation}\label{eq:E_datum}
f_E:= -2(\im\eta\Ppar+2\oH)\gamma\elf^i - 2\Pper\p_\nu\elf^i. 
\end{equation}
\label{eq:ecfie}\end{subequations}
In the sequel, we refer to this equation as the \emph{Electric Combined-Field-Only Integral Equation} (ECFOIE).

Note that although the fact that the operator $\im\eta\Id + 2\oH$  is invertible on $C^{0,\alpha}(\Gamma,\C^3)$ might suggest that the above formulation is a second-kind integral equation, it does not satisfy the criteria of Fredholm alternative due to the presence of the hypersingular operator $\oT$, which is non-compact. This issue is addressed below in Theorem~\ref{thm:well_poss_E}, where an appropriate regularization is introduced to establish the well-posedness of~\eqref{eq:ecfie}.  We now establish the uniqueness of the solution to~\eqref{eq:ecfie} for all $k>0$.


\begin{lemma}\label{lem:unique}
Let $\xi\in\C$ such that the matrix $
  (\Id_J+\xi\Xi)\in\C^{J\times J} $ is invertible, where $ \Xi $ is  defined in~\eqref{eq:constant_vector} and~$\Id_J$ is the identity matrix in $ \mathbb{C}^{J\times J} $. Then, the electric combined-field-only integral equation~\eqref{eq:ecfie} has at most one solution $ \varphi \in C^{1,\alpha}(\Gamma, \mathbb{C}^3) $ for all $ k > 0 $.
\end{lemma}
\begin{proof} Suppose the statement does not hold. Then, there exists a non-trivial~$ \varphi_0\in C^{1,\alpha}(\Gamma,\C^3)$ such that
 \begin{equation}\label{eq:homoge}
\left(\im\eta\Id + 2\oH\right)\varphi_0 + (\im\eta\Ppar+2\oH)\oA\varphi_0 -\frac12 \Pper\oB\varphi_0 =0.
 \end{equation}
  Let us the define
 $$
 F(x) = (\mathcal D-\im\eta \mathcal S)\varphi_0(x)+\xi \sum_{j=1}^J \nabla \Phi_j(x) \ell_j(\varphi_0),\quad x\in\R^3\setminus\Gamma.
$$
 
By the mapping properties of the single- and double-layer potentials introduced in~\eqref{eq:SL_pot} and~\eqref{eq:DL_pot}, it follows directly that $F \in C^{2}(\mathbb{R}^3 \setminus \Gamma, \mathbb{C}^3)$. 

Letting $F_e$ denote the restriction of $F$ to the exterior domain $\ED$, we clearly have that $F_e$ satisfies both  $\Delta F_e + k^2 F_e = 0$ in $
\ED$ (i.e.,~\eqref{eq:vec_Helm_eqn}) and the Sommerfeld radiation condition~\eqref{eq:sommerfeld_condition}. From the jump conditions of the potentials~\eqref{eq:jump_Dir} and the mapping properties of the single- and double-layer operators~\eqref{eq:single_op} and~\eqref{eq:double_op}, respectively, it follows that $\gamma^+ F_e \in C^{1,\alpha}(\Gamma, \mathbb{C}^3)$. Similarly, from the mapping properties of the adjoint double-layer and hypersingular operators in~\eqref{eq:adj_op} and~\eqref{eq:hyper_op}, respectively, we get $\p_\nu^+ F_e\in C^{0,\alpha}(\Gamma,\C^3)$. Consequently, we have that $F_e \in C^{2}(\ED, \mathbb{C}^3) \cap C^1(\oED, \mathbb{C}^3)$.  Furthermore, from~\eqref{eq:homoge}, $F_e$ also satisfies the boundary conditions~\eqref{eq:BC_elf} with $\gamma E^i=\p_\nu E^i=0$, i.e., 
\begin{equation}\label{eq:homoge_bcs}
\Ppar\gamma^+ F_e = 0\quad\text{and}\quad \Pper(\p_\nu F_e + 2\oH\gamma^+ F_e) = 0.
\end{equation}
  Therefore, in view of the fact that $\curl F_e\in C^{1}(\ED,\C^3)\cap C(\R^3\setminus\Omega,\C^3)$, and since Theorem~\ref{thm:elf} guarantees that  $(F_e,(\im\omega\mu)^{-1}\curl F_e)$ is a solution of the Maxwell PEC scattering problem~\eqref{eq:scatt_problem} with the incident field $\elf^i = 0$, by the uniqueness of solutions to that problem~\cite[Thm. 3.35]{kirsch2016mathematical}, we conclude that $F_e = 0$ in $\mathbb{R}^3 \setminus \Omega$. Consequently, both of its exterior traces vanish, i.e.,  $\gamma^+ F_e = 0$ and $\p^+_\nu F_e = 0$. 
 
Now, turning our attention to the functionals $\ell_i(\varphi_0)$, $i\in\{1,\ldots,J\}$, we note that from the definitions~\eqref{eq:functional},~\eqref{eq:int_source}, and~\eqref{eq:constant_vector}, it  follows that
$$
0=\int_{\Gamma_i}\nu\cdot\gamma^+ F_e\de s = \ell_i(\varphi_0) +\xi\sum_{j=1}^J\Xi_{i,j}\ell_j(\varphi_0). $$
Therefore, the vector $v\in\C^{J}$, $(v_j) =\ell_j(\varphi_0)$, satisfies $(\Id_J+\xi\Xi)v=0$,t and since by assumption $(\Id_J+\xi\Xi)\in \C^{J\times J}$ is invertible,  we conclude that $(v_j)=\ell_j(\varphi_0)=0$ for all $j\in\{1,\ldots,J\}$.

Having established that all the functional $\ell_j$ vanish at $\varphi_0$, we now examine the restriction $F_i$ of $F=(\mathcal D-\im\eta\mathcal S)\varphi_0$ to the interior domain $\Omega$, which satisfies $\Delta F_i + k^2 F_i = 0$ in $\Omega$. From the jump conditions~\eqref{eq:jumps_LP}, we obtain 
\begin{equation}\label{eq:jumps}
-\varphi_0 = \gamma^- F_i - \gamma^+ F_e = \gamma^- F_i, \quad \text{and} \quad -\im\eta\varphi_0 = \p^-_\nu F_i - \p^+_\nu F_e = \p^-_\nu F_i.
\end{equation}
(Note that we used the interior traces $\gamma^-$ and $\p_\nu^-$ introduced in Corollary~\ref{cor:neumann_trace_interior}.)  From ~\eqref{eq:jumps} we get that each coordinate component of the vector field $F_i$ satisfies the homogeneous Robin boundary value problem: Find $u \in C^2(\Omega) \cap C^1(\overline\Omega)$ such that $\Delta u + k^2 u = 0$ in $\Omega$ and $\p_\nu^- u - \im\eta \gamma^- u = 0$ on $\Gamma$. However, it is well known that this problem admits only the trivial solution $u = 0$. Consequently, we have $F_i = 0$ in $\overline{\Omega}$, and therefore, $\varphi_0 = -\gamma^- F_i = 0$. This completes the proof, as we have reached a contradiction. \end{proof}

\begin{remark}  
It follows from the proof above that if $\Gamma$ is connected, i.e., $J=1$, then the uniqueness of solutions to~\eqref{eq:ecfie} holds for all $k > 0$, provided that $\xi \Xi_{1,1} \neq -1$.  
\end{remark}

Unlike the case when $k > 0$, the integral equation~\eqref{eq:ecfie} in the limiting case $k = 0$ admits infinitely many solutions for any choice of $\xi \in \mathbb{C}$. This fundamental lack of uniqueness leads to the following rather unfavorable result:

\begin{proposition} \label{prop:lack_uniqueness}The electric combined-field integral equation~\eqref{eq:ecfie} for $k = 0$ lacks uniqueness, admitting infinitely many solutions $\varphi \in C^{1,\alpha}(\Gamma, \mathbb{C}^3)$ for any choice of the parameter $\xi \in \mathbb{C}$.
\end{proposition}
\begin{proof} Without loss of generality, we assume that $\Gamma$ is connected (i.e., $J=1$) and simply connected. Note that we have $\Xi_{1,1} =-1$. 

Consider first the special case  $\xi = 1$. Let $\varphi_0 \in C^{1,\alpha}(\Gamma, \mathbb{C}^3)$ be such that:
\begin{equation}\label{eq:non_unique_1}
\left(\frac{1}{2}\Id + \oK_0 - \im\eta \oS_0\right)\varphi = \gamma(\nabla\Phi_1),\quad \nabla\Phi_1(x):=-\frac{(x-x^\star_1)}{4\pi |x-x^\star_1|^3},\quad x^\star_1\in\Omega.\end{equation}
Such a density exists and it is not identically zero by virtue of  fact that $\frac12\Id+\oK_0-\im\eta \oS_0$ is invertible on $C^{0,\alpha}(\Gamma,\C^3)$.  
Then, from the uniqueness of the exterior Dirichlet problem for the Laplace equation~\cite[Thm. 6.12]{kress2012linear} it follows that
\begin{equation*}
\nabla\Phi_1(x) = (\mathcal{D}_0 - \im\eta \mathcal{S}_0)\varphi(x) ,\quad x \in \ED.
\end{equation*}
Moreover, for this particular  density we have
$$
\ell_1(\varphi) =  \int_{\Gamma} \nu\cdot \left( \frac{1}{2}\Id + \oK_0 - \im\eta \oS_0 \right)\varphi \de s = \int_{\Gamma}  \p_\nu \Phi_1 \de s = -1\neq0.
$$
This shows that
$$
F(x) := (\mathcal{D}_0 - \im\eta \mathcal{S}_0)\varphi(x) +\xi \nabla\Phi_1(x)\ell_1(\varphi) = 0,\quad x\in\ED. 
$$
Therefore, since $\gamma^+ F = \p_\nu^+ F = 0$, we can deduce, using the identities $\oP_t(\gamma^+ F) = 0$ and $\oP_\nu(\p^+_\nu F + 2\oH \gamma^+ F) = 0$, along with the jump conditions~\eqref{eq:jumps_LP}, that the density $\varphi$ in~\eqref{eq:non_unique_1} satisfies the homogeneous BIE~\eqref{eq:homoge}.

Next, consider the case $\xi \neq 1$. The objective is to construct a non-trivial density $\varphi \in C^{1,\alpha}(\Gamma, \mathbb{C}^3)$ such that the corresponding potential of the form~\eqref{eq:cfie_e_ansatz} yields a non-trivial solution to the exterior vector Laplace problem~\eqref{eq:vector_laplace} with zero boundary conditions. As shown in~\cite{werner1963perfect_reflection}, such a vector Laplace equation solution exists and can be used to construct the desired density.  

Let $\phi_0 \in C^{2}(\ED,\C) \cap C^1(\R^3 \setminus \Omega, \C)$ be the unique solution of the Laplace equation $\Delta \phi_0 = 0$ in $\ED$, with the boundary condition $\phi_0 = 1$ on $\Gamma$ and the asymptotic decay condition $\phi_0(x) = O(|x|^{-1})$ as $|x| \to \infty$, uniformly in $x/|x|$. Define $c := \int_\Gamma \p_\nu \phi_0 \de s \neq 0$. Note that $\nabla \phi_0$ satisfies the exterior vector Laplace problem~\eqref{eq:electrostatic_problem_elf} with zero boundary conditions~\cite{werner1963perfect_reflection}. Moreover, by Theorem~\ref{rm:zero_frequency}, it also solves~\eqref{eq:vector_laplace} with zero boundary conditions. 

Consider then the density $\varphi \in C^{1,\alpha}(\Gamma, \mathbb{C}^3)$ given by the (unique and nontrivial) solution of:
\begin{equation}\label{eq:non_unique_2}
\left(\frac{1}{2}\Id + \oK_0 - \im\eta \oS_0\right)\varphi = \gamma(\nabla\phi_0)-\frac{c\xi}{1-\xi} \gamma(\nabla\Phi_1).\end{equation}
For this density we have 
$$
\ell_1(\varphi) =  \int_{\Gamma} \nu\cdot \left( \frac{1}{2}\Id + \oK_0 - \im\eta \oS_0 \right)\varphi \de s =\int_{\Gamma}
\left(\p_\nu\phi_0- \frac{c\xi}{1-\xi}  \p_\nu \Phi_1 \right)\de s =\frac{c}{1-\xi}.
$$
Then, the potential 
$$
F(x) = (\mathcal{D}_0 - \im\eta \mathcal{S}_0)\varphi(x) + \xi\nabla\Phi_1(x)\ell_1(\varphi),\quad x\in\ED,
$$
satisfies
$$
\gamma F = \gamma(\nabla\phi_0)-\frac{c\xi}{1-\xi}\gamma(\nabla\Phi_1) + \frac{c\xi}{1-\xi}\gamma(\nabla\Phi_1) = \gamma(\nabla\phi_0).
$$

Therefore, again, by the uniqueness of the exterior Dirichlet problem for the Laplace equation~\cite[Thm. 6.12]{kress2012linear}, we conclude that $F = \nabla \phi_0$ in $\ED$. Since,  by Theorem~\ref{rm:zero_frequency}, $F$  solves~\eqref{eq:vector_laplace} with zero boundary conditions, we get $\oP_t(\gamma^+ F) = 0$ and $\oP_\nu(\p^+_\nu F + 2 \oH \gamma^+ F) = 0$. Using the properties of the layer potentials, we obtain the desired result: $\varphi \neq 0$ in~\eqref{eq:non_unique_2} solves the homogeneous  BIE~\eqref{eq:homoge}.

 \end{proof}

We now proceed to present the following lemma, which is instrumental in establishing the existence of solutions to the integral equation~\eqref{eq:ecfie} via the Fredholm alternative:
\begin{lemma} \label{lem:comm_bounded} Let $\oB$ the operator defined in~\eqref{eq:op_B}. Then, the commutator 
\begin{equation}\label{eq:commutator}
[\oP_\nu,\oB] :=\oP_\nu \oB-\oB\oP_\nu= -4[\oP_\nu,\oT]+4\im\eta[\oP_\nu,\oK'] -4\xi\sum_{j=1}^J[\oP_\nu, \p_\nu(\nabla\Phi_j)\ell_j] 
\end{equation}
involving $[\oP_\nu,\oT]$, which is defined by the Cauchy principal value integral 
\begin{equation}\label{eq:principal_value}
[\oP_\nu,\oT]\varphi(x) =\, {\rm p.v.}\!\int_{\Gamma}\frac{\p^2 G(x,y)}{\p\nu(x)\p\nu(y)}\left\{\varphi(y)\cdot\nu(x)\nu(x)-\varphi(y)\cdot\nu(y)\nu(y)\right\}\de s(y),\quad x\in\Gamma,
\end{equation}
 is bounded on $C^{0,\alpha}(\Gamma,\C^3)$.
\end{lemma}

\begin{proof} Formula~\eqref{eq:commutator} follows directly from the definition $\oB = -4(\oT - \im\eta \oK’ +\xi\sum_{j=1}^J\p_\nu (\nabla\Phi_j)\ell_j)$ in~\eqref{eq:op_B}. Given that the operators $\oP_\nu$ and~$\oK’$, introduced in~\eqref{eq:P_perp_op} and~\eqref{eq:adj_op}, respectively, are bounded from $C^{0,\alpha}(\Gamma, \mathbb{C}^3)$ to $C^{0,\alpha}(\Gamma, \mathbb{C}^3)$ under the assumed $C^{2,\alpha}$-regularity of $\Gamma$, it follows that the commutator $[\oP_\nu, \oK’] = \oP_\nu \oK’ - \oK’ \oP_\nu$ is bounded on $C^{0,\alpha}(\Gamma, \mathbb{C}^3)$.

Similarly, since the functionals $\ell_j: C^{0,\alpha}(\Gamma, \mathbb{C}^3) \to \mathbb{C}$, $j\in\{1,\ldots,J\}$, defined in~\eqref{eq:functional} are bounded, and $\p_\nu(\nabla \Phi_j) \in C^{0,\alpha}(\Gamma, \mathbb{C}^3)$, we conclude that the term $\sum_{j=1}^J [\oP_\nu, \p_\nu(\nabla \Phi_j) \ell_j]$ is also bounded on $C^{0,\alpha}(\Gamma, \mathbb{C}^3)$.

To prove the statement, it remains to show that the commutator $[\oP_\nu,\oT] =\oP_\nu \oT - \oT \oP_\nu: C^{0,\alpha}(\Gamma, \mathbb{C}^3) \to C^{0,\alpha}(\Gamma, \mathbb{C}^3)$ is well-defined and bounded, where $\oT$ denotes the hypersingular operator introduced in~\eqref{eq:hyper_op}.

To do so, we first point out the fact that $ \oT - \oT_0: C^{0,\alpha}(\Gamma,\C^3) \to C^{0,\alpha}(\Gamma,\C^3)$ is compact and therefore bounded, with~$\oT_0$ defined in~\eqref{eq:hyper_op}  using  $k = 0$. Consequently, the desired result follows from the fact that the commutator  $[\oP_\nu, \oT_0]: C^{0,\alpha}(\Gamma,\C^3) \to C^{0,\alpha}(\Gamma,\C^3)$, $[\oP_\nu, \oT_0] := \oP_\nu \oT_0 - \oT_0 \oP_\nu$, given by the Cauchy principal value  integral
$$
[\oP_\nu, \oT_0]\varphi(x) = {\rm p.v.} \int_{\Gamma} \frac{\p^2 G_0(x,y)}{\p \nu(x) \p \nu(y)} 
\left\{ \varphi(y) \cdot \nu(x) \nu(x) - \varphi(y) \cdot \nu(y) \nu(y) \right\} \de s(y), \quad x \in \Gamma,
$$
is well-defined and bounded, which is what we proceed to prove in the sequel. 

Clearly, the coordinate components of $[\oP_\nu, \oT_0]\varphi(x)$ can be expressed in terms of functions of the form
\begin{equation}\label{eq:q_2nd_proof}
\oQ_0\phi(x) := {\rm p.v.}\!\!\int_{\Gamma} p(x,y) q(x,y) \phi(y) \de s(y),\quad x\in\Gamma,
\end{equation}
where
\begin{subequations}\begin{align}
p(x,y) &:=\frac{\p^2 G_0(x,y)}{\p \nu(x) \p \nu(y)}= \frac{\nu(x)\cdot\nu(y)}{4\pi|x-y|^3} - \frac{3}{4\pi} \frac{\nu(y)\cdot(x-y)\ \nu(x)\cdot(x-y)}{|x-y|^5},\label{eq:def_p}\\
q(x,y) &:= \nu_i(x)\nu_j(x) - \nu_i(y)\nu_j(y)=[\nu_i(x)-\nu_i(y)]\nu_j(x)+\nu_i(y)[\nu_j(x)-\nu_j(y)],\label{eq:def_q}
\end{align}\end{subequations}
and $\phi\in C^{0,\alpha}(\Gamma)$ could be any of the components of the vector density $\varphi$,
with $\nu_i \in C^{1,\alpha}(\Gamma)$, $i=1,2,3$, denoting the coordinate components of the unit normal vector $\nu \in C^{1,\alpha}(\Gamma,\mathbb{R}^3)$. Therefore, to prove the theorem, it suffices to show that the integral operator $\oQ_0$, given by~\eqref{eq:q_2nd_proof}, is well defined and bounded on $C^{0,\alpha}(\Gamma)$. To do so, we write
$\oQ_0\phi(x) = u(x)+v(x)$ in terms of
$$
u(x):=\phi(x)\ {\rm p.v.}\!\int_{\Gamma}p(x,y)q(x,y)\de s(y)\quad\text{and}\quad v(x):=\int_{\Gamma}p(x,y)q(x,y)\{\phi(y)-\phi(x)\}\de s(y),
$$
and in what follows we show that $u,v\in C^{0,\alpha}(\Gamma)$.

In view of the identities $u=\phi\,\oQ_0\mathbbm{1}$, where $\mathbbm{1}(x):=1$, $x\in\Gamma$,  and 
\begin{equation}\label{eq:Q_1_id}
\oQ_0\mathbbm{1}(x) =\nu_i(x)\nu_j(x) (\oT_0\mathbbm{1})(x)-\oT_0(\nu_i\nu_j)(x)=-\oT_0(\nu_i\nu_j)(x),
\end{equation}
the fact that $\nu_i\nu_j,\mathbbm{1}\in C^{1,\alpha}(\Gamma)$, and the mapping properties of $\oT_0$~\eqref{eq:hyper_op}, it readily follows that $u\in C^{0,\alpha}(\Gamma)$.

On the other hand, to show that $v\in C^{0,\alpha}(\Gamma)$, we resort to~\cite[Lem. 2.10 \& Rem. 2.11]{COLTON:1983}, which states that for $v$ to be H\"older continuous with exponent $\alpha \in (0,1)$, it suffices to establish the existence of  positive constants $c_1$, $c_2$ and $c_3$ such that:
\begin{enumerate}[label=(\alph*)]
\item $|\oQ_0\mathbbm{1}(x)|\leq c_1$ for all $x\in \Gamma$;
\item  $|p(x,y)q(x,y)| \leq c_2|x-y|^{-2}$  for all  $x,y\in\Gamma$;
\item  $|p(x_1,y)q(x_1,y) - p(x_2,y)q(x_2,y)| \leq c_3|x_1-y|^{-3}|x_1-x_2|$,  for all $x_1, x_2 \in \Gamma$, $y \in \Gamma$, satisfying $2\left|x_1-x_2\right| \leq |x_1-y|$. 
\end{enumerate}
Clearly, (a) follows directly from~\eqref{eq:Q_1_id}, which combined with the mapping properties of $\oT_0$,  ensures that $\oQ_0\mathbbm{1} \in C^{0,\alpha}(\Gamma)$. Thus, we can define $c_1 := \|\oQ_0\mathbbm{1}\|_{0,\alpha}$.

 To establish (b), we use the fact that there exists a constant $c_\nu>0$ such that $\nu(y)\cdot(x-y) \leq c_\nu|x-y|^2$ and $\nu(x)\cdot(x-y) \leq c_\nu|x-y|^2$ for all $x, y \in \Gamma$ (see, e.g.,~\cite[Thm. 2.2]{COLTON:1983}). From the definition of $p$ in~\eqref{eq:def_p} and the property $|\nu_i(x)| \leq |\nu(x)| = 1$ for $i = 1, 2, 3$ and all $x \in \Gamma$, there exists a constant $c_p > 0$ such that
\begin{equation}\label{eq:bound_p}
|p(x, y)| \leq \frac{1}{4\pi |x-y|^3} + \frac{3 c_\nu^2}{4\pi |x-y|} \leq \frac{c_p}{|x-y|^3}, \quad \text{for all } x, y \in \Gamma.
\end{equation}
Moreover, since $\nu$ is Lipschitz continuous, there exists a constant $\ell_\nu > 0$ such that $|\nu(x) - \nu(y)| \leq \ell_\nu|x-y|$ for all $x, y \in \Gamma$. Consequently, we obtain
\begin{equation}\label{eq:bound_prod}
|q(x, y)| \leq |\nu_i(x) - \nu_i(y)||\nu_j(x)| + |\nu_i(y)||\nu_j(x) - \nu_j(y)| \leq 2\ell_\nu|x-y|, \quad \text{for all } x, y \in \Gamma.
\end{equation} 
Combining~\eqref{eq:bound_p} and~\eqref{eq:bound_prod}, we can select $c_2 := 2c_p\ell_\nu > 0$ such that
$$
|p(x, y)q(x, y)| \leq \frac{c_2}{|x-y|^2}, \quad \text{for all } x, y \in \Gamma.
$$

Finally, to establish (c), we note that \begin{subequations}\begin{align}
|p(x_1, y)q(x_1, y) - p(x_2, y)q(x_2, y)| =&~|p(x_1, y)q(x_1, y) - p(x_1, y)q(x_2, y) \nonumber \\
&~+ p(x_1, y)q(x_2, y) - p(x_2, y)q(x_2, y)| \nonumber \\
\leq&~|q(x_1, y) - q(x_2, y)||p(x_1, y)| \label{eq:diff_q} \\
&~+ |p(x_1, y) - p(x_2, y)||q(x_2, y)|. \label{eq:diff_p}
\end{align}\label{eq:bound_for_c}\end{subequations}
Using the Lipschitz continuity of $\nu$ we readily find that 
$|q(x_1, y) - q(x_2, y)| = |\nu_i(x_1)\nu_j(x_1) - \nu_i(x_2)\nu_j(x_2)| \leq 2\ell_\nu|x_1 - x_2|$.
It then follows from~\eqref{eq:bound_p} that the term in~\eqref{eq:diff_q} can be bounded as
\begin{equation}\label{eq:1st_bound}
|q(x_1, y) - q(x_2, y)||p(x_1, y)| \leq 2c_p\ell_\nu|x_1 - y|^{-3}|x_1 - x_2|.
\end{equation}

Similarly, to bound the term in~\eqref{eq:diff_p}, we apply the triangle inequality to obtain
\begin{subequations}\begin{align}
|p(x_1,y)-p(x_2,y)|\leq&  \frac{1}{4\pi}\left| \frac{\nu(x_1)\cdot\nu(y)}{|x_1-y|^3}-\frac{\nu(x_2)\cdot\nu(y)}{|x_2-y|^3}\right|+ \label{eq:diff_1}\\
&\frac{3}{4\pi}\left| \frac{\nu(y)\cdot(x_1-y)\ \nu(x_1)\cdot(x_1-y)}{|x_1-y|^5}- \frac{\nu(y)\cdot(x_2-y)\ \nu(x_2)\cdot(x_2-y)}{|x_2-y|^5}\right|\label{eq:diff_2}
\end{align}\label{eq:diff_ps}\end{subequations}

To derive suitable upper bounds for the terms on the right-hand side of~\eqref{eq:diff_ps}, we begin by analyzing the expression  
$\left| |x_1-y|^{-m} - |x_2-y|^{-m} \right|$, $m > 0$, under the condition $2|x_1 - x_2| \leq |x_1 - y|$. We apply the Mean Value Theorem to the function \( x \mapsto |x-y|^{-m} \), which guarantees the existence of a point \(\xi = \sigma x_1 + (1-\sigma)x_2\), with \(\sigma \in (0,1)\), such that  
\begin{equation}\label{eq:MVT}
 \frac{1}{|x_1-y|^{m}} - \frac{1}{|x_2-y|^{m}} = -m \frac{ (\xi-y) \cdot (x_1 - x_2) }{|\xi - y|^{m+2}}.
\end{equation}
Next, we use the triangle inequality to bound \( |x_1 - y| \) in terms of \( |\xi - y| \). Note that  
$$
|x_1-y| \leq |x_1 - \xi| + |\xi - y| = |(1-\sigma)(x_1 - x_2)| + |\xi - y| = (1-\sigma)|x_1 - x_2| + |\xi - y|.
$$
Given the condition \( 2|x_1 - x_2| \leq |x_1 - y| \), we obtain  
$|x_1 - y| \leq \frac{1-\sigma}{2}|x_1 - y| + |\xi - y|,$
which implies
$\frac{1}{2}|x_1 - y| \leq \left( \frac{1 + \sigma}{2} \right)|x_1 - y| \leq |\xi - y|.$
Combining this inequality with~\eqref{eq:MVT}, we conclude that
\begin{equation}\label{eq:bound_diff_MVT}
\left|\frac{1}{|x_1-y|^{m}} - \frac{1}{|x_2-y|^{m}} \right| \leq m \frac{|x_1 - x_2|}{|\xi - y|^{m+1}} \leq m 2^{m+1} \frac{|x_1 - x_2|}{|x_1 - y|^{m+1}}.
\end{equation}

Now, using~\eqref{eq:bound_diff_MVT} with $m=3$ for the term on the right-hand side of~\eqref{eq:diff_1}, we obtain
\begin{equation*}\begin{split}
\left| \frac{\nu(x_1) \cdot \nu(y)}{|x_1 - y|^3} - \frac{\nu(x_2) \cdot \nu(y)}{|x_2 - y|^3} \right| 
\leq&~ \frac{|(\nu(x_1) - \nu(x_2)) \cdot \nu(y)|}{|x_1 - y|^3} 
+ \left| \frac{1}{|x_1 - y|^3} - \frac{1}{|x_2 - y|^3} \right| |\nu(x_2) \cdot \nu(y)| \\
\leq&~  \ell_\nu\frac{ |x_1 - x_2| }{|x_1 - y|^{3}}
+ 48\frac{|x_1-x_2|}{|x_1-y|^4}.
\end{split}\end{equation*}

Similarly,  for the more involved term on the right-hand side of~\eqref{eq:diff_2}, we have
\begin{subequations}\begin{align}
\left|\frac{\nu(y)\cdot(x_1-y)\ \nu(x_1)\cdot(x_1-y)}{|x_1-y|^5}
-\frac{\nu(y)\cdot(x_2-y)\ \nu(x_2)\cdot(x_2-y)}{|x_2-y|^5}\right|\hspace{4cm}\nonumber\\
\leq  \frac{\left|\nu(y)\cdot(x_1-y)\ \nu(x_1)\cdot(x_1-y)-\nu(y)\cdot(x_2-y)\ \nu(x_2)\cdot(x_2-y)\right|}{|x_1-y|^5}+\label{eq:long_term}\\
|\nu(y)\cdot(x_2-y)\ \nu(x_2)\cdot(x_2-y)|\left|\frac{1}{|x_1-y|^5}
-\frac{1}{|x_2-y|^5}\right|.\hspace{2cm}\label{eq:short_term}
\end{align}\end{subequations}
Adding and subtracting $\nu(y)\cdot(x_1-y)\ \nu(x_2)\cdot(x_2-y)$ inside the absolute value in the numerator of  the expression in~\eqref{eq:long_term}, it follows that it can be bounded above by
\begin{align*} \frac{|\nu(y)\cdot(x_1-y)| |\nu(x_1)\cdot x_1-\nu(x_2)\cdot x_2-y\cdot(\nu(x_1)-\nu(x_2))|}{|x_1-y|^5}+
\frac{|\nu(y)\cdot(x_1-x_2)|| \nu(x_2)\cdot(x_2-y)|}{|x_1-y|^5}\\
\leq c_\nu(\ell_\nu(1+d)+2)\frac{|x_1-x_2|}{|x_1-y|^3},\end{align*}
where $d:=\max_{y\in\Gamma}|y|$.

Next, using~\eqref{eq:bound_diff_MVT} with $m=5$, we can bound~\eqref{eq:short_term} as 
$$
|\nu(y)\cdot(x_2-y)\ \nu(x_2)\cdot(x_2-y)|\left|\frac{1}{|x_1-y|^5}
-\frac{1}{|x_2-y|^5}\right|\leq 320\,c_\nu^2\frac{|x_1-x_2|}{|x_1-y|^2}.
$$

Therefore, finally, multiplying~\eqref{eq:diff_ps} by $|q(x_2, y)|$, which in view of~\eqref{eq:bound_prod} and the inequality $|x_2 - y| \leq |x_2 - x_1| + |x_1 - y| \leq \frac{3}{2}|x_1 - y|$, it satisfies $|q(x_2, y)| \leq 2\ell_\nu|x_2 - y| \leq 3\ell_\nu|x_1 - y|$, we get that there exists $\tilde c_3>0$ such that
\begin{align}
|p(x_1, y) - p(x_2, y)||q(x_2, y)| 
\leq& \left(\frac{36\ell_\nu}{\pi|x_1 - y|^{3}}  
+ \frac{3\ell_\nu^2 + 9c_\nu\ell_\nu(\ell_\nu + 2 + d\ell_\nu)}{4\pi |x_1 - y|^2}  
+ \frac{720\ell_\nu c_\nu^2}{\pi |x_1 - y|}\right)|x_1-x_2| \nonumber \\
\leq& \tilde c_3 \frac{|x_1 - x_2|}{|x_1 - y|^3}.\label{eq:2nd_bound}
\end{align}

Using~\eqref{eq:1st_bound} and~\eqref{eq:2nd_bound} to bound the expression on the right-hand side of~\eqref{eq:bound_for_c}, we conclude that we can select \( c_3 := 2c_p \ell_\nu + \tilde{c}_3 \) to establish condition~(c).

With the verification of conditions (a), (b), and (c) required to apply~\cite[Lem.~2.10 \& Rem.~2.11]{COLTON:1983}, we conclude that $v \in C^{0,\alpha}(\Gamma)$. This shows that the operator $\oQ_0:C^{0,\alpha}(\Gamma) \to C^{0,\alpha}(\Gamma)$ in~\eqref{eq:q_2nd_proof} is bounded. The proof is now complete.

\end{proof}

The following corollary for the tangent plane projector $\Ppar$ is used in the next section to establish the well-posedness of a combined-field integral equation for the magnetic field problem~\eqref{eq:mgf_equiv}.

\begin{corollary}\label{cor:comm_par_bounded} The commutator 
$$
[\Ppar,\oB]:C^{0,\alpha}(\Gamma,\C^3)\to C^{0,\alpha}(\Gamma,\C^3),
$$
is bounded.
\end{corollary}
\begin{proof}The proof follows immediately from the algebraic identity $[\Ppar, \oB] = -[\Pper, \oB]$ combined with the result of Lemma~\ref{lem:comm_bounded}.\end{proof}

In order to establish well-posedness and, incidentally, derive a preconditioned integral equation amenable to iterative linear algebra solvers, we resort to the following Calderón identity (see e.g. \cite[Eq.~3.48]{COLTON:1983})
\begin{equation}\label{eq:calderon}\begin{aligned}
\oS\oT=\oK^2-\frac{1}{4}\Id.
\end{aligned}\end{equation}

We are now ready to state and prove the main result of this section:
\begin{theorem}\label{thm:well_poss_E}  
Under the conditions for uniqueness of  Lemma~\ref{lem:unique}, the electric combined-field integral equation~\eqref{eq:ecfie} admits a unique solution $\varphi \in C^{1,\alpha}(\Gamma, \mathbb{C}^3)$.  
\end{theorem}  
\begin{proof} 

We begin by introducing the orthogonal decomposition $\varphi = \varphi_t + \varphi_\nu$, where $\varphi_t = \Ppar\varphi$ and $\varphi_\nu = \Pper \varphi$. The linearly independent equations~\eqref{eq:tan_e_cfie} and~\eqref{eq:per_e_cfie}, which were combined to produce~\eqref{eq:ecfie}, can then be written as:
\begin{subequations}\begin{eqnarray}\label{eq:par_split}
\varphi_t + \oP_t \oA \varphi &=& f_t, \\
\label{eq:perp_split_0}
-2(\im\eta \Id + 2\oH)\varphi_\nu + \oP_\nu (\oB - 4\oH\oA)\varphi &=& f_\nu,
\end{eqnarray}\label{eq:without_reg}\end{subequations}
where the operators $\oA$ and $\oB$ are defined in~\eqref{eq:op_A} and~\eqref{eq:op_B}, respectively, and where we have set 
$$
f_t := -2 \oP_t \gamma E^i\quad\text{and}\quad f_\nu := 4 \oP_\nu (\p_\nu E^i + 2 \oH \gamma E^i).
$$ 
Note that, under the assumed $C^{1,\alpha}(U,\C^3)$-regularity of the incident field $E^i$ (where $\Gamma\subset U$, $U$ open), and leveraging the mapping properties of the traces~\eqref{eq:ext_traces}, projector operators~\eqref{eq:projectors}, and the mean-curvature multiplication operator~$\oH$, it follows that $f_t \in C^{1,\alpha}(\Gamma,\C^3)$ and $f_\nu \in C^{0,\alpha}(\Gamma,\C^3)$.

Using the identities $\oP_\nu \oB\varphi_t = [\oP_\nu, \oB]\varphi_t$ and $\oP_\nu \oB\varphi_\nu = [\oP_\nu, \oB]\varphi_\nu + \oB\varphi_\nu$, we further manipulate equation~\eqref{eq:perp_split_0} to express it as  
\begin{equation}\label{eq:perp_split}  
\oB\varphi_\nu - 2(\im\eta \Id + 2\oH)\varphi_\nu + ([\oP_\nu, \oB] - 4 \oP_\nu \oH \oA)\varphi = f_\nu.  
\end{equation}

To establish the existence of a solution for the coupled system of equations~\eqref{eq:par_split}–\eqref{eq:perp_split}, we rely on the Fredholm alternative. To this end, we proceed to regularize the system by applying the vector single-layer operator $\oS_0: C^{0,\alpha}(\Gamma, \mathbb{C}^3) \to C^{1,\alpha}(\Gamma, \mathbb{C}^3)$ to~\eqref{eq:perp_split}. Specifically, $\oS_0$ denotes the operator defined in~\eqref{eq:single_op}, corresponding to a wavenumber $k=0$. 
In what follows, we denote by $\oK_0$, $\oK'_0$, and $\oT_0$ the corresponding double-layer, adjoint double-layer, and hypersingular operators, respectively, as defined in~\eqref{eq:BIOs}. 

Now, applying $\oS_0$ to~\eqref{eq:perp_split} and using the identities  
\begin{equation}\label{eq:reg_B}
\oS_0\oB = -4\oS_0\oT_0+\oS_0(\oB+4\oT_0) = \Id - 4 \oK^2_0  + \oS_0(\oB+4\oT_0),
\end{equation}
which follow from~\eqref{eq:calderon}, we arrive at  
\begin{equation}\label{eq:perB}
\varphi_\nu + \left(-4 \oK_0^2\oP_\nu + \oS_0\left\{(\oB+4\oT_0)\oP_\nu - 2(\im\eta \Id + 2\oH)\oP_\nu + [\oP_\nu, \oB] - 4 \oP_\nu \oH \oA \right\}\right)\varphi = \oS_0f_\nu.
\end{equation}  
Adding~\eqref{eq:par_split} to~\eqref{eq:perB}, we obtain the following regularized integral equation to be regarded in the space $C^{0,\alpha}(\Gamma, \C^3)$:
\begin{equation}\label{eq:cfie_reg}
\varphi + \oL\varphi = f_t + \oS_0f_\nu,
\end{equation}
where $f_t + \oS_0f_\nu \in C^{1,\alpha}(\Gamma, \C^3)$ and the operator  
\begin{equation}\label{eq:compact_op}
\oL := \oP_t \oA - 4 \oK^2_0\oP_\nu + \oS_0\left\{(\oB+4\oT_0)\oP_\nu - 2(\im\eta \Id + 2\oH)\oP_\nu + [\oP_\nu, \oB] - 4 \oP_\nu \oH \oA \right\}
\end{equation}
is bounded from $C^{0,\alpha}(\Gamma, \C^3)$ to $C^{1,\alpha}(\Gamma, \C^3)$ and therefore compact on $C^{0,\alpha}(\Gamma, \C^3)$, by virtue of the fact that $C^{1,\alpha}(\Gamma, \C^3)$ is compactly embedded in $C^{0,\alpha}(\Gamma, \C^3)$.

Indeed, since both $\oA : C^{0,\alpha}(\Gamma, \C^3) \to C^{1,\alpha}(\Gamma, \C^3)$ and $\oP_t : C^{1,\alpha}(\Gamma, \C^3) \to C^{1,\alpha}(\Gamma, \C^3)$ are bounded, the composition $\oP_t \oA : C^{0,\alpha}(\Gamma, \C^3) \to C^{1,\alpha}(\Gamma, \C^3)$ is also bounded. Similarly, as $\oP_\nu : C^{0,\alpha}(\Gamma, \C^3) \to C^{0,\alpha}(\Gamma, \C^3)$ and $\oK_0 : C^{1,\alpha}(\Gamma, \C^3) \to C^{1,\alpha}(\Gamma, \C^3)$ are bounded, the composition $\oK_0^2 \oP_\nu : C^{0,\alpha}(\Gamma, \C^3) \to C^{1,\alpha}(\Gamma, \C^3)$ is also bounded. On the other hand, the operator inside the curly brackets in~\eqref{eq:compact_op} is bounded because each of the operators involved, namely $\oB +4\oT_0$, $\im\eta \Id + 2\oH$, $\oP_\nu$, $\oH$, $[\oP_\nu, \oB]$, and $\oA$, are bounded in $C^{0,\alpha}(\Gamma, \C^3)$. Specifically, we use the following facts: $\oB + 4\oT_0 = -4(\oT - \oT_0) +4\im\eta\oK' -4\xi\sum_{j=1}^J\p_\nu(\nabla\Phi_j)\ell_j$ is bounded, which follows from the fact that the difference of hypersingular operators $\oT - \oT_0$ has an (improperly) integrable kernel; and the commutator $[\oP_\nu, \oB]$ is bounded, as established in Lemma~\ref{lem:comm_bounded}.  Therefore, since $\oS_0 : C^{0,\alpha}(\Gamma, \C^3) \to C^{1,\alpha}(\Gamma, \C^3)$ is bounded, we conclude that  
$$
\oS_0\left\{(\oB +4\oT_0)\oP_\nu - 2(\im\eta \Id + 2\oH)\oP_\nu + [\oP_\nu, \oB] - 4 \oP_\nu \oH \oA \right\},
$$
being the composition of bounded operators, is itself bounded from $C^{0,\alpha}(\Gamma, \C^3)$ to $C^{1,\alpha}(\Gamma, \C^3)$.

Having established that the regularized integral equation~\eqref{eq:cfie_reg} is of the second kind, the Fredholm alternative (see, e.g.,~\cite[Cor. 317]{COLTON:1983}) implies that to prove well-posedness, it suffices to establish uniqueness. To this end, suppose for the sake of contradiction that uniqueness fails. Then, there exists a nonzero function $\varphi_0 \in C^{0,\alpha}(\Gamma, \C^3)$ such that  
\begin{align}\label{eq:homo_eq}
\varphi_0 + \oL\varphi_0 = 0.  
\end{align}  
By the mapping properties of $\oL$, it follows that $\varphi_0 = -\oL\varphi_0 \in C^{1,\alpha}(\Gamma, \C^3)$.

On other hand, from the definition of the operator $\oL$ in~\eqref{eq:compact_op} and~\eqref{eq:without_reg}, we can express~\eqref{eq:homo_eq} as 
\begin{equation}\label{eq:homo_eq_2}
\mu_{0,t} + \oS_0\mu_{0,\nu} =0
\end{equation}
where
\begin{subequations}\label{eq:homo_dens}
    \begin{align}
\mu_{0,t} :=&\oP_t(\Id+\oA)\varphi_0\in C^{1,\alpha}(\Gamma,\C^3)\\
\mu_{0,\nu} :=&\oP_\nu(-2\im\eta \Id - 4\oH +  \oB - 4\oH\oA)\varphi_0\in C^{0,\alpha}(\Gamma,\C^3).
\end{align}
\end{subequations}
Then, clearly, since $\mu_{0,t}$ and $\mu_{0,\nu}$ are orthogonal to each other and they are related via~\eqref{eq:homo_eq_2},  we arrive at 
$$
0=\int_\Gamma\mu_{0,t}\cdot\overline{\mu_{0,\nu}}\de s=-\int_\Gamma(\oS_0\mu_{0,\nu})\cdot\overline{\mu_{0,\nu}}\de s.
$$ 
Now, using the fact that  $\oS_0: H^{-1/2}(\Gamma,\C^3)\to H^{1/2}(\Gamma,\C^3)$ is coercive~\cite[Thm. 7.6]{Mclean2000Strongly} and $\mu_{0,\nu}\in C^{0,\alpha}(\Gamma,\C^3)\subset H^{-1/2}(\Gamma,\C^3)$, we get 
$$
0=\real\int_\Gamma (\oS_0\mu_{0,\nu})\cdot\overline{\mu_{0,\nu}}\de s=\real\langle\oS_0\mu_{0,\nu},\mu_{0,\nu}\rangle_{\Gamma}\gtrsim \|\mu_{0,\nu}\|^2_{H^{-1/2}(\Gamma,\C^3)}\gtrsim \|\mu_{0,\nu}\|^2_{H^{1/2}(\Gamma,\C^3)}\geq 0.
$$
We conclude then that $\mu_{0,\nu} = 0$, and since $\mu_{0,t} = -\oS_0 \mu_{0,\nu}$, it follows that $\mu_{0,t} = 0$. Given that $\mu_{0,\nu}=0$ and $\mu_{0,t}=0$, we find from~\eqref{eq:homo_dens} and~\eqref{eq:without_reg} that  $\varphi_0$ is solution to the (non-regularized) homogeneous equation~\eqref{eq:homoge}. However, Lemma~\ref{lem:unique} guarantees that the only solution to~\eqref{eq:homoge} is $\varphi_0 = 0$. This contradiction demonstrates that $\varphi_0$ must be the trivial solution, thereby proving uniqueness. 

Finally, to complete the proof, we need to show that the unique solution $\varphi \in C^{0,\alpha}(\Gamma, \C^3)$ of~\eqref{eq:cfie_reg} is actually an element of $C^{1,\alpha}(\Gamma, \C^3)$. However, since $\oL\varphi \in C^{1,\alpha}(\Gamma, \C^3)$ and $f_t + \oS_0 f_\nu \in C^{1,\alpha}(\Gamma, \C^3)$, it readily follows from itself~\eqref{eq:cfie_reg} that 
$$
\varphi = -\oL\varphi + f_t + \oS_0 f_\nu \in C^{1,\alpha}(\Gamma, \C^3).
$$

The proof is now complete.

\end{proof}

Based on the proof of Theorem~\ref{thm:well_poss_E}, we define the second-kind integral equation  
\begin{subequations}\begin{equation}
\oP_t (\Id + \oA)\varphi + \oS_0 \oP_\nu (-2i \eta \Id - 4 \oH +\oB - 4 \oH \oA ) \varphi = f,
\end{equation}
where  
\begin{equation}
f := -2 \oP_t \gamma E^i + 4 \oS_0 \oP_\nu \big( \p_\nu E^i + 2 \oH \gamma E^i \big),
\end{equation}\label{eq:R-ECFOIE}\end{subequations}
as the \emph{Regularized Electric Combined-Field-Only Integral Equation (R-ECFOIE)}. This equation is well-posed in $C^{0,\alpha}(\Gamma, \C^3)$ under the conditions of Theorem~\ref{thm:well_poss_E} and is equivalent to the ECFOIE~\eqref{eq:ecfie}.

\begin{remark}[Role and Selection of $\xi$]\label{rem:purpose_of_xi}
As anticipated, the lack of uniqueness in the zero-frequency limit of both ECFOIE~\eqref{eq:ecfie} and R-ECFOIE~\eqref{eq:R-ECFOIE} (see Proposition~\ref{prop:lack_uniqueness}) leads to increasingly ill-conditioned linear systems when discretized as $k \downarrow 0$. This phenomenon is how the low-frequency breakdown manifests itself in our approach. However, as demonstrated in Section~\ref{sec:numerics}, this ill-conditioning primarily degrades solution accuracy rather than significantly affecting the convergence of the iterative linear solver (e.g., GMRES), particularly in the R-ECFOIE case.

In practice, however, the parameter choice $\xi = 1$ plays a crucial role in mitigating these issues. In the zero-frequency limit, this choice ensures the well-posedness of the electrostatic boundary value problem by enforcing the integral condition~\eqref{eq:int_condition} (see Remark~\ref{rem:integral_cond_ansatz}). Consequently, it significantly alleviates the low-frequency breakdown, enabling the numerical solution even at the zero-frequency limit for simply connected surfaces. This rationale underlies the modification of the classical combined field potentials in our ansatz~\eqref{eq:cfie_e_ansatz}.
\end{remark}

\subsection{Magnetic combined-field-only formulation}\label{sec:magnetic_field_BIE}
In this section, we develop a combined-field formulation for the exterior Helmholtz boundary value problem governing the magnetic field, as described by equation~\eqref{eq:mgf_equiv}. To this end, we adopt the classical acoustic combined-field ansatz~\cite{panich1965question,leis1965dirichletschen,brakhage1965dirichletsche,Burton1971Application}:
\begin{equation}\label{eq:cfie_h_ansatz}
H^s(\nex) = (\mathcal D - \im\eta\mathcal S)\psi(\nex), \quad \nex \in \ED, \quad \eta \neq 0,
\end{equation}
where $\psi \in C^{1,\alpha}(\Gamma, \C^3)$. Note that this is simply the vector form of the classical acoustic combined-field potential. It differs from the one used for the electric field in that it does not include the term involving the parameter $\xi \in \C$. 

As in the equivalence established in Theorem~\ref{lem:equiv_int}, the electric field is obtained from the magnetic field~\eqref{eq:cfie_h_ansatz} via the relation:
$$ E^s = -(\im\omega\epsilon)^{-1} \curl H^s\quad (\omega>0). $$

From the boundary condition $\Pper\gamma H^s=-\Pper\gamma H^i$ in~\eqref{eq:BC_mgf} and the jump conditions for the layer potentials~\eqref{eq:jumps_LP}, we obtain 
\begin{equation}\label{eq:tan_h_cfie}
\Pper\psi + \Pper\oA^{(0)}\psi = -2\Pper\gamma H^i,
\end{equation}
where $\oA^{(0)} := \oA(0)$, with the operator $\oA$ defined in~\eqref{eq:op_A}.

On the other hand, from the boundary condition for the normal component $\Ppar(\p_\nu H^s + \mathscr R \gamma H^s) = -\Ppar(\p_\nu H^i + \mathscr R \gamma H^i)$ in~\eqref{eq:BC_mgf} and the jump condition~\eqref{eq:jumps_LP}, we get 
\begin{equation}\label{eq:per_h_cfie}
\left(\im\eta\Id + \oR\right)\Ppar\psi + \Ppar\left(-\frac12\oB^{(0)}+ \oR\oA^{(0)}\right)\psi = -2\Ppar(\p_\nu H^i + \oR \gamma H^i).
\end{equation}
where we used the (multiplication) operator $\oR$ defined as
\begin{align}
\oR:C^{m,\alpha}(\Gamma,\C^3)\to C^{0,\alpha}(\Gamma,\C^3),& \quad \oR\varphi= \mathscr R\varphi,\label{eq:R_op}
\end{align}
for $m\in\N_0$ and $\alpha\in[0,1)$, which is bounded by virtue of the fact that $\mathscr{R} \in C^{0,\alpha}(\Gamma,\R^{3\times 3})$ since $\Gamma$ is of class $C^{2,\alpha}$, $\alpha\in(0,1)$. Note we have used above the fact that $\oR$ commutes with the projector $\Ppar$, i.e., $\oP_t \oR=\oR\oP_t$. On the other hand, the operator $\oB^{(0)}$ in~\eqref{eq:per_h_cfie} is given by  $\oB^{(0)} := \oB(0)$, with the operator $\oB$ defined in~\eqref{eq:op_B}.

We can then combine~\eqref{eq:tan_h_cfie} and~\eqref{eq:per_h_cfie} into a single equation by  applying to~\eqref{eq:tan_h_cfie} the invertible operator $\eta\Id + \oR$ and then add the resulting expression to~\eqref{eq:per_h_cfie}. This yields 
\begin{subequations}\begin{equation}
\begin{aligned}
\left(\im\eta\Id +  \oR\right)\psi + (\im\eta\Pper+\oR)\oA^{(0)}\psi  -\frac12\Ppar\oB^{(0)}\psi 
= g_H,
\end{aligned}\label{eq:mcfie_eq}
\end{equation}
where
\begin{equation}
g_H:=-2(\im\eta\Pper+\oR)\gamma H^i -2\Ppar\p_\nu H^i.
\end{equation}
\label{eq:mcfie}\end{subequations}
In the sequel we refer to this equation as the \emph{Magnetic Combined-Field-Only Integral Equation} (MCFOIE).

As in the case of the electric field equation~\eqref{eq:ecfie}, the combined-field approach renders the magnetic field equation~\eqref{eq:mcfie} uniquely solvable for all $k>0$, as we show in what follows. Interestingly, in contrast to the ECFOIE~\eqref{eq:ecfie}, we can establish a uniqueness result for~\eqref{eq:mcfie} in the limit case $k=0$, provided that the additional assumption of $\Gamma$ being simply connected is satisfied.
 
\begin{lemma}\label{lem:unique_mgf} 
The magnetic combined-field integral equation~\eqref{eq:mcfie} admits at most one solution $\psi \in C^{1,\alpha}(\Gamma, \mathbb{C}^3)$ for all $k > 0$. Furthermore, assuming that $\Gamma$ is simply connected, the equation also admits at most one solution $\psi \in C^{1,\alpha}(\Gamma, \mathbb{C}^3)$ when $k = 0$.
\end{lemma}
\begin{proof}

Suppose there exists a nontrivial vector density $\psi_0 \in C^{1,\alpha}(\Gamma,\C^3)$ that solves~\eqref{eq:mcfie} for $k\geq0$ and $g_H = 0$. Let 
$$
F(x) = (\mathcal{D} - \im\eta \mathcal{S})\psi_0(x), \quad x \in \R^3 \setminus \Gamma.
$$
Suppose $k>0$.  Since the restriction $F_e\in C^{2}(\ED,\C^3)\cap C^{1,\alpha}(\R^3\setminus\Omega,\C^3)$ of $F$  to the exterior domain satisfies all the conditions in~\eqref{eq:mgf_equiv}, it follows from Theorem~\ref{lem:equiv_int} that $(-\frac{1}{\im\omega\epsilon}\curl F_e,F_e)$ is a solution of the PEC scattering problem~\eqref{eq:scatt_problem} with incident field $E^i = 0$. Then, from the uniqueness of solutions to the PEC scattering problem~\cite{kirsch2016mathematical}, we have that $F$ vanishes in the exterior domain $\R^3 \setminus \Omega$, and hence 
$\gamma^+ F = 0$ and $\p_\nu^+ F = 0.$ 

Suppose now that $k=0$. We then have that $F$ satisfies the condition $F = O(|x|^{-1})$ as $|x| \to \infty$, uniformly in $x/|x|$, and it also satisfies the homogeneous conditions in~\eqref{eq:H_vec}, i.e.,
$$
\Delta F = 0 \quad \text{in } \ED, \quad \oP_\nu\gamma^+F = 0, \quad \oP_t(\p^+_\nu F + \mathscr{R} \gamma^+ F) = 0.
$$
By Theorem~\ref{rm:mgf_zero_frequency}, we conclude that $F$ also satisfies the static conditions of the magnetostatic problem~\eqref{eq:magnetostatic_problem_mgf} with $H^i_0 = 0$, i.e.,
$$
\curl F = 0, \quad \dive F = 0 \quad \text{in } \ED, \quad \nu \cdot F = 0 \quad \text{on } \Gamma.
$$
Since this static problem admits at most one solution that decays as $O(|x|^{-1})$ when $\Gamma$ is simply connected~\cite{werner1966small_frequencies}, we conclude that $F = 0$ in  $\ED$. From this, it follows that $\gamma^+ F = 0$ and $\p_\nu^+ F = 0$ in the case $k=0$

Finally, using the jump relations of the layer potentials~\eqref{eq:jumps_LP}, we obtain that $\gamma^- F = -\psi_0$ and $\p_\nu^- F = -\im\eta \psi_0$, implying that $\p_\nu^- F - \im\eta \gamma^- F = 0$. Therefore, since $\Delta F + k^2 F = 0$ in $\Omega$, by the uniqueness of the interior Robin problem for both the Helmholtz and Laplace equations, we conclude that $F = 0$ in $\overline{\Omega}$. Consequently, $\psi_0 = 0$, which leads to a contradiction. The proof is complete.
\end{proof}

As in the case of the electric field equation in Section~\ref{sec:electric_field_BIE}, we now proceed to regularize~\eqref{eq:mcfie_eq} via the application of the single-layer operator $\oS_0$, this time applied to the tangential field components.
\begin{theorem}\label{thm:well_poss_H} There exists a unique solution $\psi \in C^{1,\alpha}(\Gamma, \C^3)$ to the magnetic combined-field-only integral equation~\eqref{eq:mcfie} for all $k>0$. Furthermore, assuming that $\Gamma$ is simply connected, the equation has a unique solution $\psi \in C^{1,\alpha}(\Gamma, \mathbb{C}^3)$ when $k = 0$.
\end{theorem}
\begin{proof} 
We start by introducing the orthogonal decomposition $\psi=\psi_t+\psi_\nu$ where $\psi_t=\oP_t\psi$ and $\psi_\nu=\oP_\nu\psi$. Equations~\eqref{eq:tan_h_cfie} and~\eqref{eq:per_h_cfie}  that make up~\eqref{eq:mcfie} can then be written as 
\begin{subequations}\begin{eqnarray}\label{eq:par_split_h}
\psi_\nu + \Pper\oA^{(0)}\psi &=& g_\nu,\\
\label{eq:perp_split_0_h}
-2(\im\eta \Id+\oR)\psi_t + \Ppar (\oB^{(0)}-2 \oR \oA^{(0)})\psi  &=& g_t,
\end{eqnarray}\end{subequations}
where we have set 
$$
g_\nu:=-2\Pper \gamma H^i\quad\text{and}\quad g_t:=4\Ppar(\p_\nu H^i+\oR \gamma H^i),$$
which satisfy $g_\nu\in C^{1,\alpha}(\Gamma,\C^3)$ and $g_t\in C^{0,\alpha}(\Gamma,\C^3)$ by the properties of the incident field.

From the identities  $\Ppar\oB^{(0)}\psi_\nu=[\Ppar,\oB^{(0)}]\psi_\nu$ and $\Ppar \psi_t= [\Ppar,\oB^{(0)}]\psi_t+\oB^{(0)}\psi_t$, equation~\eqref{eq:perp_split_0_h} can be recast~as
\begin{equation}\label{eq:perp_split_h}
\oB^{(0)}\psi_t-2(\im\eta \Id+ \oR)\psi_t +  ([\Ppar,\oB^{(0)}]-2\Ppar\oR  \oA^{(0)})\psi = g_t.
\end{equation}

We proceed to regularize the hypersingular operator $\oB^{(0)}$ in~\eqref{eq:perp_split_h} by applying the Laplace single-layer operator $\oS_0: C^{0,\alpha}(\Gamma, \mathbb{C}^3) \to C^{1,\alpha}(\Gamma, \mathbb{C}^3)$. Them, from the identity~\eqref{eq:reg_B}
and the definition of the commutator $[\oP_t,\oB]$ in~Corollary~\ref{cor:comm_par_bounded} we get 
\begin{equation}\label{eq:perB_h}
\psi_t + \left(-4 \oK_0^2\oP_t + \oS_0\left\{(\oB^{(0)}+4\oT_0)\oP_t - 2(\im\eta \Id + \oR)\oP_t + [\oP_t, \oB^{(0)}] - 2 \oP_t \oR \oA^{(0)} \right\}\right)\psi = \oS_0g_t.
\end{equation}
 Adding~\eqref{eq:par_split_h} to~\eqref{eq:perB_h}, we obtain the following equation, which should be considered in the space $C^{0,\alpha}(\Gamma,\C^3)$:
\begin{equation}\label{eq:cfie_reg_h}
\psi + \oM \psi = g_\nu + \oS_0 g_t,
\end{equation}
where $g_\nu + \oS_0 g_t \in C^{1,\alpha}(\Gamma,\C^3)$, and the operator 
\begin{equation}\label{eq:compact_op_h}
\oM :=\oP_\nu \oA^{(0)} - 4 \oK^2_0\oP_t + \oS_0\left\{(\oB^{(0)}+4\oT_0)\oP_t - 2(\im\eta \Id + \oR)\oP_t + [\oP_t, \oB^{(0)}] - 2 \oP_t \oR \oA^{(0)} \right\}
\end{equation}
is bounded from $C^{0,\alpha}(\Gamma,\C^3)$ to $C^{1,\alpha}(\Gamma,\C^3)$ and compact on $C^{0,\alpha}(\Gamma,\C^3)$. This follows directly from Corollary~\ref{cor:comm_par_bounded} and the mapping properties of the other operators, in a manner analogous to the discussion in the proof of Theorem~\ref{thm:well_poss_E}.

Being~\eqref{eq:cfie_reg_h} of the second kind, the Fredholm alternative (see, e.g.,~\cite[Cor. 317]{COLTON:1983}) implies that to prove well-posedness, it suffices to establish uniqueness.  To this end, suppose, for contradiction, that there exists a non-trivial function $\psi_0 \in C^{0,\alpha}(\Gamma, \C^3)$ such that  
\begin{align}
\psi_0 + \oM\psi_0 = 0.  
\end{align}  
By the mapping properties of $\oM$, we obtain that $\psi_0 = -\oM \psi_0 \in C^{1,\alpha}(\Gamma, \C^3)$. However, as in the proof of Theorem~\ref{thm:well_poss_E}, the coercivity of $\oS_0$ implies that this contradicts the uniqueness result established in Lemma~\ref{lem:unique_mgf}. This contradiction then implies that $\psi_0$ must be zero, proving uniqueness. 

Finally, from the identity $\psi = -\oM \psi + g_\nu + \oS_0 g_t$ we conclude that the unique solution $\psi \in C^{0,\alpha}(\Gamma, \C^3)$ of~\eqref{eq:cfie_reg_h} belongs to $C^{1,\alpha}(\Gamma,\C^3)$. The proof is now complete.
\end{proof}

Based on the proof of Theorem~\ref{thm:well_poss_H}, we define the second-kind integral equation  
\begin{subequations}\begin{equation}
\oP_\nu (\Id + \oA^{(0)})\psi - \oS_0 \oP_t \big( 2i \eta \Id + 2 \oR - \oB^{(0)} + 2 \oR \oA^{(0)} \big) \psi = g,
\end{equation}
where  
\begin{equation}
g := -2 \oP_\nu \gamma H^i + 4 \oS_0 \oP_t \big( \p_\nu H^i +  \oR \gamma H^i \big),
\end{equation}\label{eq:R-MCFOIE}\end{subequations}
as the \emph{Regularized Magnetic Combined-Field-Only Integral Equation (R-MCFOIE)}. This equation is well-posed in $C^{0,\alpha}(\Gamma, \C^3)$ and is equivalent to the MCFOIE~\eqref{eq:mcfie} for all $k>0$.

\begin{remark}\label{rem:use_mag_BIE}
A clear advantage of the MCFOIE~\eqref{eq:mcfie} and the associated R-MCFOIE~\eqref{eq:R-MCFOIE} over their electric field counterparts,\eqref{eq:ecfie} and~\eqref{eq:R-ECFOIE}, is that the magnetic BIEs do not suffer from the so-called low-frequency breakdown when the surface $\Gamma$ is simply connected, as established in Theorem~\ref{thm:well_poss_H}. This allows for the use of the simpler ansatz~\eqref{eq:cfie_h_ansatz}, leading to more efficient implementations and numerical BIE discretizations that remain well-conditioned for all $k \geq 0$.
\end{remark}
\begin{remark}\label{rem:use_stratton_chu}
Yet another potential advantage of the MCFOIE~\eqref{eq:mcfie} and the associated R-MCFOIE \eqref{eq:R-MCFOIE} over their electric field counterparts,~\eqref{eq:ecfie} and~\eqref{eq:R-ECFOIE}, respectively, is that the induced surface electric currents $J = \nu \times \gamma H$, which are relevant in many application, can be obtained directly from the vector surface density $\psi \in C^{1,\alpha}(\Gamma, \mathbb{C}^3)$ solution of the BIEs. In fact, using the representation formula~\eqref{eq:cfie_h_ansatz} and the mapping properties of the layer potentials~\eqref{eq:jump_Dir}, we readily obtain  
\begin{equation}\label{eq:elec_currents}
    \gamma H^s = \psi + 2(\oK - \im\eta \oS)\psi,
\end{equation}
from which it follows that the induced electric currents are given by 
$$
J = \nu \times H = \nu \times \big\{\psi + 2(\oK - \im\eta \oS)\psi + \gamma H^i\big\}.
$$
Moreover, using~\eqref{eq:elec_currents}, both the electric and magnetic scattered fields can be evaluated everywhere in $\ED$ via the Stratton–Chu representation formula~\eqref{eq:stratton-chu}.
\end{remark}

We finish this section with the following corollary:
\begin{corollary}\label{eq:reg_solution}Let the incident electromagnetic field $E^i, H^i \in C^{1,\alpha}(U, \C^3)$ satisfies~\eqref{eq:maxwell_inc}, where $U \subset \R^3$ is an open set containing $\Gamma = \partial\Omega$, a closed surface of class $C^{2,\alpha}$ for some $\alpha \in (0,1]$. Then, the unique solution to the PEC scattering problem~\eqref{eq:scatt_problem} satisfies $E^s, H^s \in C^2(\R^3 \setminus \overline{\Omega}, \C^3) \cap C^{1,\alpha}(\R^3 \setminus \Omega, \C^3)$.
\end{corollary}
\begin{proof}
Let $\varphi \in C^{1,\alpha}(\Gamma, \C^3)$ be the unique solution of the electric integral equation~\eqref{eq:R-ECFOIE}. Then the electric field $E^s \in C^{2}(\ED, \C^3) \cap C^{1,\alpha}(\Gamma, \C^3)$, defined by the combined-field potential~\eqref{eq:cfie_e_ansatz}, satisfies the electric vector Helmholtz problem~\eqref{eq:elf_equiv}. Consequently, by Theorem~\ref{thm:elf}, the pair $(E^s, \widetilde H^s)$, with $\widetilde H^s := (\im\omega\mu)^{-1} \curl E^s \in C^1(\ED, \C^3) \cap C(\R^3 \setminus \Omega, \C^3)$, forms a solution of the PEC scattering problem~\eqref{eq:scatt_problem}.

Similarly, let $\psi \in C^{1,\alpha}(\Gamma, \C^3)$ be the unique solution of~\eqref{eq:R-MCFOIE}. Then the magnetic field $H^s \in C^{2}(\ED, \C^3) \cap C^{1,\alpha}(\R^3 \setminus \Omega, \C^3)$, given by the potential~\eqref{eq:cfie_h_ansatz}, solves the magnetic vector Helmholtz problem~\eqref{eq:mgf_equiv}. Therefore, by Theorem~\ref{lem:equiv_int}, the pair $(\widetilde E^s, H^s)$, with $\widetilde E^s := -(\im\omega\epsilon)^{-1} \curl H^s \in C^1(\ED, \C^3) \cap C^{0,\alpha}(\R^3 \setminus \Omega, \C^3)$, also solves the PEC scattering problem~\eqref{eq:scatt_problem}.

By the uniqueness of $C^1(\ED,C^3)\cap C(\R^3\setminus\Omega,C^2)$-regular solutions to the PEC scattering problem (cf.~\cite[Thm. 3.35]{kirsch2016mathematical}), we conclude that the fields $E^s, H^s \in C^{2}(\ED, \C^3) \cap C^{1,\alpha}(\R^3 \setminus \Omega, \C^3)$ are the unique solution pair and possess the desired regularity. This completes the proof.
\end{proof}

\section{Numerical examples\label{sec:numerics}}
In this section, we present numerical examples that demonstrate the solution of the proposed BIEs using a high-order Nyström method for globally smooth surfaces based on the Density Interpolation Method (DIM)~\cite{HDI3D,perez2019planewave,perez2020planewave,faria2021general}, which provides accurate treatment of both singular and nearly singular integrals. In particular, we employ the General Purpose variant of DIM~\cite{faria2021general}, as implemented in the open-source Julia package \texttt{Inti.jl}.

Meshing is performed entirely using Gmsh~\cite{geuzaine2009gmsh}, which supports high-order curved surface elements. In our experiments, all surface meshes consist of curved triangles defined by fourth-degree polynomials, offering sufficient geometric accuracy to reliably evaluate both the mean curvature ($\mathscr H$) and the full curvature operator $(\mathscr R)$ at the interior quadrature nodes employed by our Nystr\"om method. Although Gmsh does not enforce the global $C^{2,\alpha}$-smoothness of the surface $\Gamma$ assumed in our analysis, the resulting curved meshes exhibit sufficient regularity in practice to ensure convergence of the overall discretization scheme, as confirmed by the numerical examples presented in the sequel.

To accelerate the solution of the resulting linear systems, our Nyström method implementation in \texttt{Inti.jl} supports both $\mathcal{H}$-matrix compression and the FMM, via integration with the \texttt{HMatrices.jl} and \texttt{FMM3D.jl} packages, respectively. In all examples presented in this section, the linear systems are solved using GMRES~\cite{saad1986gmres}.

 \subsection{Field corrections to achieve exact divergence-free conditions}
 
Before delving into our numerical examples, we address a potential drawback of our BIE approach. A distinguishing feature of the numerical solutions obtained from the formulations presented in this work—unlike those produced by standard electromagnetic BIE methods (e.g., EFIE, MFIE, CFIE, etc.~\cite{volakis2012integral})—is that the electric and magnetic field solutions produced by the ECFOIE~\eqref{eq:ecfie} and MCFOIE~\eqref{eq:mcfie} (or their regularized versions, R-ECFOIE~\eqref{eq:R-ECFOIE} and R-MCFOIE~\eqref{eq:R-MCFOIE}) are not, in general, exact solutions to Maxwell’s equations in the exterior domain $\ED$. This discrepancy arises because the divergence-free condition for the scattered fields is not exactly imposed throughout the domain, but rather enforced only through the boundary conditions. As a result, the computed fields may exhibit nonzero divergence, typically on the order of the discretization error associated with the numerical BIE solution.

To address this issue, various correction procedures can be devised to recover fields that satisfy Maxwell’s equations exactly, while preserving the accuracy of the boundary conditions, which are only approximately satisfied due to discretization errors in the BIE solution. Among these strategies, perhaps the most straightforward is presented below.
(A different approach can, for example, be devised based on the Stratton--Chu formula, as discussed in Remark~\ref{rem:use_stratton_chu}.)

Let $\widetilde E^s$ denote the approximate scattered electric field, given by
\begin{equation}\label{eq:approx_scattered}
\widetilde E^s(x) = (\widetilde{\mathcal D} - \im\eta\widetilde{\mathcal S})\widetilde \varphi(x)+\xi\sum_{j=1}^J\nabla \Phi_j(x)\tilde\ell_j(\widetilde \varphi), \quad x \in \ED,
\end{equation}
where $\widetilde{\mathcal D}$, $\widetilde{\mathcal S}$, and $\tilde\ell_j$ denote, respectively, the double-layer potential~\eqref{eq:DL_pot}, the single-layer potential~\eqref{eq:SL_pot}, and the functional~\eqref{eq:functional}, approximated using the underlying surface-integral quadrature rule employed by our DIM-based Nystr\"om method. We recall that the surface $\Gamma$ is assumed to encompass $J$ maximal connected components. Here, $\widetilde\varphi$ represents the approximate numerical solution of~\eqref{eq:ecfie} or~\eqref{eq:R-ECFOIE} at the quadrature nodes. For $k \neq 0$, we then define the corrected electric scattered field as
\begin{equation}\label{eq:corrected_scattered}
\widetilde E^s_c(x):= \widetilde E^s(x) + \frac{1}{k^2}\nabla\dive \widetilde E^s(x), \quad x \in \ED.
\end{equation}

It is easy to verify, using the fact that $\widetilde E^s$ is an exact solution of the Helmholtz equation, that the corrected field satisfies Maxwell's equation
$\curl^2 \widetilde E^s_c - k^2 \widetilde E^s_c = 0$ in  $\ED$, so, in particular, $\curl \widetilde E_c^s=0$ in $\ED$. As we show in the numerical results presented below, being  $\dive \widetilde E^s$ of the order of the numerical discretization errors, the corrections term $k^{-2}\nabla\dive\widetilde E^s$ effectively enforces the divergence-free condition while not affecting the overall accuracy of the approximate electric field, that, our course, for suffciently large $k>0$ values. 

In a completely analogous manner, the magnetic field obtained by solving the MCFOIE~\eqref{eq:mcfie} or its regularized counterpart R-MCFOIE \eqref{eq:R-MCFOIE} takes the form
\begin{equation}\label{eq:approx_scattered_H}
\widetilde H^s(x) = (\widetilde{\mathcal D} - \im\eta,\widetilde{\mathcal S}),\widetilde \psi(x), \quad x \in \ED,
\end{equation}
where $\widetilde \psi$ denotes the approximate numerical solution of~\eqref{eq:mcfie} or~\eqref{eq:R-MCFOIE}.
This field can be corrected as
\begin{equation}\label{eq:corrected_scattered_H}
\widetilde H^s_c(x) := \widetilde H^s(x) + \frac{1}{k^2} \nabla \dive \widetilde H^s(x), \quad x \in \ED.
\end{equation}

\subsection{Error measures}
For the sake of definiteness, for the experiments presented in this section we measure numerical errors directly and indirectly using the following measures:
\begin{equation}\label{eq:error_measures}
e_{F} := \max_{j\in\{1,\ldots,100\}} \frac{|\widetilde F(x_j) - F_{\rm ref}(x_j)|}{|F_{\rm ref}(x_j)|}\quad\text{and} \quad 
e_{\dive F} := \max_{j\in\{1,\ldots,100\}} \frac{|\dive \widetilde F(x_j)|}{|\widetilde F(x_j)|},    
\end{equation}
where the reference field $F_{\rm ref}$ is an exact (electric or magnetic) solution to the problem under consideration, and $\widetilde F$ denotes an approximate solution evaluated using~\eqref{eq:approx_scattered} and/or~\eqref{eq:corrected_scattered}, as well as the corresponding formulae for the magnetic field. The target points $x_j$, $j\in\{1,\ldots,100\}$, are (approximately-) uniformly distributed on a sphere of radius 5 centered at the origin, ensuring the sphere encloses $\Gamma$ for all surfaces considered.  

Additionally, in some of our numerical experiments---particularly those at low frequencies---we report approximate values of the surface charge integrals
\begin{equation}\label{eq:surf_charge}
q_s := \sum_{j=1}^J q_s^{(j)},\qquad q_s^{(j)} := \int_{\Gamma_j} \nu \cdot \gamma E^s , \mathrm{d}s,
\end{equation}
which, under the assumptions on our incident fields, should theoretically vanish for all wavenumbers, including the static case $k = 0$~\cite{werner1963perfect_reflection}.

\subsection{Planewave scattering by a PEC sphere}
In our first set of examples, we consider the classical problem of electromagnetic wave scattering by a PEC unit sphere, $\Gamma = \mathbb{S}^2$, for which an exact solution can be obtained via the Mie series expansion~\cite{jackson2021classical}.  In detail, letting $p$ and $d$ respectively denote the polarization and propagation direction vectors, which are assumed to satisfy $d\cdot p =0$, the problem data used are given by 
\begin{subequations}\begin{equation}
        E^i(x):=p\e^{\im kx\cdot d}\text{ for } x\in\R^3;\qquad \p_\nu E^i(x) = \im k(d\cdot\nu(x))\gamma E^i(x)\text{ for } x\in\Gamma,
    \end{equation}
in the case of the electric field equations~\eqref{eq:ecfie} and~\eqref{eq:R-ECFOIE}, and by 
\begin{equation}
    H^i(x):=\sqrt\frac{\epsilon}{\mu}(d\times p)\e^{\im kx\cdot d}\text{ for } x\in\R^3;\qquad \p_\nu H^i(x) = \im k(d\cdot\nu(x))\gamma H^i(x)\text{ for } x\in\Gamma,
\end{equation}\label{eq:planewave}\end{subequations}
in the case of the magnetic field equations~\eqref{eq:mcfie} and~\eqref{eq:R-MCFOIE}. In this section we use the vectors $p = [1, 0, 0]$ and $d = [0, 0, 1]$,  and we set $\epsilon=\mu=1$.

Table~\ref{tab:planewave_sphere} reports the relative numerical errors defined in~\eqref{eq:error_measures}, computed with respect to the exact Mie series solution, which is available in Julia via the \texttt{SphericalScattering.jl} package~\cite{hofmann2023sphericalscattering}. Relative errors are presented for both the fields $\widetilde{E}^s$ and $\widetilde{H}^s$, as well as their divergence-free counterparts $\widetilde{E}^s_c$ and $\widetilde{H}^s_c$. The results correspond to a fixed wavenumber $k = \pi$ (i.e., wavelength $\lambda = 2$) and various nearly uniform surface discretizations characterized by the target mesh size $h$.

Besides demonstrating the validity of the proposed Maxwell BIE formulations and solvers , these experiments illustrate the preconditioning effect of the Calderón-type operator regularization in reducing the number of GMRES iterations by transforming the original ECFOIE and MCFOIE formulations into Fredholm second-kind BIEs. However, despite this reduction in iteration count, the results also suggest that regularization does not always yield a net performance benefit, as it can increase the computational cost per GMRES iteration.

This trade-off is reflected in the reported wall-clock times for solving the linear systems arising from the Nyström discretization of the respective BIEs, using both $\mathcal{H}$-matrix and FMM acceleration. In particular, while regularization improves GMRES convergence at this wavenumber, it may lead to longer overall solving times depending on the chosen acceleration strategy and problem size. Notably, some entries for the $\mathcal{H}$-matrix-based solver are missing due to lack of memory, showcasing a key advantage of the FMM, its lower memory footprint, which enabled it to successfully handle all problem sizes tested.

It is worth noting that all numerical experiments in this section were conducted on a machine equipped with an Apple M3 Max (14-core CPU) and 36 GB of RAM. No code performance optimization was pursued in these tests.

\begin{table}[htbp]
\centering
 \footnotesize
  \renewcommand{\arraystretch}{1.2}
\begin{tabular}{c|c|ccc|cc||ccc|cc}
\multicolumn{12}{c}{Single Sphere}\\
\toprule
\multicolumn{2}{c|}{} & \multicolumn{5}{c||}{ECFOIE~\eqref{eq:ecfie}}         & \multicolumn{5}{c}{R-ECFOIE~\eqref{eq:R-ECFOIE}}        \\ \hline
$\lambda/h$ & $N$  & $e_{E^s}$ & $e_{\dive E^s}$ & $e_{E^s_c}$ & iter. & time  (s)& $e_{E^s}$ & $e_{\dive E^s}$ & $e_{E^s_c}$ &iter. &time (s)\\ \hline
4           & 972  & 7.90e-3   & 1.92e-2        & 4.90e-4    & 108    &1$\mid$11 & 7.76e-3   & 2.00e-2          & 6.36e-3     & 32 &  2$\mid$10   \\
8           & 3228 & 9.20e-4   & 1.49e-3        & 9.40e-4     & 128    &6$\mid$91& 1.03e-3   & 1.70e-3        & 9.91e-4     & 28 &   4$\mid$70 \\
12          & 6792 & 3.19e-4   & 5.75e-4        & 2.44e-4     & 142    &16$\mid$139& 3.30e-4   & 6.58e-4        & 3.42e-4     & 27    &11$\mid$80\\
16          & 12576 & 1.04e-4  & 1.99e-4        & 8.53e-5     & 156    &36$\mid$284& 1.06e-4   & 2.09e-4        & 8.98e-5     & 26 &  23$\mid$148 \\
20          & 18912 & 5.38e-5  & 1.01e-4    & 4.87e-5     &  170 & 67$\mid$545& 6.55e-5   & 1.38e-4        & 5.73e-5     & 25 &  30$\mid$274\\
24          & 26700 & 3.16e-5  & 4.13e-5        & 3.23e-5     & 179    &99$\mid$597& 3.10e-5   & 4.74e-5        & 3.16e-5     & 24  & 52$\mid$280 \\ 
28          & 35652 & 2.03e-5  & 4.12e-5        & 1.55e-5     & 185    &153$\mid$802& 2.17e-5   & 4.31e-5        & 1.61e-5     & 24  &  164$\mid$342 \\ 
32          & 48180 & 9.10e-6  & 1.66e-5        & 9.23e-6     & 193    &233$\mid$1117& 9.28e-6   & 1.95e-5        & 8.52e-6     & 24  &| $\mid$484 \\ 
\toprule
\multicolumn{2}{c|}{} & \multicolumn{5}{c||}{MCFOIE~\eqref{eq:mcfie}}         & \multicolumn{5}{c}{R-MCFOIE~\eqref{eq:R-MCFOIE}}        \\ \hline
$\lambda/h$ & $N$  & $e_{H^s}$ & $e_{\dive H^s}$ & $e_{H^s_c}$ & iter.& time (s)& $e_{H^s}$ & $e_{\dive H^s}$ & $e_{H^s_c}$ & iter. &time (s)\\ \hline
4           & 972  & 4.34e-3   & 6.77e-3        & 3.94e-3     & 75   & 1$\mid$10& 5.64e-3   & 8.78e-3        & 5.11e-3     & 32   &  1$\mid$12 \\
8           & 3228 & 5.33e-4   & 1.03e-3        & 5.08e-4     & 92   & 4$\mid$65& 6.25e-4   & 1.67e-3        & 5.96e-4     & 30   & 5$\mid$73 \\
12          & 6792 & 1.62e-4   & 2.40e-4        & 1.57e-4     & 101    & 12$\mid$96& 1.69e-4   & 4.18e-4        & 1.59e-4     & 28   & 11$\mid$81  \\
16          & 12576 & 5.18e-5  & 9.60e-5        & 4.88e-5     & 112  &  26$\mid$196 & 5.09e-5   & 1.02e-4        & 4.72e-5     & 27  & 21$\mid$153 \\
20          & 18912 & 3.02e-5  & 5.43e-5        & 2.77e-5     & 131   & 47$\mid$416 & 3.13e-5   & 6.59e-5        & 3.09e-5     & 26   & 33$\mid$278  \\
24          & 26700 & 2.54e-5  & 4.52e-5        & 2.38e-5     & 144  & 79$\mid$485  & 2.41e-5   & 4.49e-5        & 2.36e-5     & 26  &  59$\mid$327\\ 
28& 35652 & 1.40e-5& 2.19e-5&1.35e-5&156& 123$\mid$648 &1.08e-5&2.44e-5&1.03e-5&26&175$\mid$366\\
32&48180&8.62e-6&1.32e-5&7.77e-6&163&176$\mid$923 &7.16e-6&1.15e-5&6.66e-6&25&| $\mid$488\\
\bottomrule
\end{tabular}
\caption{ 
Accuracy, GMRES iteration counts (with a relative tolerance of $10^{-8}$), and GMRES execution times for the numerical solution of the proposed BIE formulations for electromagnetic scattering by a PEC sphere $\Gamma$ of diameter $\lambda = 2$, illuminated by a plane wave~\eqref{eq:planewave} with wavenumber $k = \pi$. The results are computed for various mesh sizes $h > 0$, leading to surface discretizations involving $N$ quadrature nodes, which are used in the Nyström discretization of the BIEs and correspond to linear systems with $3N$ unknowns. For all examples, the parameters $\eta = k$ and $\xi = 0$ were used. GMRES execution times are  provided for implementations accelerated via both $\mathcal H$-matrix and FMM techniques (in that order). The R-ECFIE and R-MCFIE implementations employ the exact Calderón identity---i.e., it is based on the discretization of the operators $\oL$ in~\eqref{eq:compact_op} and $\oM$ in~\eqref{eq:compact_op_h}---with $\mathcal H$-matrix/FMM acceleration using a compression tolerance of $10^{-8}$.
}
\label{tab:planewave_sphere}
\end{table}

In our second numerical experiment, we evaluate the accuracy and performance of the proposed BIE formulations and their associated solvers in the low-frequency regime. Table~\ref{tab:low_frequency} reports the relative errors defined in~\eqref{eq:error_measures} for a fixed surface discretization and a range of wavenumbers, including values approaching zero within machine precision. 

As predicted by our analysis in Section~\ref{sec:magnetic_field_BIE}, the results confirm that the R-MCFOIE does not suffer from low-frequency breakdown in this setting, which aligns with the fact that the surface $\Gamma = \mathbb{S}^2$ is simply connected. In contrast, the R-ECFOIE exhibits a pronounced loss of accuracy when used with the parameter value $\xi = 0$, reflecting the manifestation of low-frequency breakdown---again, in agreement with the theoretical insights from Section~\ref{sec:electric_field_BIE}. This loss of accuracy is largely mitigated by choosing $\xi = 1$ in the combined-field ansatz~\eqref{eq:cfie_e_ansatz}. As expected, the nearly vanishing divergence values appear largely unaffected by the low-frequency breakdown, in contrast to the surface charge values $q_s$, which are more sensitive. Similar accuracy levels are observed for the electric and magnetic BIEs (ECFOIE and MCFOIE, respectively) when used without the single-layer regularization. However, in the ECFOIE case, achieving the prescribed relative error tolerance ($10^{-8}$) requires significantly more GMRES iterations, further highlighting the preconditioning effect of the Calderón-type operator regularization. We note that relative errors for the divergence-free corrected fields $\widetilde{E}^s_c$~\eqref{eq:corrected_scattered} and $\widetilde{H}^s_c$~\eqref{eq:corrected_scattered_H} are not reported in this experiment, as the correction term diverges in the limit $k \downarrow 0$, which amplifies numerical errors and renders the correction unreliable at extremely low frequencies.

\begin{figure}[hbt!]
  \centering
  \includegraphics[scale=0.5]{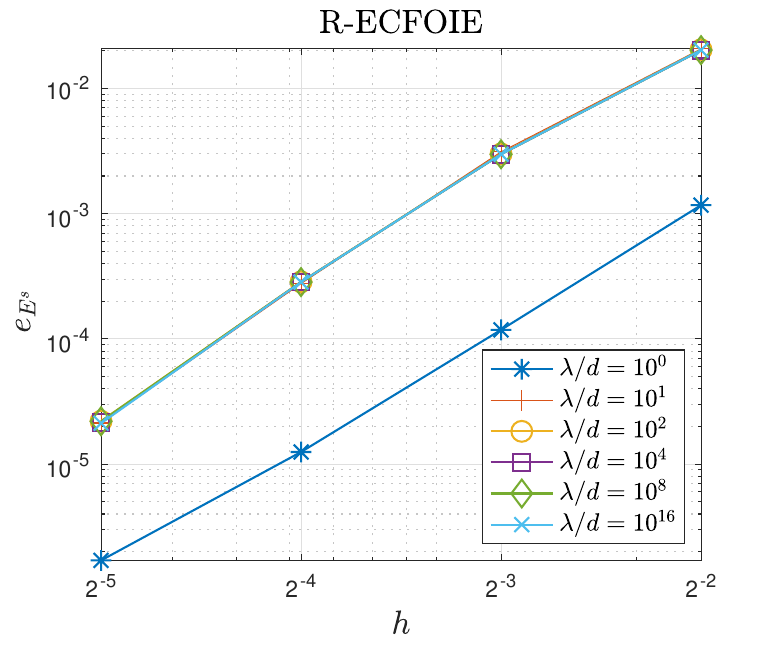}\qquad
  \includegraphics[scale=0.5]{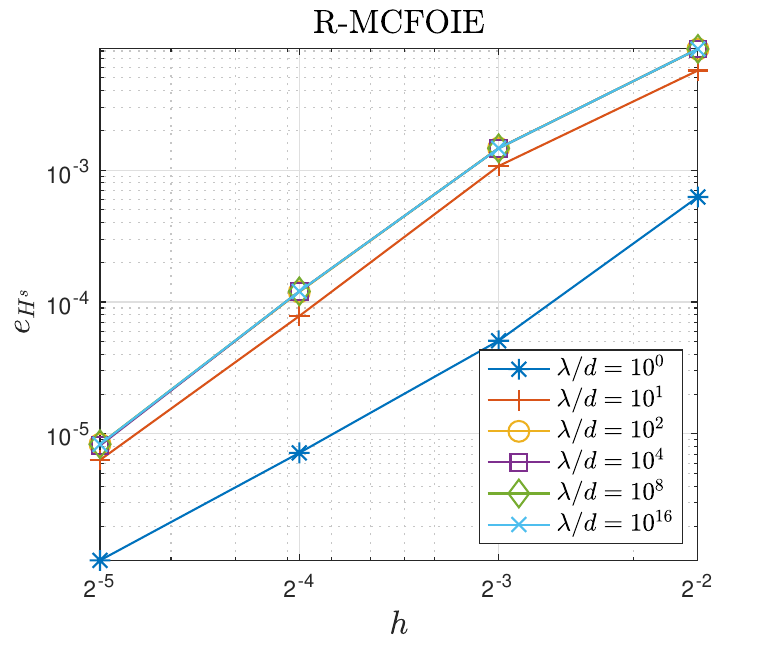}
  \caption{Error convergence plots at low frequencies. Left: Relative errors for various mesh sizes in the low-frequency regime, for the problem described in Table~\ref{tab:low_frequency}, using the R-ECFOIE~\eqref{eq:R-ECFOIE} with parameter values $\xi = 1$ and $\eta = \pi$. Right: Corresponding relative errors obtained with the R-MCFOIE~\eqref{eq:R-MCFOIE}. While convergence is evident, a degradation in accuracy is observed.}
  \label{fig:low_freq_convergence}
\end{figure}

Despite the significant accuracy improvements observed for the (R-)ECFOIE with $\xi = 1$, the errors are not uniformly small across all frequencies. In particular, a noticeable deterioration in accuracy is observed when $\lambda/d > 1$. To further investigate this undesired numerical behavior, Figure~\ref{fig:low_freq_convergence} examines whether the expected error convergence is recovered as the mesh size $h$ decreases, both in the low and extremely low frequency regimes. The results confirm that convergence is indeed achieved; however, substantial variations in accuracy persist across frequencies for both the R-ECFOIE and the R-MCFOIE. Since this behavior is also observed in the latter---which is not subject to low-frequency breakdown---we speculate that the accuracy variations may stem from non-uniformly accurate evaluations of the boundary integral operators across different frequencies. 

In our third and final experiment of this section, we examine the relative errors and GMRES iteration counts resulting from the discretization of the (R-)ECFOIE and (R-)MCFOIE formulations for high-frequency scattering problems. This includes a challenging case where the scatterer has a diameter of 32 wavelengths, leading to a linear system with over two million unknowns. As expected, the number of GMRES iterations increases with frequency. Interestingly, both regularized formulations exhibit comparable performance, and Calderón-type preconditioning appears to lose effectiveness at higher frequencies. In fact, the iteration count for the unregularized formulations grows more slowly than for the regularized ones in this regime. It is also worth noting that, across all examples, the divergence errors remain on par with the relative errors, indicating that the approximate divergence-free condition is consistently preserved.

\begin{table}[h!]
\centering
 \footnotesize
 \renewcommand{\arraystretch}{1.2}
\begin{tabular}{c|ccc|c|ccc|c||cc|c}
\multicolumn{12}{c}{Single Sphere}\\
\toprule
 & \multicolumn{4}{c|}{R-ECFOIE~\eqref{eq:R-ECFOIE}; $\xi = 0$}         & \multicolumn{4}{c||}{R-ECFOIE~\eqref{eq:R-ECFOIE}; $\xi=1$}     &\multicolumn{3}{c}{R-MCFOIE~\eqref{eq:R-MCFOIE}}    \\ \hline
$\lambda/d$    & $e_{E^s}$ & $e_{\dive\! E^s}$ &$|q_s|$&  \#iter & $e_{E^s}$ & $e_{\dive\! E^s}$ &$|q_s|$ &  \#iter& $e_{H^s}$ & $e_{\dive\! H^s}$ &   \#iter \\ \hline
$10^{16}$   &  {\bf 1.30e1}  &    2.99e-4        & {\bf 3.24e1} &  25  & 2.98e-3   & 3.43e-4             &6.30e-15& 25 &1.46e-3&1.48e-4&13    \\
$10^{08}$          &   {\bf 1.30e1} &     2.99e-4        & {\bf 3.24e1} &  25   &     2.98e-3       &   3.43e-4   &2.62e-14&   25   &1.46e-3&1.48e-4&13\\
$10^{04}$          &  {\bf  1.30e1} &      3.00e-4        &{\bf 3.26e1}& 25    &     2.98e-3       &   3.43e-4   &2.03e-14&   25   &1.46e-3&1.48e-4&13\\
$10^{02}$         &  {\bf 1.52e-1} &    8.80e-4         & {\bf 3.73e-1}&22     &   3.00e-3        &   3.45e-4   &1.23e-14&25    &1.47e-3&1.47e-4&14  \\
$10^{01}$         & 4.75e-3  &    8.92e-4    &4.07e-3& 20    &      3.11e-3      &   3.35e-4   &9.63e-15& 26   &1.07e-3&8.02e-5&16  \\
$10^{00}$         & 1.32e-4  & 2.50e-4         &5.52e-5   & 26    & 1.18e-4           & 2.24e-4     &4.87e-15& 38  &5.09e-5&1.02e-4&26   \\ 
\toprule
 & \multicolumn{4}{c|}{ECFOIE~\eqref{eq:ecfie}; $\xi = 0$}         & \multicolumn{4}{c||}{ECFOIE~\eqref{eq:ecfie}; $\xi=1$}     &\multicolumn{3}{c}{MCFOIE~\eqref{eq:mcfie}}    \\ \hline
$\lambda/d$    & $e_{E^s}$ & $e_{\dive\! E^s}$&$|q_s|$&  \#iter & $e_{E^s}$ & $e_{\dive\! E^s}$ & $|q_s|$&   \#iter& $e_{H^s}$ & $e_{\dive\! H^s}$ &  \#iter \\ \hline
$10^{16}$   &  {\bf 2.57e1}  &    7.82e-4        & {\bf 6.46e1} &   287  &       1.69e-3     &    2.43e-4 &4.71e-14&   270  &6.31e-4&8.36e-5&82\\
$10^{08}$   & {\bf 2.57e1}   &     7.82e-4        &{\bf 6.46e1} &    287  & 1.69e-3   &   2.43e-4           &4.55e-14& 270 &6.31e-4&8.36e-5& 82   \\
$10^{04}$    & {\bf 2.53e1}  &   7.80e-4          & {\bf 6.35e1}& 287    &       1.69e-3     &   2.43e-4   &1.63e-14&   270   &6.31e-4&8.36e-5&85\\
$10^{02}$         &{\bf 1.43e-1}   &   7.95e-4          & {\bf 3.52e-1}& 224     &        1.69e-3   &2.45e-4      &2.45e-14&  278  &6.29e-4&8.26e-5&89  \\
$10^{01}$         &  3.56e-3 &    9.11e-4    &3.86e-3 &   203 &  1.69e-3          &2.68e-4      &5.34e-14& 285   &4.25e-4&6.18e-5& 99 \\
$10^{00}$         &1.33e-4   &   3.05e-4          &6.42e-5& 163    &    1.17e-4        & 2.38e-4     &8.65e-15&249   &5.18e-5&9.60e-5&  111 \\ 
\bottomrule
\end{tabular}
\caption{Relative errors, induced surface charge~\eqref{eq:surf_charge}, and GMRES iteration counts (with an absolute tolerance of $10^{-8}$) for all four BIE formulations applied to the electromagnetic scattering problem by the unit sphere $\Gamma = \mathbb{S}^2$ (of diameter $d=2$) under planewave illumination~\eqref{eq:planewave} at low frequencies. The results are obtained using a fixed mesh size of $h = 0.125$. In all cases, the parameter $\eta = \pi$ is used, and two different values of the parameter $\xi$---introduced to mitigate the low-frequency breakdown---are shown, demonstrating its effectiveness in the case of the (R-)ECFOIE.}
\label{tab:low_frequency}
\end{table}

\begin{table}[h!]
\centering
 \footnotesize
 \renewcommand{\arraystretch}{1.2}

\begin{tabular}{c|c|ccc|ccc}
\multicolumn{8}{c}{Single Sphere}\\
\toprule
 && \multicolumn{3}{c|}{ECFOIE~\eqref{eq:ecfie}}         & \multicolumn{3}{c}{R-ECFOIE~\eqref{eq:R-ECFOIE}}    \\\hline
$\lambda/d$   & $N$& $e_{E^s}$ &   \#iter &time (s)& $e_{E^s}$ &   \#iter&time (s)\\ \hline
2   & 3228   &    9.59e-3         &   96   &  78  &    1.24e-2          & 47 &125 \\
4   &12576   &    1.06e-3         &   92 &  185  &  9.08e-4     & 59      &370\\
8   & 48180  &    5.26e-4        &  103   &   745&   6.45e-4    & 79     &1936\\
16   &  185796 &    4.94e-4         &  111   &     3812      & 6.76e-4     &111    &12885\\
32         &  738612 &   9.73e-4    &  121   &     21110       & |    & |  & |\\\bottomrule
 &&\multicolumn{3}{c|}{MCFOIE~\eqref{eq:mcfie}}&\multicolumn{3}{c}{R-MCFOIE~\eqref{eq:R-MCFOIE}}    \\ \hline
 $\lambda/d$   & $N$&  $e_{H^s}$  &   \#iter&time (s)& $e_{H^s}$  &   \#iter &time (s) \\ \hline
2& 3228&1.33e-2&73&60&   1.59e-2&43&116\\
4&12576&9.75e-4&72&144&8.14e-4&59&376\\
8&48180&5.59e-4&88&663&8.44e-4&80&2024\\
16&185796&5.97e-4&98 &3437&7.07e-4&112&12671 \\
32&738612&1.51e-3&110&18447 &|&|&| \\
\bottomrule
\end{tabular}
\caption{Relative errors and GMRES iteration counts (with an absolute tolerance of $10^{-8}$) are reported for the regularized BIE formulations applied to the electromagnetic scattering problem by the unit PEC sphere $\Gamma = \mathbb{S}^2$ (with diameter $d = 2$) under plane wave illumination~\eqref{eq:planewave} at higher frequencies. The results are obtained using mesh sizes of defined by $h = \lambda/4$. In all cases, the parameters $\eta = k$ and $\xi = 0$ are used. The linear systems are solved using GMRES accelerated by the FMM with a tolerance of $10^{-8}$.}
\label{tab:high_frequency}
\end{table}

Finally, we mention that the Julia code that can (partially) generate the results presented in this section is publicly available at \url{https://github.com/caperezar/Maxwell-a-la-Helmholtz}.

\subsection{Dipole fields as manufactured solutions}
In our second set of examples, we consider electromagnetic fields constructed as linear combinations of electric and magnetic dipole fields. In detail, the electric dipole fields are given by:
\begin{subequations}\begin{equation}\label{eq:dipole_field_electric}\begin{aligned}
H_e(x) &:= \curl_x\sum_{j=1}^Jp_jG(x,y_j) = \sum_{j=1}^J\nabla_x G(x,y_j) \times p_j, \qquad
E_e(x):= -\frac{1}{\im\omega\epsilon} \curl H_e(x),
\end{aligned}\end{equation}
while the magnetic dipole fields are defined as
\begin{equation}\label{eq:dipole_field_magnetic}\begin{aligned}
E_m(x) &:= \curl_x\sum_{j=1}^Jp_jG(x,y_j) =\sum_{j=1}^J \nabla_x G(x,y_j) \times p_j, \qquad
H_m(x) := \frac{1}{\im\omega\mu} \curl E_m(x),
\end{aligned}\end{equation}\label{eq:dipole_fields}\end{subequations}
for $x\in\R^3\setminus\bigcup_{j=1}^J\{y_j\}$,  where $p_j \in \mathbb{R}^3$ and $y_j \in \Omega_j$, for $j=1,\ldots,J$, with $\Omega_j$ denoting the maximal connected components of $\Omega$.

Clearly, both fields represent exact solutions to Maxwell’s equations~\eqref{eq:maxwell} in the exterior domain $\ED$, and they satisfy the Silver--Müller radiation condition~\eqref{eq:rad_cond_sm}.
Therefore, by Theorems~\ref{thm:elf} and~\ref{lem:equiv_int}, they provide exact solutions $E^s$ and $H^s$ to the exterior vector Helmholtz problems~\eqref{eq:elf_equiv} and~\eqref{eq:mgf_equiv} when the incident fields are given by $E^i = -E_e$ and $H^i = -H_m$, respectively. Moreover, in the limit $k \downarrow 0$, the fields $E_e$ and $H_m$ reduce to the exact solutions $E_0$ and $H_0$ of the static problems~\eqref{eq:vector_laplace} and~\eqref{eq:H_vec}, respectively. In this section, we use these fields as reference solutions to validate our BIE approach on more complex geometries, including disconnected and multiply connected surfaces $\Gamma$. 

Analogous to Table~\ref{tab:low_frequency}, Table~\ref{tab:low_frequency_2_spheres} presents numerical errors for low-frequency problems involving a disconnected surface $\Gamma$ composed of two disjoint unit spheres, $\Omega_1$ and $\Omega_2$, centered at $c_1$ and $c_2$, respectively. In this case, the surfaces and the dipole fields~\eqref{eq:dipole_field_electric} and~\eqref{eq:dipole_field_magnetic} are constructed using the following numerical values:
\begin{equation}\label{eq:param_two_spheres}
\begin{aligned}
&c_1 = [-1.5, 0, 0], \quad &&y_1 = c_1 + [0, 0.1, -0.2], \quad &&p_1 = \frac{1}{\sqrt 3}[1, 1, 1], \\
&c_2 = [1.5, 0, 0], \quad &&y_2 = c_2 + [0, -0.2, 0.1], \quad &&p_2 = \frac{1}{\sqrt 3}[-1, 1, -1].
\end{aligned}
\end{equation}

\begin{table}[h!]
\centering
\footnotesize
 \renewcommand{\arraystretch}{1.2}
\begin{tabular}{c|ccc|c|ccc|c||cc|c}
\multicolumn{12}{c}{Two Spheres}\\
\toprule
 & \multicolumn{4}{c|}{R-ECFOIE~\eqref{eq:R-ECFOIE}; $\xi = 0$}         & \multicolumn{4}{c||}{R-ECFOIE~\eqref{eq:R-ECFOIE}; $\xi=1$}     &\multicolumn{3}{c}{R-MCFOIE~\eqref{eq:R-MCFOIE}}    \\ \hline
$\lambda/d$    & $e_{E^s}$ & $e_{\dive\! E^s}$ &$\displaystyle\max_{j=1,2}|q^{(j)}_s|$&  \#iter & $e_{E^s}$ & $e_{\dive\! E^s}$ &$\displaystyle\max_{j=1,2}|q^{(j)}_s|$ &  \#iter& $e_{H^s}$ & $e_{\dive\! H^s}$ &   \#iter \\ \hline
$10^{16}$   &  {\bf 4.72e0}  &   2.73e-4         & {\bf 4.69e-1} & 55    & 2.53e-3   & 2.69e-4             &3.96e-10& 54 &1.26e-3&1.36e-4&27    \\
$10^{08}$          &   {\bf 4.72e0} &    2.73e-4        & {\bf 4.69e-1} &   55  &  2.53e-3          &    2.69e-4  &3.96e-10&   54   &1.26e-3&1.36e-4&27\\
$10^{04}$          &  {\bf  4.67e0} &              2.75e-4& {\bf 4.65e-1} &55   &     2.53e-3       &   2.69e-4   &2.28e-8&      54&1.26e-3&1.36e-4&27\\
$10^{02}$         &  {\bf 4.16e-2} &    5.02e-4         & {\bf 5.09e-3}&  46   &   2.39e-3        &    2.82e-4  &2.16e-4&  54  &1.22e-3&1.32e-4&27  \\
$10^{01}$         & 1.49e-3  &   3.85e-4    &1.32e-4&  39   &      1.48e-3      &   3.84e-4   &1.32e-4& 46   &9.42e-4&1.78e-4&27  \\
$10^{00}$         &  2.03e-5 &   5.63e-5       &  9.73e-6 &  47   &    2.03e-5        &   5.63e-5   &7.31e-6& 48  &1.68e-5&4.00e-5&41   \\ 
\toprule
 & \multicolumn{4}{c|}{ECFOIE~\eqref{eq:ecfie}; $\xi = 0$}         & \multicolumn{4}{c||}{ECFOIE~\eqref{eq:ecfie}; $\xi=1$}     &\multicolumn{3}{c}{MCFOIE~\eqref{eq:mcfie}}    \\ \hline
$\lambda/d$    & $e_{E^s}$ & $e_{\dive\! E^s}$&$\displaystyle\max_{j=1,2}|q^{(j)}_s|$&  \#iter & $e_{E^s}$ & $e_{\dive\! E^s}$ & $\displaystyle\max_{j=1,2}|q^{(j)}_s|$&   \#iter& $e_{H^s}$ & $e_{\dive\! H^s}$ &  \#iter \\ \hline
$10^{16}$   &  {\bf 1.99e0}  & 1.51e-4           & {\bf 2.64e-1} & 616     &1.85e-3            & 1.33e-4    &7.56e-10&     631&9.86e-4&1.00e-4&184\\
$10^{08}$   & {\bf 1.99e0}   & 1.51e-4           &{\bf 2.64e-1} &  616    &1.85e-3    & 1.33e-4            &7.56e-10& 631 &9.86e-4&1.00e-4&  184  \\
$10^{04}$    & {\bf  1.98e0}  &   1.51e-4          & {\bf 2.62e-1 }&  616   &  1.85e-3          & 1.33e-4     &3.97e-8&   631   &9.86e-4&1.00e-4&184\\
$10^{02}$         &{\bf 3.98e-2}   &  2.63e-4           & {\bf 5.61e-3}& 551     & 1.93e-3          &1.34e-4      &3.52e-4& 630  &9.56e-4&9.94e-5&184 \\
$10^{01}$         &  9.62e-4 & 2.73e-4       &1.57e-4 & 436  &  9.76e-4          &2.22e-4      &1.83e-4& 535   &7.62e-4&1.40e-4&188 \\
$10^{00}$         &1.81e-5   &     5.00e-5       &8.91e-6&274& 1.79e-5               &5.07e-5      &6.73e-6& 278  &1.09e-5&2.32e-5&153  \\ 
\bottomrule
\end{tabular}
\caption{Relative errors, approximate induced surface charge~\eqref{eq:surf_charge}, and GMRES iteration counts (with an absolute tolerance of $10^{-8}$) reported for the BIE formulations applied to the PEC electromagnetic boundary value problem, with boundary data defined by the electric and magnetic dipole fields in equations~\eqref{eq:dipole_field_electric} and~\eqref{eq:dipole_field_magnetic}, respectively, using the numerical parameters specified in~\eqref{eq:param_two_spheres}. The surface $\Gamma$ is disconnected consisting of two disjoint spheres of equal diameter $d=2$. The mesh size $h=0.125$ and the parameter $\eta = \pi$ are used throughout. The results for $\xi = 0$ and $\xi = 1$ illustrate the effectiveness of the modified combined-field ansatz~\eqref{eq:cfie_e_ansatz} in mitigating the low-frequency breakdown in the (R-)ECFOIE formulation.}
\label{tab:low_frequency_2_spheres}
\end{table}

These results demonstrate the effectiveness of the combined-field ansatz~\eqref{eq:cfie_e_ansatz} with $\xi=1$ in mitigating the effects of low-frequency breakdown in the (R-)ECFOIE formulation in the case disconnected surfaces, while also providing numerical evidence of the frequency robustness of the (R-)MCFOIE. It is worth noting that the invertibility condition for the matrix $\Id_2 + \xi\Xi$, which is required to ensure the uniqueness of solutions to the (R-)ECFOIE (see Lemma~\ref{lem:unique}), was verified in all the examples reported in this table. In fact, the condition number of this matrix remained below $2$ in every case.

\begin{figure}[hbt!]
   \centering
   \includegraphics[width=0.35\textwidth]{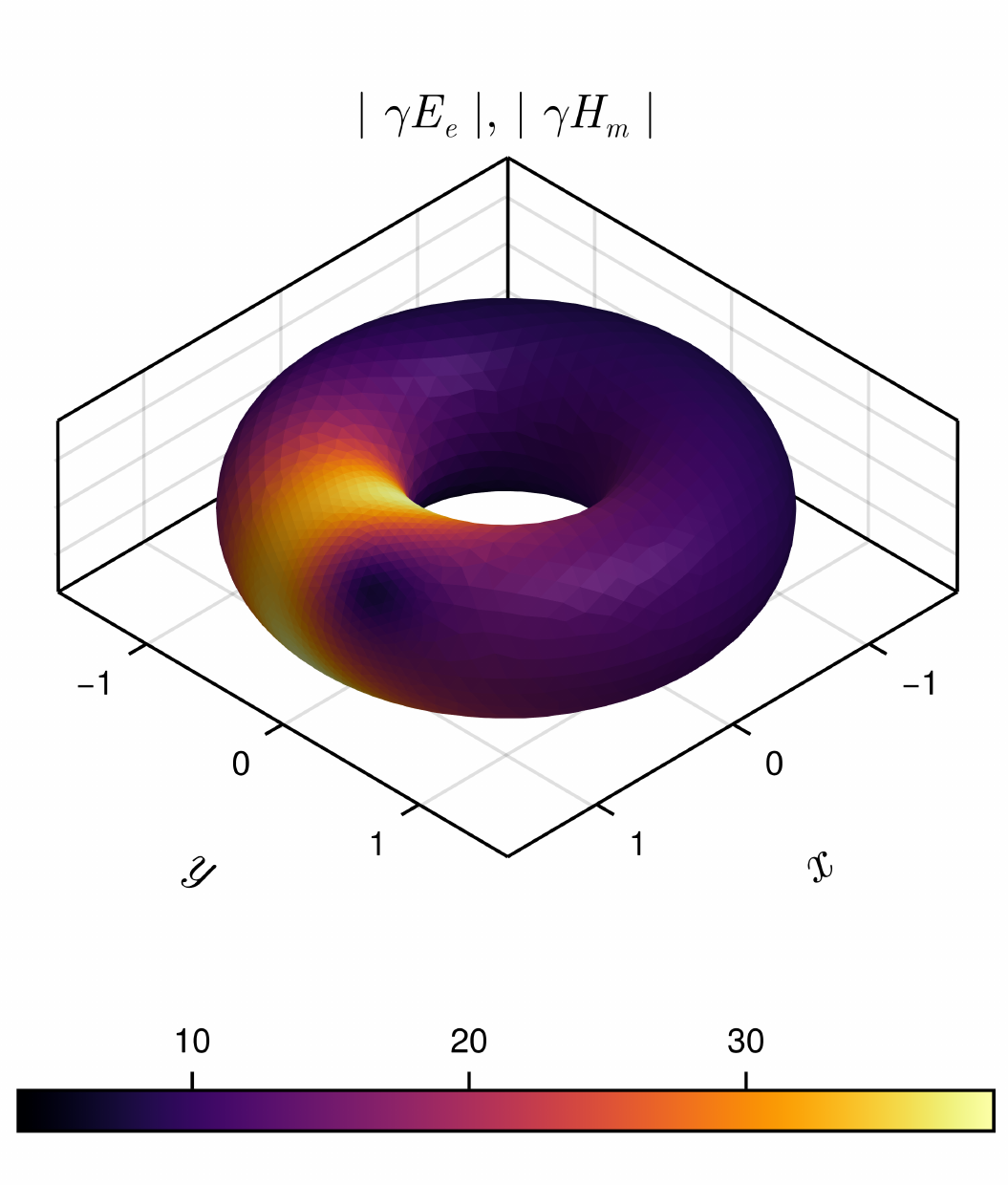}\hspace{2cm}
   \includegraphics[width=0.35\textwidth]{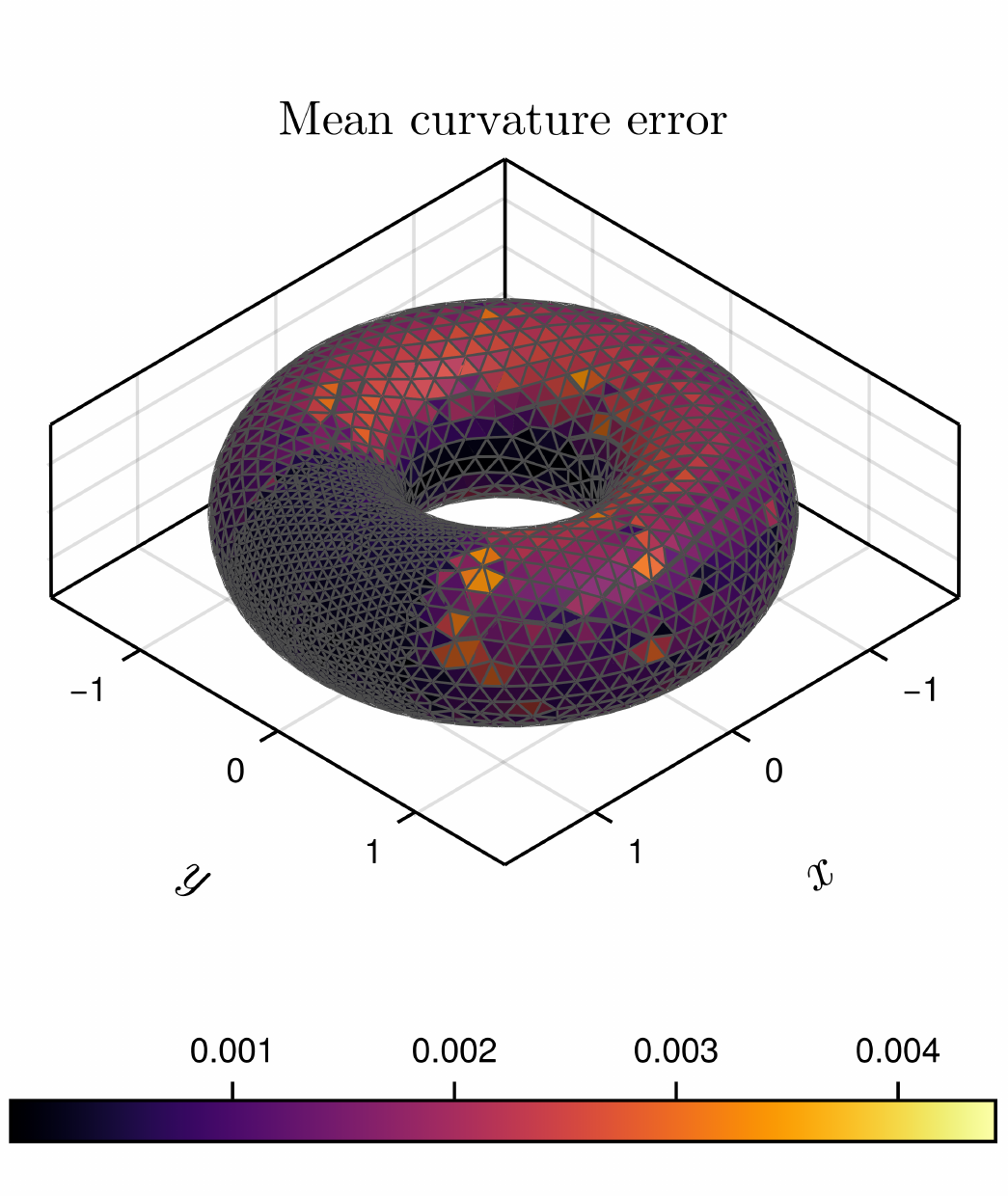}
   \caption{Left: Magnitude of the surface trace of dipole fields used in the numerical experiments reported in Table~\ref{tab:low_frequency_torus}.
Right: Surface mesh used in the same experiments, colored according to the absolute error in the mean curvature computed using \texttt{Inti.jl}.}
   \label{fig:torus_geometry}
 \end{figure}

Similarly, Table~\ref{tab:low_frequency_torus} presents the numerical errors at low frequencies for the case of a single torus surface with minor and major radii of $1/2$ and $1$, respectively, centered at the origin. The results show that our approach remains valid in the context of a multiply connected geometry. In this experiment, we use dipole fields~\eqref{eq:dipole_fields} with
\begin{equation}\label{eq:data_torus}
y_1 = [1, 0, 0] \quad \text{and} \quad p_1 = \frac{1}{\sqrt3}[1, 1, 1].
\end{equation}
The left panel of Figure~\ref{fig:torus_geometry} shows the magnitude of the dipole fields on the torus surface. Due to the large field gradients, the mesh used in these examples, shown in the right panel of Figure~\ref{fig:torus_geometry}, is adaptively refined based on both the proximity to the dipole source and the local mean curvature. In this panel, the surface is colored according to the absolute error in the mean curvature computed using \texttt{Inti.jl}.
These results highlight the capability of \texttt{Inti.jl} to generate sufficiently smooth surface representations, thereby enabling the accurate geometric computations required for the robust implementation of the proposed BIEs. (For completeness, Appendix~\ref{app:curv_tens} provides an explicit representation of the curvature operator $\mathscr{R}$ in terms of the coefficients of the first and second fundamental forms and a basis of the tangent plane to the surface $\Gamma$. These formulae were used in all geometric computations throughout the paper.)

Interestingly, and in contrast with theoretical expectations (see Remark~\ref{rem:uniqueness_magnetostatic}), in practice, the magnetic BIEs do not exhibit signs of low-frequency breakdown at extremely low frequencies in the non simply connected setting. In fact, reasonably small values of $\dive H^s$ and $\curl H^s$---consistent with the discretization---are also obtained for both plane wave incidence and for incident fields generated by dipoles distributed along a circle in the $yz$-plane, enclosing one of the cross-sectional circles of the torus. These positive results, however, may well be due to the fact that the incident fields considered do not excite components lying in the nearly non-trivial kernel of the resulting linear systems. Additional contour integral conditions---similar to those used in~\cite{epstein2010debye}---may need to be incorporated to fully address this potential issue. Further investigation into this particular aspect of the (R-)MCFOIE will be the subject of future work.

\begin{table}[h!]
\centering
\footnotesize
 \renewcommand{\arraystretch}{1.2}
\begin{tabular}{c|ccc|c|ccc|c||cc|c}
\multicolumn{12}{c}{Torus}\\
\toprule
 & \multicolumn{4}{c|}{R-ECFOIE~\eqref{eq:R-ECFOIE}; $\xi = 0$}         & \multicolumn{4}{c||}{R-ECFOIE~\eqref{eq:R-ECFOIE}; $\xi=1$}     &\multicolumn{3}{c}{R-MCFOIE~\eqref{eq:R-MCFOIE}}    \\ \hline
$\lambda/d$    & $e_{E^s}$ & $e_{\dive\! E^s}$ &$|q_s|$&  \#iter & $e_{E^s}$ & $e_{\dive\! E^s}$ &$|q_s|$ &  \#iter& $e_{H^s}$ & $e_{\dive\! H^s}$ &   \#iter \\ \hline
$10^{16}$          &   {\bf 8.47e-2} &    4.95e-4        & {\bf 1.04e-2} & 47    & 5.18e-3   & 6.28e-4  &  6.54e-16 & 51 & 5.68e-3 & 6.28e-4 & 44 \\
$10^{08}$          &   {\bf 8.47e-2} &    4.95e-4        & {\bf 1.04e-2} & 47 &  5.18e-3         &    6.28e-4  & 5.20e-16 &   51   & 5.68e-3 & 6.28e-4 & 44 \\
$10^{04}$          &  {\bf  8.47e-2} &  4.95e-4 & {\bf 1.04e-2} &  47  &     5.18e-3  &   6.28e-4   & 4.99e-16 &  51  & 5.68e-3 & 6.28e-4 & 44 \\
$10^{02}$         &  {\bf 5.74e-2} &    5.50e-4         & {\bf 6.85e-3}&  46   &   5.17e-3        &    6.40e-4  &  6.20e-16  &  51  & 5.53e-3 & 6.39e-4 & 44 \\
$10^{01}$         & 2.23e-3  &   5.65e-4    & 2.07e-4 &  43  &      1.77e-3      &   4.31e-4   & 4.93e-16 & 52   & 1.66e-3 & 5.29e-4 & 41 \\
$10^{00}$         &  2.12e-4 &  1.49e-4  & 1.25e-5  &  56   &    2.24e-4        &   1.37e-4   & 3.32e-15 & 72  & 1.42e-4 & 3.30e-4 & 48 \\ 
\toprule
 & \multicolumn{4}{c|}{ECFOIE~\eqref{eq:ecfie}; $\xi = 0$}         & \multicolumn{4}{c||}{ECFOIE~\eqref{eq:ecfie}; $\xi=1$}     &\multicolumn{3}{c}{MCFOIE~\eqref{eq:mcfie}}    \\ \hline
$\lambda/d$    & $e_{E^s}$ & $e_{\dive\! E^s}$&$|q_s|$&  \#iter & $e_{E^s}$ & $e_{\dive\! E^s}$ & $|q_s|$&   \#iter& $e_{H^s}$ & $e_{\dive\! H^s}$ &  \#iter \\ \hline
$10^{16}$   & {\bf 5.68e-2} & 5.62e-4 &{\bf 6.88e-3} & 775 & 5.45e-3 & 5.71e-4 & 2.09e-16 & 843& 7.21e-3 & 2.34e-4 & 435 \\
$10^{08}$   & {\bf 5.68e-2} & 5.62e-4 &{\bf 6.88e-3} & 775 & 5.45e-3 & 5.71e-4 & 4.27e-16 & 843 & 7.21e-3 & 2.34e-4 & 435 \\
$10^{04}$    & {\bf 5.68e-2}  & 5.63e-4 & {\bf 6.88e-3} &  775 &  5.45e-3 & 5.71e-4 & 1.15e-16 & 843 & 7.22e-3 & 2.34e-4 & 435\\
$10^{02}$         &{\bf 3.87e-2}   &  6.14e-4 & {\bf 4.49e-3}& 764 & 5.49e-3 & 5.82e-4 & 8.73e-16 & 843 & 6.61e-3 & 2.39e-4 & 434 \\
$10^{01}$ & 2.26e-3 & 6.17e-4 & 1.53e-4 & 702 &  2.19e-3 & 5.07e-4 & 2.55e-16 & 851 & 7.45e-4 & 1.67e-4 & 404 \\
$10^{00}$         & 9.85e-5   &  2.11e-4  & 5.42e-5 & 506 & 1.52e-4 & 3.85e-4 & 3.88e-15 & 743  & 6.55e-5 & 1.69e-4 & 276  \\ 
\bottomrule
\end{tabular}
\caption{Relative errors, approximate induced surface charge~\eqref{eq:surf_charge}, and GMRES iteration counts (with an absolute tolerance of $10^{-8}$) reported for the BIE formulations applied to the PEC electromagnetic boundary value problem, with boundary data defined by the electric and magnetic dipole fields in equations~\eqref{eq:dipole_field_electric} and~\eqref{eq:dipole_field_magnetic}, respectively, using the numerical parameters specified in~\eqref{eq:data_torus}. The surface $\Gamma$ consists of a torus with radii $1$ and $1/2$, which is displayed in Figure~\ref{fig:torus_geometry}. The mesh has been locally refined near the dipole location with minimum and maximum mesh sizes of 1/16 and 1/8, respectively. The parameter $\eta = \pi$ is used throughout. The results for $\xi = 0$ and $\xi = 1$ illustrate the effectiveness of the modified combined-field ansatz~\eqref{eq:cfie_e_ansatz} in mitigating the low-frequency breakdown in the (R-)ECFOIE formulation.}
\label{tab:low_frequency_torus}
\end{table}


\subsection{Scattering by more general smooth surfaces}

In this final section, we consider examples involving more complex, smooth surfaces $\Gamma$. We focus on PEC scattering problems under plane wave incidence~\eqref{eq:planewave}. Since exact solutions are not available for these geometries and incidences, we once again assess the accuracy of our numerical results by also solving the same problems using manufactured boundary data derived from dipole fields placed inside the surface, following the approach described in the previous section.

In particular, we focus on two smooth geometries: the bean-shaped surface~\cite{Bruno:2001ima} and the flower-shaped surface~\cite{wildman2004accurate}. These surfaces are respectively parametrized by 
\begin{equation}\begin{aligned}
\text{bean:}\quad&  \left[\frac{4}{5}\sqrt{1-\frac{1}{10}\cos(\pi \cos(\theta))}\cos(\phi)\cos(\theta),-\frac{3}{10}\cos(\pi \cos(\theta))+\right.\\
&\hspace{7cm}\left[\frac{4}{5}\sqrt{1-\frac{2}{5}\cos(\pi \cos\theta)}\sin(\phi)\sin(\theta),\cos(\theta)\right],\\
\text{flower:}\quad&  
\sqrt{0.8+0.5(\cos(2\phi)-1)(\cos(4\theta)-1)}\left[\cos(\phi)\sin(\theta),\sin(\phi)\sin(\theta),\cos(\theta)\right],
\end{aligned}
\label{eq:smooth_surfaces}\end{equation}
for $(\theta,\phi)\in[0,\pi]\times[0,2\pi] $.

For each of these surfaces and for both the electric and magnetic regularized BIE formulations, we solve two scattering problems corresponding to the wavenumbers $k = \pi \cdot 10^{-8}$ and $k = 5\pi$. Accuracy in each case is assessed using the error measures $e_{E^s}$ and $e_{H^s}$, defined as in~\eqref{eq:error_measures}, with the dipole fields $E_e$ and $H_m$~\eqref{eq:dipole_fields} serving both as reference solutions and as sources of the manufactured boundary data; as in the above section. The dipoles are placed at $ y_1 = [0,-0.25,0]$ in the case of the bean-shaped surface, and at $y_1=[0,0.1,-0.2]$ in the case of the flower-shaped surface, and both with $p_1=\frac{1}{\sqrt3}[1,1,1]$. The results of these test problems are reported in Table~\ref{tab:general_surfaces}. In all cases, we used GMRES with an absolute tolerance of $10^{-4}$ and $\mathcal{H}$-matrix acceleration. All surface meshes were generated using Gmsh with a target mesh size of $h = 0.07$.

\begin{table}[h!]
\centering
\small 
 \renewcommand{\arraystretch}{1.2}
\begin{tabular}{c|cc|cc||cc|cc}
 & \multicolumn{4}{c||}{Bean}         & \multicolumn{4}{c}{Flower}         \\ \hline
$k$    & $e_{E^s}$ &  \#iter & $e_{H^s}$ & \#iter& $e_{E^s}$ &    \#iter&$e_{H^s}$ &    \#iter\\ \hline
$10^{-8}\pi$&7.32e-3&42&6.97e-3&23&8.54e-3&45&8.27e-3&27\\ 
$5\pi$&1.88e-3&122&3.41e-3&108&2.12e-3&144&3.21e-3&141\\
\bottomrule
\end{tabular}
\caption{Accuracy and GMRES iteration counts for the dipole field test problems involving the bean- and flower-shaped surfaces defined in~\eqref{eq:smooth_surfaces} and displayed in Figures~\ref{fig:electric_scattering_general_surfaces} and~\ref{fig:magnetic_scattering_general_surfaces}.}
\label{tab:general_surfaces}
\end{table}

\begin{figure}[hbt!]
   \centering
   \includegraphics[width=0.78\textwidth]{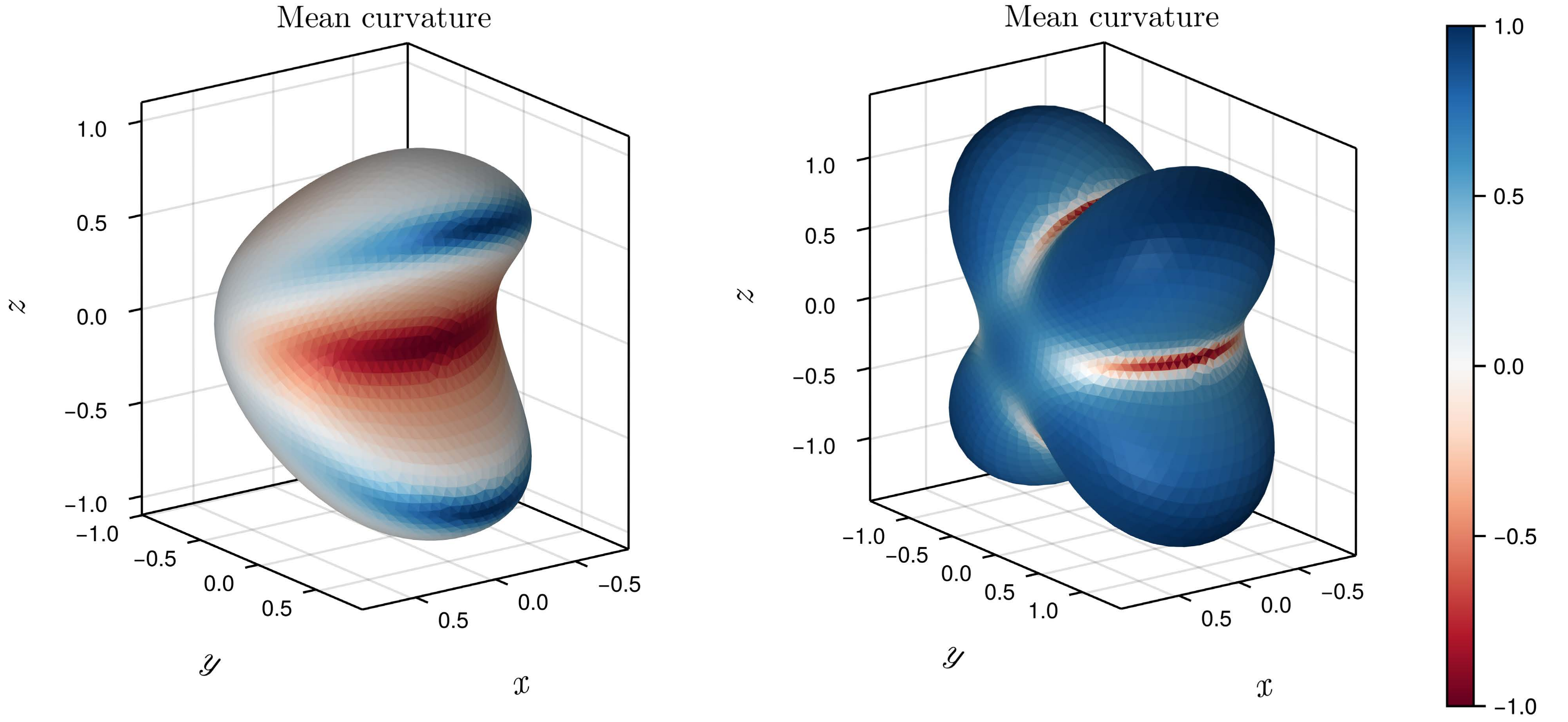}\\[0.2mm]  
    \includegraphics[width=0.78\textwidth]{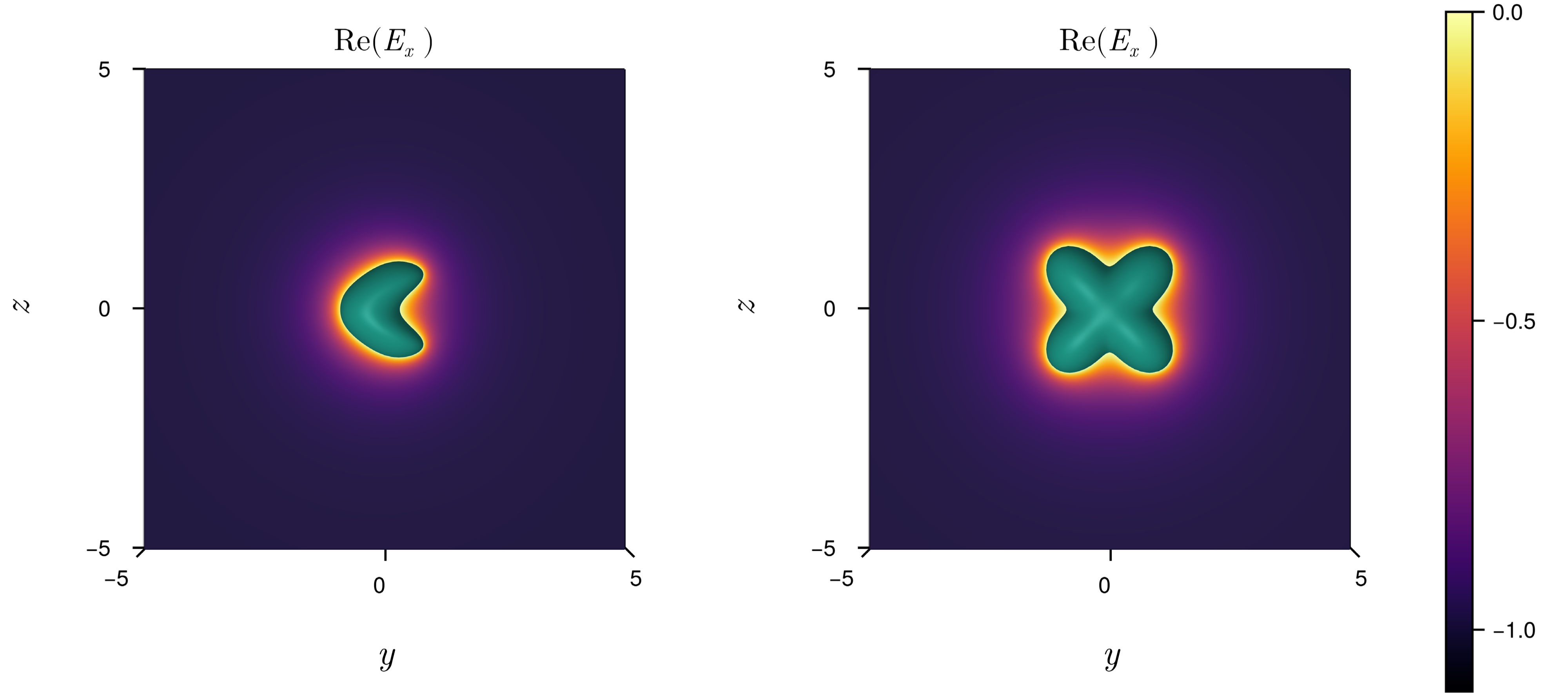}\\[0.2mm]
   \includegraphics[width=0.78\textwidth]{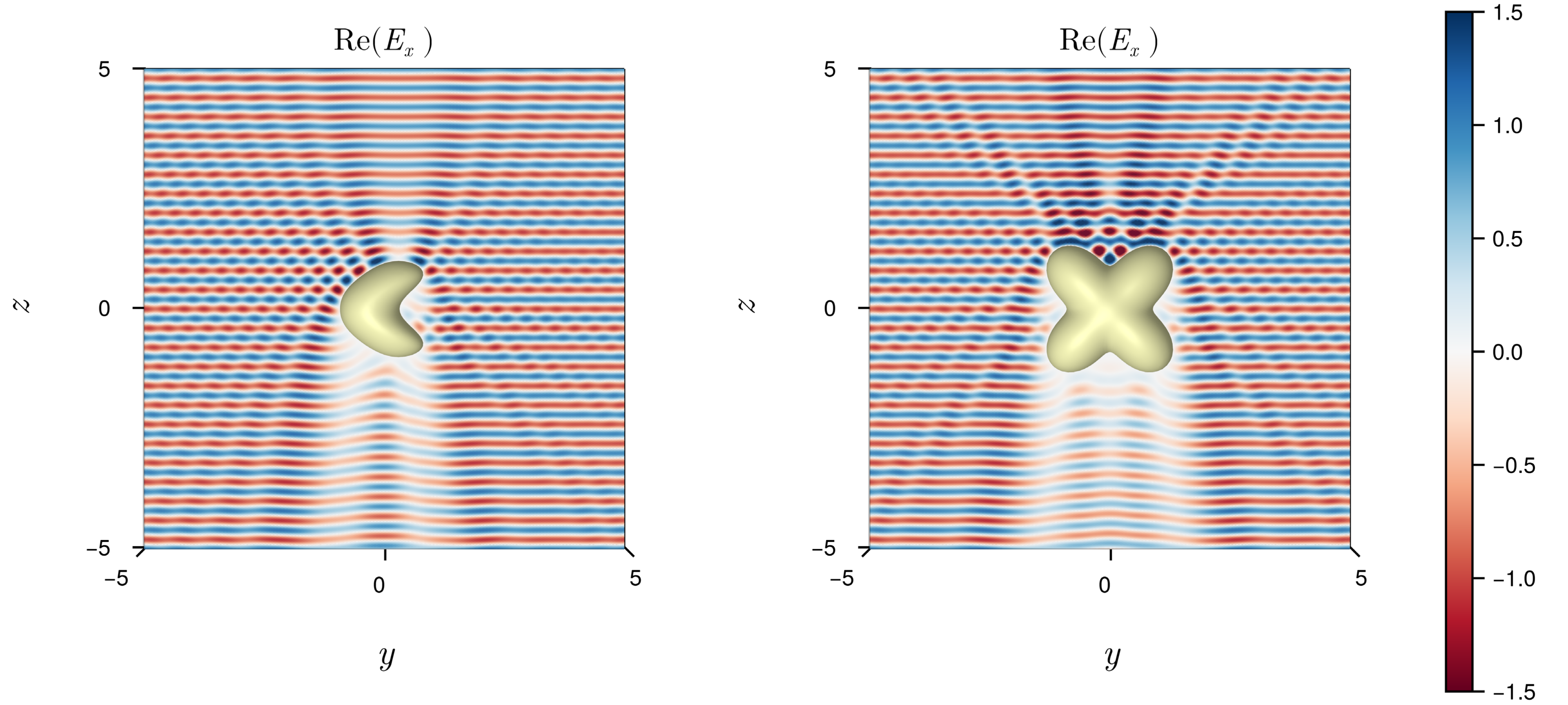}   
   \caption{Electric fields resulting from planewave scattering by the bean- and flower-shaped PEC surfaces~\eqref{eq:smooth_surfaces}. Upper panels: Mean curvature on surface meshes. Middle and lower panels: Real part of the $x$-component of the total electric field on the $yz$-plane at $k=10^{-8}\pi$ (middle) and $k=5\pi$ (bottom). }
   \label{fig:electric_scattering_general_surfaces}
 \end{figure}

\begin{figure}[hbt!]
   \centering
   \includegraphics[width=0.78\textwidth]{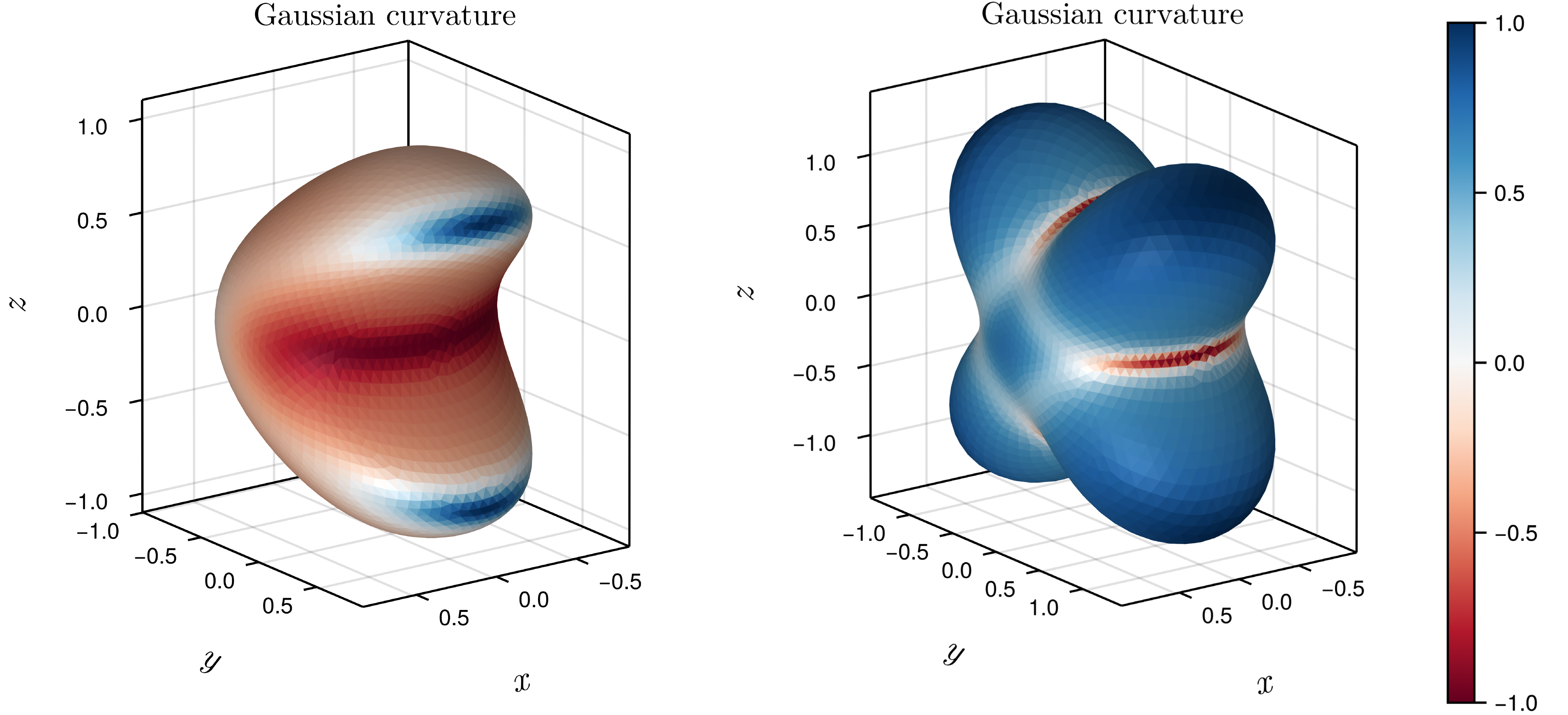}\\[0.2mm]
   \includegraphics[width=0.78\textwidth]{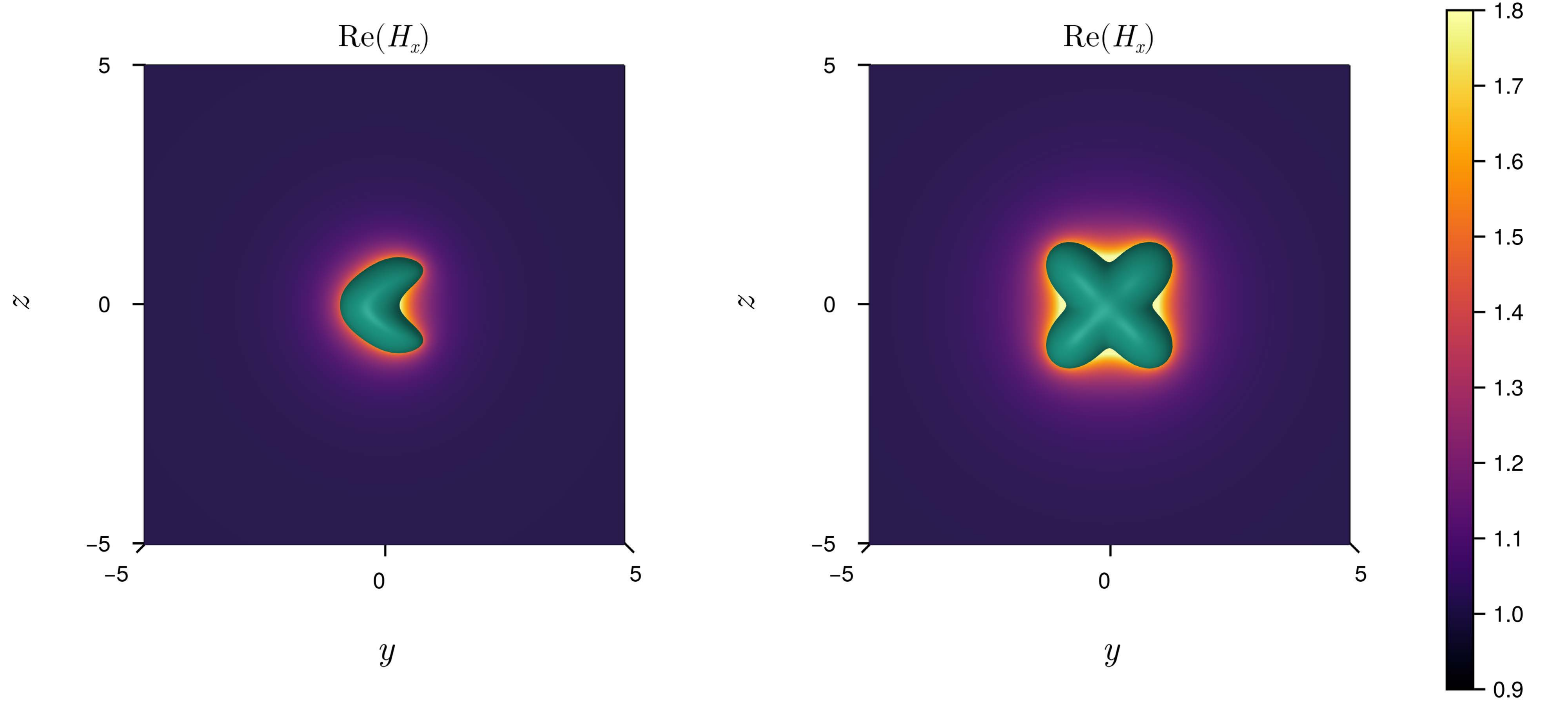}\\[0.2mm]
   \includegraphics[width=0.78\textwidth]{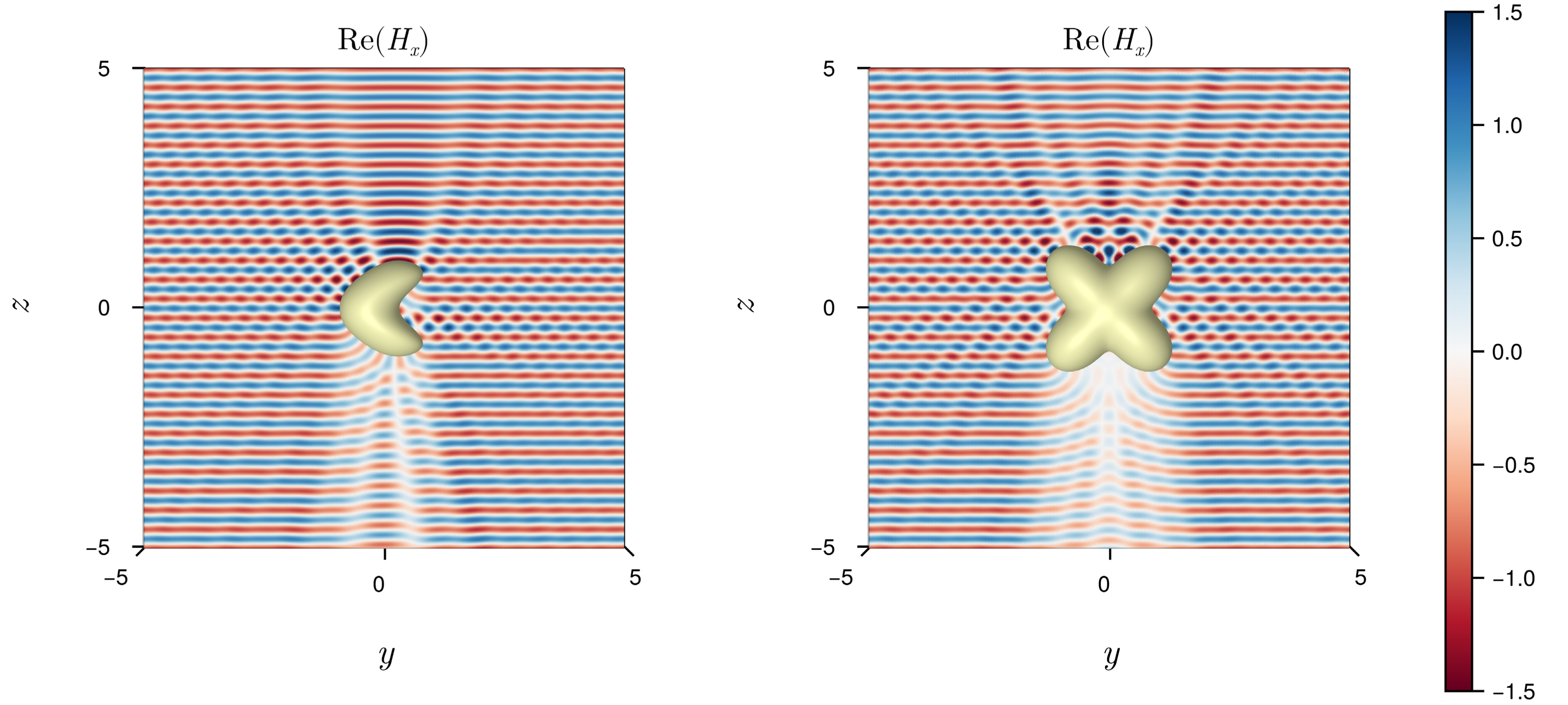}   
   \caption{Magnetic fields resulting from planewave scattering by the bean- and flower-shaped PEC surfaces~\eqref{eq:smooth_surfaces}. Upper panels: Gaussian curvature on surface meshes. Middle and lower panels: Real part of the $x$- component of the total magnetic field on the $yz$-plane at $k=10^{-8}\pi$ (middle) and $k=5\pi$ (bottom).}
   \label{fig:magnetic_scattering_general_surfaces}
 \end{figure}

Subsequently, we solve the corresponding scattering problems under $x$-polarized plane wave illumination~\eqref{eq:planewave}, with $p = [1, 0, 0]$, impinging from above on the PEC structures (i.e. with $d = [0, 0, -1]$). The results of these experiments are shown in Figure~\ref{fig:electric_scattering_general_surfaces}, which shows both the geometry and the computed scattered fields obtained by solving the electric scattering problem~\eqref{eq:elf_equiv} using the R-ECFOIE~\eqref{eq:R-ECFOIE}. The surfaces in the upper panels are colored by their mean curvature $\mathscr H$, required for the solution of the R-ECFOIE. The middle and bottom panels show the $x$-component of the total electric field ($E_x$) on a portion of the $yz$-plane, coinciding with a plane of symmetry of the surfaces. Since the plane normal is tangent to the PEC surfaces, the $x$-component of the total field should vanish on them---a feature accurately captured by our numerical solution.
Similarly, Figure~\ref{fig:magnetic_scattering_general_surfaces} presents the results of the numerical solution of R-MCFOIE~\eqref{eq:R-MCFOIE} used to solve the magnetic scattering problem~\eqref{eq:mgf_equiv}. In this case, the upper panels display the Gaussian curvature of each surface, which is indirectly used for the solution of the R-MCFOIE. The remaining panels report the $x$-component of the total magnetic field ($H_x$) on the same $yz$-plane cross section. It should be noted that the accuracy achieved in resolving the fields near the surfaces is a remarkable capability of the DIM~\cite{faria2021general} implemented in \texttt{Inti.jl}.

\section{Conclusions and future work}
We presented novel boundary integral equations for the PEC electromagnetic scattering problem in three dimensions, formulated entirely in terms of Helmholtz integral operators. This approach bypasses some of the difficulties inherent in classical formulations based on electric and magnetic field integral operators.

However, several avenues for further research remain open. One is the extension of our theoretical results to Sobolev function spaces. A suitable variational formulation of the proposed BIEs, combined with integration by parts, could potentially be used to lower the high surface regularity requirements. We are also actively working on extending this framework to transmission problems, as well as developing a Lippmann–Schwinger-type approach for the numerical solution of inhomogeneous scattering problems, again using only Helmholtz integral operators in a component-wise fashion. These efforts are pursued with the longer term goal of producing fast, high-order, frequency-robust, open-source, integral equation solvers---implemented in the Julia package \texttt{Inti.jl}---based on the polynomial Density Interpolation technique~\cite{anderson2024fast}.

\appendix

\section{Curvature operator and mean curvature}\label{app:curv_tens}
In this section, we derive a matrix representation of the curvature operator $\mathscr{R} = \der \nu$, expressed in terms of the tangent vectors to the surface $\Gamma$ obtained from the derivatives of a local parametrization. Assuming that the surface $\Gamma \subset \mathbb{R}^3$ is $C^{2,\alpha}$-smooth, then for every point $x_0 \in \Gamma$, there exists an open set $U \subset \mathbb{R}^2$, a point $(u_0, v_0) \in U$, and an injective mapping $\mathrm{x} \in C^{2,\alpha}(U, \Gamma)$ such that
$\mathrm{x}: U \to \Gamma, \quad \mathrm{x}(u_0, v_0) = x_0$. 
The map $\mathrm{x}$ is a local parametrization (or chart) of the surface. Its partial derivatives $\mathrm{x}_u$ and $\mathrm{x}_v$ belong to $C^{1,\alpha}(U, \mathbb{R}^3)$ and span the tangent plane to $\Gamma$ at $x_0$. We denote these tangent vectors by ${\rm e}_1 := \mathrm{x}_u$ and ${\rm e}_2 := \mathrm{x}_v$. 
To complete the local frame, we introduce the dual vectors $\mathrm{e}^1$ and $\mathrm{e}^2$ in $\mathbb{R}^3$, defined by the relation $\mathrm{e}_i \cdot \mathrm{e}^j = \delta_{i,j}$ for $i,j \in {1,2}$. Finally, we denote by $\nu_u\in C^{0,\alpha}(U,\R^3)$ and $\nu_v\in C^{0,\alpha}(U,\R^3)$ the partial derivatives of the (well-define) unit normal vector field $\nu=({\rm e_1}\times{\rm e}_2)/|{\rm e_1}\times{\rm e}_2|\in C^{1,\alpha}(U,\R^3)$ with respect to the coordinates $u$ and $v$, respectively.

 As usual~\cite{do2016differential}, we employ the notations
$$
E = \e_1 \cdot \e_1,\quad F = \e_1 \cdot \e_2,\quad \text{and} \quad G = \e_2 \cdot \e_2,
$$
for the coefficients of the first fundamental form: 
$$
g = \begin{bmatrix}E & F \\ F & G\end{bmatrix}, 
$$
with inverse given by
$$ g^{-1} = \frac{1}{EG - F^2}\begin{bmatrix}G & -F \\ -F & E\end{bmatrix}.
$$

We also introduce the coefficients of the second fundamental form:
$$
L = -\e_1 \cdot \nu_u, \quad M = -\e_1 \cdot \nu_v = -\e_2 \cdot \nu_u, \quad N = -\e_2 \cdot \nu_v.
$$

Using these definitions, we express $\nu_u$ and $\nu_v$ in the basis $\{\e^1, \e^2\}$:
$$
\nu_u = -L\e^1 - M\e^2, \quad \nu_v = -M\e^1 - N\e^2,
$$
which can be compactly written as:
$$
[\nu_u \ \nu_v] = -[\e^1 \ \e^2] \begin{bmatrix}L & M \\ M & N\end{bmatrix},
$$
assuming $\e^i$ ($i=1, 2$) are $3 \times 1$ column vectors.

The relationship between the bases of the tangent space can be explicitly written as:
$$
[\e^1 \ \e^2] = [\e_1 \ \e_2]g^{-1},
$$
assuming now that $\e_i$ ($i=1, 2$) are $3 \times 1$ column vectors. Substituting this into the expression above, we obtain:
$$
[\nu_u \ \nu_v] = -[\e_1 \ \e_2]g^{-1}\begin{bmatrix}L & M \\ M & N\end{bmatrix}.
$$

On the other hand, applying the chain rule, we have:
$$
\nu_u = (\der\nor){\rm x}_u = \mathscr R\e_1, \quad \nu_v = (\der\nu){\rm x}_v = \mathscr R\e_2.
$$

For a vector $v \in \R^3$, which can be expressed as $v = (v \cdot \e^1)\e_1 + (v \cdot \e^2)\e_2 + \delta\nor$, it follows that:
$$
\mathscr Rv = \mathscr R\big((v \cdot \e^1)\e_1 + (v \cdot \e^2)\e_2 + \delta\nu\big) = \mathscr R[\e_1 \ \e_2]\begin{bmatrix}v \cdot \e^1 \\ v \cdot \e^2\end{bmatrix},
$$
where we used the fact that $\mathscr R\nu=0$.
Expanding further:
$$
\mathscr Rv = \mathscr R[\e_1 \ \e_2][\e^1 \ \e^2]^\top v = [\nu_u \ \nu_v][\e^1 \  \e^2]^\top v = [\nu_u \ \nu_v]g^{-1}[\e_1 \ \e_2]^\top v.
$$

It then follows from here that the curvature operator can be expressed as:
\begin{equation}\label{eq:curvature_operator}
\mathscr R = -[\e_1 \ \e_2]g^{-1}\begin{bmatrix}L & M \\ M & N\end{bmatrix}g^{-1}[\e_1 \ \e_2]^\top.
\end{equation}
Finally, the mean curvature can be directly calculated using
\begin{equation}\label{eq:mean_curvature}
\mathscr{H} = \frac12{\rm tr}\,\mathscr R.
\end{equation}

Clearly, in view of the assumed $C^{2,\alpha}$-smoothness of the surface $\Gamma$, it follows from~\eqref{eq:curvature_operator} and~\eqref{eq:mean_curvature} that $\mathscr{R} \in C^{0,\alpha}(U, \mathbb{R}^{3 \times 3})$ and  $\mathscr{H} \in C^{0,\alpha}(U)$.

\bibliographystyle{abbrv}
\bibliography{References}

\end{document}